%% file: main.tex
\documentclass[11pt, oneside]{amsart}
\usepackage{tikz-cd}
\usepackage{enumitem}

\usepackage{setspace}

\usepackage[english]{babel}
\usepackage{mathtools}
\usepackage[margin=1in]{geometry}
\usepackage{amssymb}
\usepackage[alphabetic]{amsrefs}
\usepackage{amsmath}
\usepackage{microtype}
\usepackage{graphicx}
\usepackage[colorlinks=true, allcolors=blue]{hyperref}
\usepackage{amsthm}
\usepackage{indentfirst}
\usepackage{import}
\usepackage{xifthen}
\usepackage{pdfpages}
\usepackage{transparent}
\usepackage{changepage}
\hypersetup{
	colorlinks=true,
	linkcolor=blue,
	filecolor=red,      
	urlcolor=red,
	citecolor=teal,
}

\newcommand\blfootnote[1]{%
  \begingroup
  \renewcommand\thefootnote{}\footnote{#1}%
  \addtocounter{footnote}{-1}%
  \endgroup
}

\renewcommand{\th}{^{th}}
\renewcommand{\bar}{\overline}
\newcommand{\e}{\epsilon}
\newcommand{\of}{\circ}
\renewcommand{\to}{\rightarrow}
\newcommand{\halfstrip}{[0,\infty)}
\newcommand{\Z}{\mathbb{Z}}

\newcommand{\N}{\mathbb{N}}
\newcommand{\C}{\mathbb{C}}
\newcommand{\R}{\mathbb{R}}
\newcommand{\del}{\partial}
\newcommand{\D}{\mathbb{D}}
\renewcommand{\P}{\mathbb{P}}
\newcommand{\delbar}{\overline{\partial}}

\newcommand{\Hom}{\text{Hom}}

\newcommand{\iso}{\cong}
\newcommand{\sub}{\subset}
\renewcommand{\sup}{\supset}
\newcommand{\inv}{^{-1}}
\newcommand{\CP}{\C\P}
\newcommand{\Maps}{\text{Maps}}
\newcommand{\interval}{[0,1]}
\newcommand{\wt}{\widetilde}
\renewcommand{\sup}{\text{sup }}
\renewcommand{\bar}{\overline}
\renewcommand{\to}{\rightarrow}
\renewcommand{\P}{\mathbb{P}}
\renewcommand{\sup}{\supset}
\renewcommand{\sup}{\text{sup }}

\newcommand{\coker}{\text{coker }}
\newcommand{\ol}{\overline}
\newcommand{\M}{\mathcal{M}}
\newcommand{\J}{\mathcal{J}}
\newcommand{\ul}{\underline}
\newcommand{\sbl}{\rule{0.4cm}{0.15mm}}
\DeclareMathOperator{\Ind}{Ind}
\DeclareMathOperator{\vir}{vir}
\DeclareMathOperator{\aut}{aut}
\DeclareMathOperator{\im}{im}

\newlist{listC}{enumerate}{2}
\setlist[listC]{label=[\textbf{C}\arabic*],itemindent=*}

\newenvironment{sketch}[1]{\paragraph{Sketch of Proof of #1: \newline} }{\hfill$\square$}
\newenvironment*{prooflemma*}{\paragraph{Proof:}}{}
\newtheorem{Theorem}{Theorem}[section]

\newtheorem{lemma}[Theorem]{Lemma}
\newtheorem{prop}[Theorem]{Proposition}
\newtheorem{cor}[Theorem]{Corollary}

\theoremstyle{definition}
\newtheorem{remark}[Theorem]{Remark}

\newtheorem{defn}[Theorem]{Definition}
\newtheorem{example}[Theorem]{Example}

\numberwithin{equation}{subsection}

\title{Floer Cohomology and Higher Mutations}
\author{Soham Chanda}

\begin{document}
\begin{abstract}
We extend the construction of higher mutation as introduced in \cite{wall} to local higher mutation, which is applicable to a larger class of monotone Lagrangians. In two-dimensional Lagrangians, local higher mutation is the same as performing a Lagrangian anti-surgery in the sense of \cite{haug} followed by a Lagrangian surgery. We prove that up to a change of local systems, the Lagrangian intersection Floer cohomology of a pair of Lagrangians is invariant under local mutation. This result generalizes the wall-crossing formula in \cite{wall}. For two-dimensional Lagrangians, this result agrees with the invariance result in \cite{palmerwoodwardsurg}.
\end{abstract}

\maketitle
\tableofcontents

\setlength{\parskip}{0.5em}
\import{Sections}{intro}

\import{Sections}{setup}

\import{Sections}{SFT}

\import{Sections}{modconstruct}
\import{Sections}{reg}
\import{Sections}{classif}

\import{Sections}{floerhomo}

\bibliographystyle{amsplain}
\bibliography{bib}

{Hill Center,
Rutgers University, 110 Frelinghuysen Road, Piscataway, NJ 08854-8019,
U.S.A.}  

\textbf{sc1929@math.rutgers.edu}
\end{document}

%% file: Sections/intro.tex
\section{Introduction}\label{sec:intro}

A natural question about Lagrangian submanifolds is how their invariants change if one locally modifies a Lagrangian.   We describe a local modification for Lagrangians that contain a \textit{torus segment}  and study the behavior of Lagrangian Floer cohomology under the change.  We call this modification a \textit{local higher mutation}.   Intuitively,  a local higher mutation is a surgery operation, it removes a torus segment from a Lagrangian and replaces with an opposite torus segment.  This modification of a Lagrangian is a local change inspired by Pascaleff-Tonkonog's higher mutation in \cite{wall}. 

Mutations come up naturally in the context of almost toric fibrations of symplectic 4-manifolds.  One obtains a mutation of a smooth toric fiber in an almost toric fibration by varying the location of the singular fiber so that it crosses the smooth fiber, see \cite{leungsym:almosttor}. Auroux \cite{aurTdual} describes a relation between the disc potentials of two families of Hamiltonian non-isotopic Lagrangian tori that are mutations of each other. Seidel \cite{seidelLC}  shows that the Floer homology of a pair of mutated tori in dimension two is non-zero only when the local system over the tori satisfies a certain algebraic relation.  The algebraic relations that occur between the potentials also occur in the theory of \textit{cluster algebras} as algebraic mutations.  Pascaleff-Tonkonog  \cite{wall} generalizes the mutation construction to higher dimensional tori and establishes a wall-crossing formula for their disc-potential.  

Mutations of two-dimensional Lagrangians can be viewed from the perspective of surgery.  Lagrangian surgery, introduced by Polterovich in \cite{pol:surg}, is a version of classical surgery technique.  Mutation of two-dimensional Lagrangians arise from different resolutions of a singular Lagrangian with transverse double intersection by Lagrangian surgery.  Fukaya-Oh-Ohta-Ono \cite{ch10} studies the effect of Lagrangian surgery on the moduli space of discs with boundaries on the Lagrangian.  Woodward-Palmer \cite{palmerwoodwardsurg} proves the invariance of the disc potential and Floer homology under Lagrangian surgery.

 \subsection{Main result}

The main result in this chapter is an invariance of Floer cohomology, up to a change in local system. In spirit, this result is similar to the invariance result in \cite{palmerwoodwardsurg} which proves a similar result for Lagrangian surgery.

We begin with a quick recap of Lagrangian intersection Floer cohomology with local systems. For a more detailed discussion, see Section \ref{sec:floerhomo}. Let $L',K'$ be two transversely intersecting, orientable, spin, compact Lagrangians in a compact monotone symplectic manifold $M'$. Assume that $(L',K')$ forms a monotone tuple.  Let $\mathcal{E}_{L'}, \mathcal{E}_{K'}$ be flat complex line bundles over ${L'}, {K'}$ with holonomies $\rho_{L'}, \rho_{K'} $.  Define the chain complex  $CF({({L'}, \rho_{L'})},  ({K'}, \rho_{K'}))$ as follows, 

\begin{equation*}
    CF({({L'}, \rho_{L'})},  ({K'}, \rho_{K'})) = \bigoplus_{p\in {L'} \cap {K'}} \hom(\mathcal{E}_{{L'}}|_p, \mathcal{E}_{{K'}}|_p ),
\end{equation*}

\noindent  where $\mathcal{E_{K'}}|_p$,$\mathcal{E_{L'}}|_p$ refers to the fibers of the line bundles over the point $p$. The differential $\del$  on $CF({({L'}, \rho_{L'})},  ({K'}, \rho_{K'}))$ is defined by counting rigid holomorphic strips weighted by the holonomies.  Oh \cite{oh:addendum} shows that if all Maslov two discs are simple and the disc potentials $W_{L'}$ and $W_{K'}$ match,  then $\del^2 = 0$.  Later, from structural results of images of pseudoholomorphic discs proved in \cite{kwonoh:sturc}, \cite{Lazz}, it was established that all Maslov two discs with boundary on orientable monotone Lagrangians are simple.  The cohomology of the complex $CF(({L'}, \rho_{L'}), ({K'}, \rho_{K'})$ is denoted as $HF(({L'}, \rho_{L'}), ({K'}, \rho_{K'})),$ which we call the Lagrangian Floer cohomology.   

Now we get back to stating the main result.  We describe relevant notations and assumptions in the following list. 
\begin{listC}
    \item $M$ is a compact, monotone, symplectic manifold of dimension $2n$.
    \item $L$ is a \textit{locally mutable},  oriented,  spin, compact, connected, monotone Lagrangian manifold in $M$. See Definition \ref{defn:locallymutable} for the definition of \textit{locally mutable.}
    \item $K$ is an oriented,  spin, compact, connected,  monotone Lagrangian manifold in $M$. 
    \item The Lagrangian $K$ intersects $L$ transversely but does not intersect a small \textit{mutation neighborhood} $B$ of $L$. See Definition \ref{defn:mutnbd} for  the definition of \textit{mutation neighborhood}.
    \item  $(L, K)$ forms a \textit{monotone tuple}. See Definition \ref{defn:monotup} for  the definition of \textit{monotone tuple} of Lagrangians.
    \item $\rho_L$ and $\rho_K$ are local systems on $L$ and $K$ respectively.
    \item In a mutation neighborhood $B$,  under the identification given by the symplectomorphism $\phi$,  $L$ is equal to the \textit{torus segment} $T_\gamma$. See Definitions \ref{defn:torusseg} for  the definition of \textit{torus segment}.
    \item $x_1, \dots x_{n-1}$ is the image of generators of $H_1(T_\gamma)$ under the inclusion map $i:T_\gamma \hookrightarrow L$.  From Definition \ref{defn:locallymutable}, we know that the induced map $i_* : H_1(T_\gamma) \to H_1(L)$ is injective.
    \item $(x_1, x_2, \dots, x_n,  y_1, \dots y_l)$ is a generating set of $H_1(L)$ such that $x_n \cap i_*[T_\gamma]$ where $[T_\gamma]$ is the generator of $H_{n-1}(T_\gamma)$. A Poincar\'{e} duality argument ensures the existence of such an $x_n$.
\end{listC}

Denote the evaluation of the holonomy $\rho_L$ on $H_1(L)$ as follows, 
\begin{align}
    \rho_L(x_i)=z_i  \label{eq:one}\\
    \rho_L(y_i)=w_i \nonumber
\end{align}

\begin{defn}
  The  \textit{mutated local system} $\rho_{L_\mu}$ on the mutated Lagrangian $L_\mu$ for the pair $(L, \rho_L)$ is given by 
\begin{align} 
    \rho_{L_\mu}(x_i)&=z_i \text{ for } i<n   \label{eq:two}\\
    \rho_{L_\mu}(x_n)&=z_n(1+z_1 + \dots + z_{n-1}) \nonumber \\
    \rho_{L_\mu}(y_i)&=w_i \nonumber.   
\end{align}
\end{defn}

 To obtain regularity of certain moduli spaces,  we use domain-dependent perturbations.  This paradigm of transversality was introduced in \cite{CMtrans} for closed curves and extended to discs with boundary on Lagrangians in \cite{charwoodstabil}, \cite{floerflip}.  One key step in this model is to stabilize domains by adding mark points coming from intersections with a Donaldson divisor.  We assume the existence of such a divisor, which we call \textit{nice divisor},  see Definition \ref{def:nicediv}, which helps us in the regularization of broken maps. One can prove that such divisors do exist by modifying Auroux's argument in Theorem 3.3 of \cite{wall} if the symplectic form is rational and the Lagrangian $L$ satisfies certain geometric conditions.  We now state the main result of the chapter.

\begin{Theorem}[Invariance under mutation]\label{maintheorem}
Let $M,L,K$ be manifolds as described in \textbf{C}[1-9]. Assume that there is a nice divisor $D$ in $M$. Let $L_\mu$ be a local higher mutation of $L$. Let $\rho_L$, $\rho_K$ be choices of local systems on $L, K$ such that $\rho_L$ is given as in equation (\ref{eq:one}). Let $\rho_{L_\mu}$ be the mutated local system as in equation (\ref{eq:two}).  Then,  we have an isomorphism between the disc potentials of the Lagrangians,  $$W_L(\rho_L) \cong W_{L_\mu}(\rho_{L_\mu})$$
Moreover,  if $W_L(\rho_L)=W_K(\rho_K)$, then we have isomorphism of the vector spaces
$$HF((L_\mu, \rho_{L_\mu}), (K, \rho_K)) \cong HF((L, \rho_L), (K, \rho_K))$$
\end{Theorem}

\begin{remark}
The algebraic relation between the local systems in Theorem \ref{maintheorem} matches with the mutation formulae present in the literature.  When $n=2$, this formula matches with the expected change in the local system needed for invariance under Lagrangian surgery as done in \cite{aurTdual}, \cite{seidelLC}, \cite{palmerwoodwardsurg}. For higher dimensions, it agrees with the mutation formula for monotone Lagrangian torus in \cite{wall},  see the relation (\ref{eq:pasctonkmut}). 
\end{remark}

\begin{remark} We need the existence of a divisor to stabilize the disc components in the outside piece of a broken manifold. Usually, one would not require domain dependent perturbation of $J$ when the Lagrangian pair is monotone because Maslov 2 discs are simple from results in \cite{kwonoh:sturc} and \cite{lazz:frame}. To regularize holomorphic buildings, we would require a version of structural results in \cites{lazz:frame,kwonoh:sturc} which works for manifolds with cylindrical ends. To avoid this technicality, we are assuming the existence of a divisor to regularize our holomorphic buildings.
    
\end{remark}

\begin{cor}\label{cor:chekclif}
    Let $L_\mu$ be a higher mutation of a monotone toric fiber $L$ in $(M,\omega)$ as constructed in \cite{wall}. Equip the Lagrangians with local systems $\rho_L, \rho_{L_\mu}$ which are related by the mutation formula \ref{eq:two} and denote the Lagrangian branes as $\wt L, \wt L_\mu$. Then we have, $$HF(\wt{L}_\mu,\wt{L}) \cong HF(\wt L , \wt L).$$
    In particular, by choosing $\rho$ to be a critical point of $W_L$, we know that $HF$ is non-vanishing. Thus, we have that the Clifford and Chekanov-type tori are Hamiltonian non-displaceable from each other.
\end{cor}

\subsection{Higher toric mutation}
Higher-dimensional toric mutation,  as described by Pascaleff-Tonkonog \cite{wall},  is a non-Hamiltonian transformation of the monotone Lagrangian torus fiber in a monotone toric manifold.  We describe the simplest example of such higher mutation here.  Let $M_0=(\C^n,   \omega_{0})$, where $w_0$ is the standard symplectic form $\sum_{i=1}^{n} dx_i\wedge dy_i$ on $\C^n\cong \R^{2n}$.  Consider the fibration given by the multiplication map
\begin{equation}\label{eq:pim}
    \pi_m : M_0 \to \C,  \,  \,  \,  \pi_m(z_1, \dots z_n) = z_1z_2\dots z_n. 
\end{equation}
Given a simple closed curve $\gamma : S^1 \to \C^* $,  we can define an embedded Lagrangian torus
\begin{equation}
\label{eqn:lagtoru}
T_{\gamma} = \{ (z_1, z_2, \dots z_n) \in \C^n | \, \pi_m  (z_1, z_2, \dots z_n) \in \gamma (S^1) ,   |z_1| = |z_2| = \dots |z_n|\}. 
\end{equation}

\begin{figure}[ht]
    \centering

        \def\svgscale{1}
    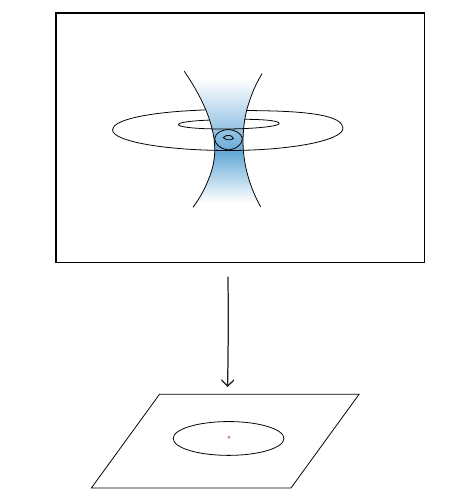

    \caption{The Lagrangian torus $T_\gamma$ }
    \label{fig:tpic}
\end{figure}

The Hamiltonian isotopy class of the Lagrangian $T_\gamma$ depends on geometric properties of the curve $\gamma$.  For a fixed area enclosed by the curve $\gamma$, with respect  to a  symplectic form $d\lambda_n$ (see Definition \ref{defn:lambdan}),  there are at most two Hamiltonian isotopy classes depending on whether the curve $\gamma$ encloses $0$ inside it.  We elaborate on constructing Hamiltonian isotopies of $T_\gamma$ in Subsection \ref{ssub:hamiso}.  When the curve  $\gamma$ encloses 0,  we say that the torus is of Clifford type,  and if it does not,  we say that it is of Chekanov type.  The process of interchanging between Clifford and Chekanov type tori is the simplest occurrence of mutation.  See \cite{chekschlenk:twisttori} for a slightly different construction of the Chekanov tori.  Pascaleff-Tonkonog use the local model $(M_0, T_\gamma)$ to define higher toric mutations in a toric manifold in \cite{wall}. We briefly discuss Pascaleff-Tonkonog's higher toric mutation now.  Assume that $X$ is a monotone toric manifold with Delzant polytope $\Delta$, whose barycenter is 0.  Let $I$ be a line joining a vertex $w$ of $\Delta$  with $0$.  The toric fibration, over a neighborhood of  the line $I$,  is symplectomorphic ( as total spaces,  not fibrations ) to the fibration $\pi_m$ over a neighborhood  of the unit circle $\gamma(\theta)=e^{i\theta} $ attached with the line segment $S$ joining the point $0$ with the point $-1$ such that the monotone torus is mapped to Clifford-type $T_\gamma$.   Choose a curve $\gamma'$ in a neighborhood of the union of curves $\gamma \cup S$, such that the area enclosed by the curve $\gamma'$ is the same as that of $\gamma$ and the Lagrangian torus $T_{\gamma'}$ is of Chekanov-type.  $T_{\gamma'}$ is the mutation of the monotone torus with respect to the vertex $w$. 

We recall the definition of monotonicity.  Let $(M, \omega)$  be a  symplectic manifold. Let $I_\omega$ be the symplectic area functional and  $I_c$ be the first Chern class defined on $\pi_2(M)$. 

\begin{align*}
    &I_\omega : \pi_2(M)\to \R \\
    &I_\omega(u) = \int_u \omega \\
    &I_c : \pi_2(M) \to \Z\\
    &I_c(u) = c_1(u)
\end{align*}

\noindent  Floer \cite{floer:symfixpoint} defines the concept of a monotone symplectic manifold.  A symplectic manifold $(M, \omega)$ is called a \textit{monotone symplectic manifold} if there exists an $\alpha>0$ such that $I_c= \alpha I_\omega$.  This relation between the symplectic area and first Chern class allowed Floer to control sphere bubbling.  Oh  \cite{ohMonotone} extends the concept of monotonicity to Lagrangians  and constructs Lagrangian intersection Floer cohomology.  Let $L$ be a Lagrangian in $(M, \omega)$ and $\mu_L$ be the Maslov index functional for discs with boundary on $L$. In addition, let $I_\omega$ be the symplectic area functional on the homotopy classes of discs with boundary on $L$. 

\begin{align*}
    & I_\omega : \pi_2(M, L) \to \R \\
    &\mu_L : \pi_2(M, L) \to \Z  
\end{align*}

\noindent A Lagrangian $L$ is said to be a \textit{monotone Lagrangian} if there is a $\lambda>0$  such that there is a area-index relation $I_\omega = \lambda \mu_L $.  This area-index relation allows control of disc bubbling.  Initially, the construction of Lagrangian intersection Floer cohomology was based on the assumption that the minimal Maslov numbers of the Lagrangians are at least three.  An addendum, \cite{oh:addendum}, extends the construction to Lagrangians with minimal Maslov index two if the \textit{disc potential} of the Lagrangians matched.  The disc potential of a Lagrangian $L$ is the number of Maslov two discs passing through a generic point on $L$.  

We follow Pascaleff-Tonkonog's exposition of disc potential with a choice of a local system.   Let $\rho:H_1(L)\to \C^*$ be a local system on $L$. We define the disc potential $W_L(\rho)$ by counting Maslov two discs with boundary on $L$.  Let $\operatorname{Tor}$ be the torsion subgroup of $H_1(L)$.  Fix a generating set $$(e_1, \dots, e_m)$$ of $H_1(L, \Z)/\operatorname{Tor} \cong \Z^m $. Let $(x_1, \dots, x_m)$ be the image of $(e_1, \dots,  e_m)$ under the holonomy map $\rho$.   The choice, $$(x_1, \dots, x_m)\in {(\C^*)}^m,$$ specifies the holonomy of a flat $\C^*$ line bundle over $L$ uniquely.  Thus we can view $({\C^*})^m$ as the space of flat line bundles over $L$.  The \textit{disc potential} is a Laurent polynomial in  $x_1, \dots, x_m$, which we interpret as a function on the space of holonomies of flat line bundles  over $L$.  $$W_L:(\C^*)^m \to \C$$ 
$$W_L(x_1, \dots, x_m) = \sum_{u|\mu(u)=2} n_{u}. x^{\del u}$$

\begin{enumerate}
\itemsep0em 
    \item $u$ is a $J$-holomorphic disc with a boundary marking.
    \item $\mu$ is the Maslov index.
    \item $n_u$ is the degree of the evaluation map at the boundary marking $ev:\mathcal{M}_J (L, u) \to L$.
    \item $\del u$ is considered as an element of $\Z^m \cong H_1(L, \Z)/\operatorname{Tor}$.
    \item $x^l$ for $l=(l_1, \dots, l_m)$ denotes $x_1^{l_1}x_2^{l_2}\dots x_m^{l_m}$.

\end{enumerate}

Now we restrict our attention to the disc potential of the Lagrangian torus $T_\gamma$ in a monotone toric manifold. The first homology $H_1(T_\gamma)$ is isomorphic to $\mathbb{Z}^n$. Choose a basis $([e_1], \dots, [e_n])$ of $H_1(T_\gamma)$ such that $e_1, \dots, e_{n-1}$ are cycles on a torus over a single fiber of $\pi_m$ and that $e_n$ is a cycle that descends to the curve $\gamma$.  Pascaleff-Tonkonog \cite[p.~62]{wall} shows that the  disc-potential of $T_{\gamma'}$ can be obtained from that of $T_\gamma$ by replacing $x_i$ in $W_L$ with the following method,  see equation (5.17) in \cite{wall},
\begin{align}
    &x_i \longrightarrow x_i \qquad\qquad\qquad \qquad  \text{ for } i<n \label{eq:pasctonkmut}\\
    &x_n \longrightarrow x_n(1+x_1\dots + x_{n-1} ) \nonumber. 
\end{align}

\subsection{Motivation}

Mutations in symplectic 4-manifolds have played a crucial role in symplectic geometry.  Vianna \cite{vianna:infinitely} proves that there are infinitely many Hamiltonian isotopy classes of monotone Lagrangians in the projective space $\CP^2$ by doing iterated mutations and comparing the disc potentials.  Mutations also have emerged to be important in mirror symmetry,  see \cite{shetruwil:combi}, \cite{hackeating:hommirrLGcal},  \cite{hackkeating:symplecto}.  

Strominger-Yau-Zaslow \cite{syz:conje} conjectures that mirror Calabi-Yau manifolds carry dual torus fibration,  thus one approach to constructing mirror manifolds is to understand torus fibrations.  Higher mutations come up naturally in examples of singular special Lagrangian fibrations.  See \cite{matessi:lag3fib}\cite{gross:spclgeom}, \cite{gross:spcltop} for examples of such singular fibrations in higher dimensions.   

The approach of obtaining the mutation formula presented in this chapter uses the paradigm of neck-stretching in Symplectic Field Theory (SFT) whereas \cite{wall} uses Liouville neighborhoods of the Lagrangian tori.  One advantage of using SFT neck-stretching is that, it allows us to obtain mutation formulas for Lagrangians that are not necessarily diffeomorphic to tori. The construction of local higher mutations enlarges the class of Lagrangian manifolds that can be mutated.  Thus,  understanding the relationship between the Floer theoretic invariants  might allow us to understand finer details about mirror symmetry. For example, we can apply Theorem \ref{maintheorem} to Lagrangian surfaces of higher genus. Although as stated, the theorem requires the ambient symplectic manifold to be compact, the same proof would be applicable to exact Lagrangians in non-compact symplectic manifolds. Thus, we can recover cluster-algebraic relations in the Floer theoretic invariants of exact Lagrangian surfaces under mutation, which are similar to the relations in the microlocal sheaves as found by Shende, Treumann and Williams \cite{shentreuwill:oncombinatorics}.

Another feature of local higher mutations is that they preserve the diffeomorphism type of Lagrangians, while potentially changing the Hamiltonian isotopy class. Thus, the change in Floer theoretic invariants under local mutations can be used to prove the existence of new Hamiltonian isotopy classes of Lagrangians.

\subsection{Outline of the chapter}

The proof of Theorem \ref{maintheorem} is modelled on the neck-stretching argument in Charest-Woodward \cite{floerflip}.  The idea of the proof is to study ``broken holomorphic maps" (also known as holomorphic buildings in the literature), which are the limits of a Symplectic Field Theory (SFT) neck-stretching process.  We perform a neck-stretching along a small contact sphere in a small Darboux neighborhood in the mutation neighborhood.  This setup of the neck-stretching process separates the mutation neighborhood from the rest of the manifold. Finally, we make a local comparison in the stretched mutation neighborhood to obtain the mutation formula.  

 In Section \ref{sec:setup} we explain the construction of local higher mutations. We prove that mutations preserve the monotonicity of a Lagrangian. 
 
 In Section \ref{sec:SFT} we explain the SFT neck-stretching process with boundary conditions and prove SFT compactness.   \cite{bhewz} proves the SFT compactness theorem for limits of closed curves with Morse-Bott Reeb orbit dynamics.  The analysis required for maps with Lagrangian boundary conditions proved so far in the literature has assumed non-degeneracy of Reeb chord dynamics,  see \cite{Abbasbook},  \cite{dylanexpconv}.  These results do not apply to the neck-stretching in our setting since Reeb chords come in families and are always degenerate.  We prove Theorem \ref{theorem:SFTCompactness},  a version of SFT compactness that applies to our setting.  The presentation here is an expanded version of the method outlined in \cite{palmerwoodwardsurg}.  In \cite{Oh-Wang2},  there is a different approach to proving exponential convergence to Reeb orbits, which is an important step for proving the SFT compactness.  It was indicated to the author by Oh that similar methods can be applied to the Morse-Bott type Reeb chord dynamics as well. 
 
 In Section \ref{sec:modcons} we explain the construction of moduli space of broken maps as the intersection of a $\delbar$ section with the zero-section of a Banach bundle.  We prove that the linearization of $\delbar$ in this setting is a Fredholm operator.  We also introduce notation for domain-dependent perturbations of the almost complex structure, which we require to prove regularity of broken maps.  
 
 In Section \ref{sec:reg} we show that the moduli space of rigid (expected dimension zero)  broken maps with boundary on $L, K$ is transversely cut out for a generic choice of domain-dependent perturbation away from the mutation neighborhood.  
 
 In Section \ref{sec:classif} we prove that the only possible rigid broken maps have no ``neck-pieces" and the only possible configuration that can occur in the stretched mutation neighborhood can be classified explicitly.  This classification result is very similar to the classification of rigid objects already found in the literature, e.g.,  Lemma 2.16 in horizontal sections in Lefschetz fibration in \cite{seidelLongEx},  classification of curves with boundary on transversely intersecting Lagrangians and surgeries as in Theorem 60.26 in \cite{ch10} and Proposition 5.13 in \cite{palmerwoodwardsurg}.

In Section \ref{sec:floerhomo} we prove that we can use broken maps to define Lagrangian intersection Floer cohomology in the monotone case and that it is the same as the usual monotone Floer theory as introduced in \cite{ohMonotone}.  Finally  we prove  the main theorem at the end of this section.

%% file: pics/tpic.pdf_tex
\begingroup%
  \makeatletter%
  \providecommand\color[2][]{%
    \errmessage{(Inkscape) Color is used for the text in Inkscape, but the package 'color.sty' is not loaded}%
    \renewcommand\color[2][]{}%
  }%
  \providecommand\transparent[1]{%
    \errmessage{(Inkscape) Transparency is used (non-zero) for the text in Inkscape, but the package 'transparent.sty' is not loaded}%
    \renewcommand\transparent[1]{}%
  }%
  \providecommand\rotatebox[2]{#2}%
  \newcommand*\fsize{\dimexpr\f@size pt\relax}%
  \newcommand*\lineheight[1]{\fontsize{\fsize}{#1\fsize}\selectfont}%
  \ifx\svgwidth\undefined%
    \setlength{\unitlength}{225bp}%
    \ifx\svgscale\undefined%
      \relax%
    \else%
      \setlength{\unitlength}{\unitlength * \real{\svgscale}}%
    \fi%
  \else%
    \setlength{\unitlength}{\svgwidth}%
  \fi%
  \global\let\svgwidth\undefined%
  \global\let\svgscale\undefined%
  \makeatother%
  \begin{picture}(1,1.06666667)%
    \lineheight{1}%
    \setlength\tabcolsep{0pt}%
    \put(0,0){\includegraphics[width=\unitlength,page=1]{tpic.pdf}}%
    \put(0.68689576,0.82997559){\color[rgb]{0,0,0}\makebox(0,0)[lt]{\lineheight{1.25}\smash{\begin{tabular}[t]{l}$T_\gamma$\end{tabular}}}}%
    \put(0.4317199,0.58931433){\makebox(0,0)[lt]{\lineheight{1.25}\smash{\begin{tabular}[t]{l}\small{$\pi_m^{-1}(p)$}\\\\\end{tabular}}}}%
    \put(0.76500002,0.54283758){\makebox(0,0)[lt]{\lineheight{1.25}\smash{\begin{tabular}[t]{l}\Large{$\C^n$}\end{tabular}}}}%
    \put(0.50666665,0.37669884){\makebox(0,0)[lt]{\lineheight{1.25}\smash{\begin{tabular}[t]{l}$\pi_m$\end{tabular}}}}%
    \put(0.55062502,0.05050467){\makebox(0,0)[lt]{\lineheight{1.25}\smash{\begin{tabular}[t]{l}\Large{$\C$}\end{tabular}}}}%
    \put(0.33393784,0.10915633){\makebox(0,0)[lt]{\lineheight{1.25}\smash{\begin{tabular}[t]{l}$\gamma$\end{tabular}}}}%
    \put(0.45812501,0.07253066){\makebox(0,0)[lt]{\lineheight{1.25}\smash{\begin{tabular}[t]{l}\small {$p$}\end{tabular}}}}%
    \put(0,0){\includegraphics[width=\unitlength,page=2]{tpic.pdf}}%
  \end{picture}%
\endgroup%

%% file: Sections/setup.tex
\section{Local Mutation}\label{sec:setup}

In this section, we begin by constructing local higher mutations. Subsequently, we gather some symplectic properties of the objects involved in the construction and prove that local mutations preserve monotonicity of Lagrangians.

\subsection{Construction of Local Higher Mutation}

Similar to defining a Lagrangian torus over a closed loop in $\C^*$ as in (\ref{eqn:lagtoru}) we can define a Lagrangian over a simple open path in $\C^*$.  

\begin{defn}\label{defn:torusseg}
Let $\gamma : [0, 1] \to \C^*$ be a path that does not intersect itself. Assume $\gamma(0) \neq \gamma(1)$.  We define $T_\gamma$,  a \textit{torus segment} over $\gamma$,  as follows, 

\begin{equation}
\label{eq:lagcyl}
    T_{\gamma} =\bigg \{ (z_1, z_2, \dots z_n) \in \C^n \bigg|  \pi_m  (z_1, z_2, \dots z_n) \in \gamma,   |z_1| = |z_2| = \dots |z_n| \bigg \}.
\end{equation} 
\end{defn}
\noindent 
\noindent For a non-intersecting path $\gamma$,  the torus segment $T_\gamma$ is diffeomorphic to the product of a $(n-1)$ dimensional torus and a closed interval of the real line. 

\begin{figure}[ht]
    \centering

        \def\svgscale{1}
    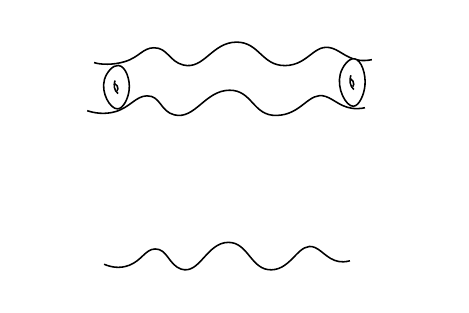

    \caption{Torus segment $T_\gamma$ over a path $\gamma$}
    \label{fig:lagoverpath}
\end{figure}

\begin{defn}[Admissible Path] We call a path $\gamma : (-t, t) \to \C^*$ \textit{admissible} if there exists an $\e >0 $ such that \begin{enumerate}
    \item  $\gamma (x) = x$ for $x\in (-t, -\e) \cup (\e, t)$,
    \item  image of $\gamma|_{(-t+\e, t-\e)}$ lies in the disc $\D(t-\e)$ of radius $t-\e$ ,
    \item $\gamma$ lies completely inside the upper half plane.
\end{enumerate}
 
\end{defn}

\begin{figure}[ht]
    \centering

        \def\svgscale{1}
    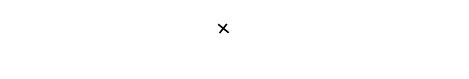

    \caption{Admissible path $\gamma$}
    \label{fig:gamma}
\end{figure}

\begin{figure}[ht]
    \centering

        \def\svgscale{1.5}
    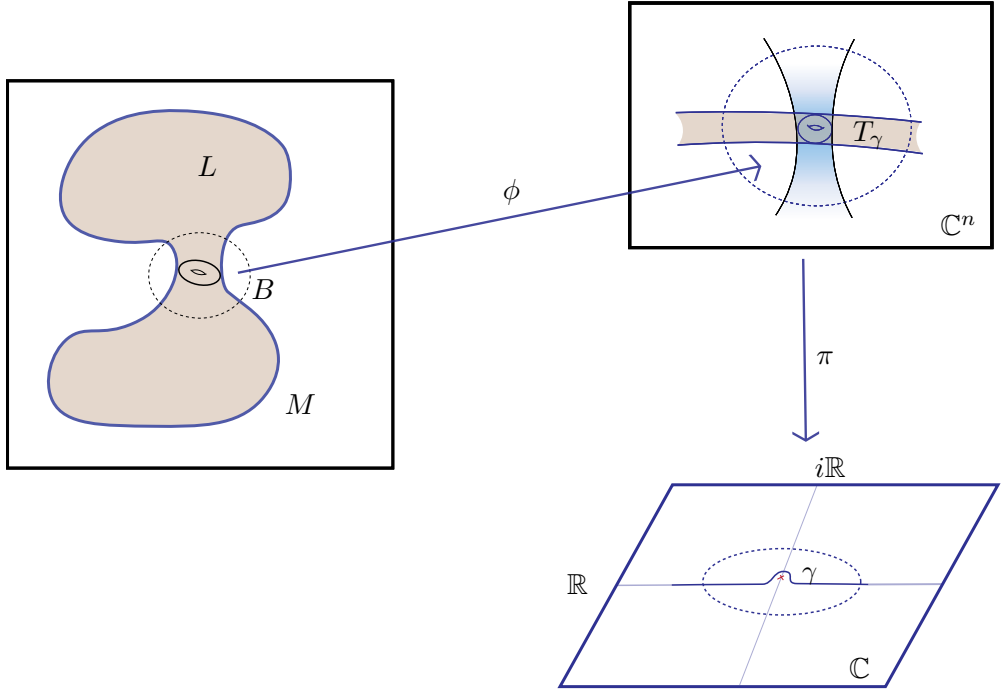

    \caption{Locally mutable $L$ in $M$}
    \label{fig:locpic}
\end{figure}

\begin{defn}[Locally mutable]\label{defn:locallymutable}
We call a Lagrangian $L \hookrightarrow (M, \omega)$ \textit{locally mutable} if the following hold:
 
 \begin{enumerate}
     \item  There  is an open set $B$ in $M$ that is symplectomorphic via $\phi$ to an open ball $B(0)$ around $0 \in \C^n$ such that   $ \phi ( B \cap L) = T_\gamma$ for some admissible path $\gamma$ in $\C$.
     \item $i_*:H_1( B \cap L) \hookrightarrow H_1(L)$ is a split injection, i.e., $0 \rightarrow H_1 (T_\gamma) \to H_1(L)$ extends to split short exact sequence. Here $i_*$ is induced from $i:\phi\inv(T_\gamma) \to L$. 
      \item $L \setminus B$ is connected.
 \end{enumerate}
\end{defn}

\begin{example}
Note that we can use a Hamiltonian isotopy (with respect to $\omega_n$) to map a circle in the complex plane that does not pass through 0 to a curve $\gamma$ that has an arc which is admissible.  Thus, a similar argument as Lemma \ref{lemm:hamiso} shows that the Clifford and Chekanov tori are locally mutable. 

\end{example}

\begin{remark}
For dimension $2$, \cite{haug} proves that if there is an isotropic surgery disc $D$ whose boundary cleanly intersects a Lagrangian $L$, then $L$ satisfies condition (1) of Definition \ref{defn:locallymutable}. 
\end{remark}

\begin{defn}[Mutation neighborhood]\label{defn:mutnbd}
The open set $B$ in Definition \ref{defn:locallymutable} is called the mutation neighborhood. 
\end{defn}

Now we can define local higher mutations.  Intuitively,  Local Higher Mutation is a cut-paste construction where we cut out a torus segment over a path and replace it with another torus segment over ``flipped" curve that goes around $0$ in the opposite direction.  We need to replace the curve while preserving some symplectic area so that eventually we preserve monotonicity of Lagrangians.  

\begin{defn}\label{defn:lambdan}
   We define $\lambda_n$ as the  1-form  on $\C^*$ given by  $\lambda_n = \frac{xdy-ydx}{2(x^2+y^2)^{\frac{n-1}{n}}}$.  The form $\lambda_n$ is obtained by pulling back the 1-form $\lambda_0 = \sum_i x_idy_i - y_idx_i,$ on $\C^n$, under the diagonal $n^{th}$-root map, $z\mapsto (z^{1/n},  \dots,  z^{1/n})$.   
\end{defn}

\begin{defn}[Local Higher Mutation]\label{defn:mutation}
 Given a locally mutable Lagrangian $L$,  the local higher mutation $L_\mu$ of $L$ is the Lagrangian obtained by removing $\phi\inv ( T_c)$  for an arc $c$ of $\gamma$ i.e.  $c:(-r, r) \to \C^*$ and $c(x)=\gamma(x)$ and gluing back $\phi \inv (T_{c'})$ for a curve $c'$  such that,
 \begin{enumerate}
    
     \item there is an $\e>0$ such that $c(x)= c'(x)$ for $x\in (-r, -r+\e) \cup (r-\e, r)$,
     \item the glued curve $c\cup c'$ is homotopic to the unit circle around 0 ,
     \item the resulting curve $\gamma' = (\gamma \setminus c )\cup c'$ is smooth,
     \item $\int_c \lambda_n = \int_{c'} \lambda_n$.
 \end{enumerate}
\end{defn}

\begin{defn}[Mutated admissible path]
    Given an admissible path $\gamma$, we call the curve $\gamma'$ as constructed in Definition \ref{defn:mutation} the mutated admissible path of $\gamma$.
\end{defn}

\begin{remark}
For dimension $2$, Lagrangian mutation is the same as performing an anti-surgery as defined in \cite{haug} followed by a Lagrangian surgery. 
\end{remark}

\begin{figure}[ht]
    \centering

        \def\svgscale{0.4}
    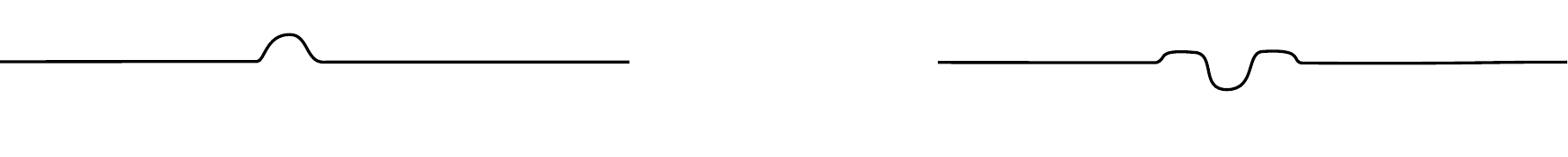

    \caption{Mutation of paths}
    \label{fig:gamm}
\end{figure}

We will show that the monodromy of the fiber torus around the origin is trivial in the mapping class group.  In Definition \ref{defn:mutation}, condition (1) implies replacing $c$ with $c'$ does not change the topology of the Lagrangian $L$.  We will also show that condition (3) implies that the symplectic areas of the discs remain the same after local higher mutation and Lemma \ref{lemm:hamiso} implies all such choices of paths $c'$ produces Hamiltonian isotopic $L_\mu$.

The construction of higher mutation is related to an example of special Lagrangian fibration that appears in \cite{joy}.  When $\gamma$ is the real axis,  $T_\gamma$ is the Harvey-Lawson cone when $k\text{ is odd}$ (see example 3.5 in \cite{joy}).   Consider the Harvey-Lawson Lagrangian fibration in dimension 3 $$\rho : \C^3 \to \R^3 $$    $$\rho(z_1, z_2, z_3) = (\text{Im} (z_1z_2z_3), |z_1^2|-|z_2^2|, |z_1^2|-|z_3^2|). $$
The singular Lagrangian $\rho \inv (0)$ is the Harvey-Lawson cone.  The Clifford and Chekanov type tori can be Hamiltonian isotoped so that they match the Harvey-Lawson fibers locally near 0.  In this setting,  higher mutation can be thought of as replacing $\rho\inv (\epsilon, 0, 0)$ with $\rho\inv(-\epsilon, 0, 0)$ locally and gluing it back appropriately.

\subsection{Monodromy of the fiber Lagrangian torus}

We will use symplectic parallel transport of the fibration $\pi_m$ (see Equation (\ref{eq:pim}) for definition of $\pi_m$) to obtain various properties of the Lagrangians we have been considering so far. One way to verify the manifold $T_\gamma$ is a Lagrangian, is by showing it is obtained by a symplectic parallel transport of a Lagrangian in a fiber. Notice that $\pi_m$ is a fibration over $\C^*$. The standard symplectic form $\omega_0$ in $\C^n$ induces a symplectic horizontal distribution, which allows us to parallel transport between fibers such that the parallel transport maps are symplectomorphisms of the fiber.

We describe how a symplectic fibration provides us a splitting of the tangent bundle of the total space of the fibration. The tangent bundle of $\C^n$ splits into vertical and horizontal components under the fibration $\pi_m$. Away from the locus of critical points of $\pi_m$, the vertical bundle $V$ of the fibration is sub bundle obtained as the kernel bundle $\ker d\pi_m$ i.e. over a point $z\in \C^n$, the fiber $V_z$ of the vertical bundle is given as $V|_z = \ker d\pi_m(z)$. The horizontal bundle $H$ is defined as the symplectic complement $V^\omega$ of the vertical bundle $V$.  We can explicitly compute $\ker d\pi_m$ at a point $z =(z_1,\dots z_n) \in \C ^n$ to obtain $V_z$. For a point $z$ such that $\pi_m(z)\neq 0$, we have that
\begin{align*}
V_z &= \text{Span}_\C(X_1(z),\dots, X_{n-1}(z)),    \\
\text{ where }& X_i(z) = (0, \dots 0, \quad \underbrace{ z_i}_{i\th \text{ spot. }} ,\quad 0, \dots ,-z_n ).
\end{align*}
A direct computation can be done to obtain symplectic complement to the vertical bundle over any $z$ such that $\pi_m(z)\neq 0$. We have the following relation.

\begin{equation}\label{eq:vertbundle}
    H_z = \text{Span}_\C \bigg \{ \big (\frac{z_1}{|z_1|^2}, \dots \frac{z_n}{|z_n|^2}\big ) \bigg\}
\end{equation}

\begin{defn}
We define $T_{\ul z }$, the Lagrangian torus in the fiber over $\ul z\in \C^*$ to be the torus given by the relations $|z_1|=\dots=|z_n|$ in the fiber $\pi_m\inv (\ul z)$.
\end{defn}

\begin{lemma}
The manifold $T_\gamma$ is a Lagrangian submanifold of $\C^n$ equipped with the standard symplectic structure and $T_\gamma$ is diffeomorphic to a torus.  
\end{lemma}

\begin{proof}
Note that when $|z_1|=\dots=|z_n|$ for some $z=(z_1,\dots,z_n) \in \C^n$, from equation (\ref{eq:vertbundle}) we have that the horizontal distribution at $z$ is just $\C .z$. The torus $T_{\ul z}$ is a Lagrangian submanifold in $\pi_m\inv (\ul z) $ by a direct computation and the manifold $T_\gamma$ is the union of tori $T_{\gamma(t)}$ for all $t$, i.e. $ \cup_t T_{\gamma(t)}$. A simple check shows that the map $z\mapsto (z^{\frac{1}{n}} z_1, \dots ,z^{\frac{1}{n}} z_n )$ is a horizontal section for $(z_1,\dots,z_n) \in \pi_m \inv (1)$ such that $|z_1|=\dots=|z_n|$. This proves that the symplectic parallel transport from the fiber over the point $a$ to the fiber over the point $ b$ maps the Lagrangian torus $T_a$ to the torus $T_b$. Thus $T_\gamma$ is a Lagrangian.

From the proof above one gets that the monodromy obtained from a circle around origin $T_1 \to T_1$ is multiplication by an $n^{th}$ root of unity $\xi$ on all factors, thus $T_\gamma$ is a trivial bundle over $\gamma$ hence a torus.
\end{proof}

\begin{remark} \label{rem:torsseglag}
The proof above also shows that torus segments are Lagrangian submanifolds.
\end{remark}

\begin{lemma}\label{lemm:Texactlag}
The torus segment $T_\gamma$ is an exact Lagrangian for any smooth, non-intersecting  path $\gamma : \interval \to \C^*$ .
\end{lemma}

\begin{proof}

The standard symplectic form $\omega_0$ has a standard primitive $\lambda_0 = \sum_{i=1}^{n}\frac{1}{2}(x_idy_i-y_idx_i)$. The restriction of $\lambda_0$ to $T_\gamma$ is closed because $T_\gamma $ is a Lagrangian. Now, to show that it is exact, it is enough to show that  line integrals of $\lambda_0$ over the generators of $H_1(T_\gamma)$ is $0$. Clearly the torus segment $T_\gamma$ retracts onto $T_p$ for any $p$ on $\gamma$ thus $H_1(T_\gamma) \cong \mathbb{Z}^{n-1}$. The generators are given by the \textit{standard }$n-1$ \textit{loops} $S^1 \to T_p$, $\theta \mapsto (p_1,\dots ,e^{i\theta}p_i,\dots,e^{-i\theta}p_n)$ where $(p_1,\dots,p_n) \in T_p$ .  Using Stokes' theorem and the fact that $|p_i| = |p_n|$, we see that the line integral $\lambda_0$ over the $n-1$ generators is $0$. Thus $T_\gamma$ is exact.

\end{proof}

\begin{lemma}\label{rem:texactlag}
For any path $c$ to the Lagrangian torus in the fiber over $p$, i.e. $c:[0,1] \to T_p$, we have $\int_c \lambda_0 = 0$.
\end{lemma}

\begin{proof}
From \ref{lemm:Texactlag}, the one form $\lambda_0$ restricted to the Lagrangian,  $T_\gamma$ is exact. Thus, we have that the line integral $\int_c \lambda_0 $ depends solely on the end points $c(0), c(1)$. Now choose a path from $c(0)$ to $c(1)$ that is a composition of arcs from the standard $n-1$ loops defined above. Since $\lambda_0$ is rotationally invariant, a direct computation shows that the integral of $\lambda_0$ over any arc of the standard $n-1$ loops is 0. Thus, we have that $\int_c \lambda_0 = 0 $.

\end{proof}

 \subsection{Hamiltonian isotopy} \label{ssub:hamiso}

In this subsection, we study the Hamiltonian isotopy classes of the family of Lagrangians $T_\gamma$ based on properties of the path $\gamma$. We follow the approach of Corollary 2.5 in \cite{lekiliRational} to provide a condition on paths $\gamma$ and $\gamma'$ that is sufficient to make the torus segments $T_\gamma$ and $T_\gamma'$ Hamiltonian isotopic. We study the case of smoothly isotoping the path $\gamma_0 : [0,1] \to \C^*$ to another path $\gamma_1$ where $\gamma_s$ is fixed in a small neighborhood of $0$ and $1$, i.e., there is an $\e$ such that $\gamma_s (t) = \gamma_0 (t)$ for $t\in [0,\e) \cup (1-\e,1]$. Call $\gamma(0)=a$ , $\gamma(1)=b$.
\begin{lemma}\label{lemm:hamiso}
Let $\gamma_s$ be an isotopy of paths that is constant on endpoints as described above, such that $\int_{\gamma_s}\lambda_n$ is constant, $T_{\gamma_0} $ is Hamiltonian isotopic to $T_{\gamma_1}$
\end{lemma}
 
\begin{proof}
The isotopy $\ol \gamma:\interval \times \interval \to \C^*$, of paths $\gamma_s:=\ol \gamma (s,\_)$, from the curve $\gamma_0$ to $\gamma_1$, lifts to an isotopy of embedded manifolds (with boundary) from the torus segment $T_{\gamma_0}$ to $T_{\gamma_1}$ via torus segments. Since we know that each torus segment is a Lagrangian manifold, we have a Lagrangian isotopy between $T_{\gamma_0}$ to $T_{\gamma_1}$. An equivalent condition to Lagrangian isotopy being a result of ambient Hamiltonian isotopy is that the Lagrangian isotopy is exact, see Section 6.1 \cite{polGeom}. To deal with an isotopy of  manifolds with boundary,  we use the relative de Rham model to define exactness. In our setting, being exact turns out to be the same as having the integral $\int_{c_s} \lambda_0$ for an isotopy (induced from the Lagrangian isotopy) of curves $c_s \in H_1(T_{\gamma_s},T_a \cup T_b)$ invariant of isotopy in the $s$-direction. Here $c$ is a homotopy of paths lying on the isotopy of the torus segments over $\ol \gamma$ with boundary on the torii $T_a,T_b$ i.e. $c_s:=c(s,\_)$ is a path (or loops) in $T_{\gamma_{s}}$ such that the boundary (if any) $\del c_s$ lies on the boundary $T_a \cup T_b$ of the torus segment. We have observed that the torus segment $T_\gamma$ is an exact Lagrangian in Lemma \ref{lemm:Texactlag} and that the integrals of $\lambda_0$ are 0 on paths that lie in a single fiber  $T_z$ in Lemma \ref{rem:texactlag}. Thus,  it is enough to check whether $\int_{c_s} \lambda_0$ is invariant of $s$ for a family of curves $c_s$ that descend to the path $\gamma_s$. A candidate for such a lift is $c_s(t) = (\gamma_s(t)^{1/n}, \dots, \gamma_s(t)^{1/n})$. Clearly $\int_{c_s}\lambda_0 = \int_{\gamma_s} \lambda_n$ from the definition of $\lambda_n$. From the hypothesis of the lemma, we see that $\int_{c_s}\lambda_0$ is invariant of $s$. Thus, we have proved the Lemma.

\end{proof}

 \begin{cor}
 If  $\gamma_0$ and $\gamma_1$  are two admissible paths that are equal near the ends and $\int_{\gamma_0} \lambda_n = \int_{\gamma_1} \lambda_n $ then $T_{\gamma_0}$ is Hamiltonian isotopic to $T_{\gamma_1}$.
 \end{cor}
 \begin{proof}
 Lemma \ref{lemm:hamiso}  reduces the problem of finding a Hamiltonian isotopy of the Lagrangian $T_\gamma$ for admissible paths $\gamma$ to finding an isotopy of $\gamma$ that keeps the line integral with respect to $\lambda_n$ constant. Assume that $\gamma_0$ and $\gamma_1$ are two admissible paths that are equal near the ends, are isotopic to each other and $\int_{\gamma_0} \lambda_n = \int_{\gamma_1} \lambda_n $. We can `close up' $\gamma_0$ and $\gamma_1$ to an embedded circle with the same area by attaching a path  $\gamma_u$ that lies in the upper half plane. Lemma 2.3 of \cite{lekiliRational} show that $\gamma_0 \cup \gamma_u$ is Hamiltonian isotopic with respect to the symplectic form  $\omega_n=d\lambda_n$ to $\gamma_1 \cup \gamma_u$ where the isotopy is constant on $\gamma_u$. Thus we get an isotopy $\gamma_s$ of $\gamma_0$ to $\gamma_1$ such that  $\int_{\gamma_s}\lambda_n$ is constant. Thus, using Lemma \ref{lemm:hamiso} we finish the proof of our corollary.
 \end{proof}

\noindent We define a cylindrical version of the Lagrangian $T_\gamma$ when the path $\gamma$ is chosen to be cylindrical. The construction depends on choosing a path on $\C$ that lies on the real line away from a compact set.

\begin{defn}[Cylindrical path]
 A cylindrical path is defined as a map $\gamma: \R \to \C^*$ such that  $\gamma(r) \in \R$  when $|r|>M$ for some $M>0$ and the following have the following limits 
 $$\lim_{r \to \pm \infty} \gamma(r) = \pm \infty.$$
 \end{defn}
 
 \begin{defn}[Asymptotically cylindrical path]
 A path $\gamma : \R \to \C^*$ is defined to asymptotically cylindrical if there is an $\e\neq 0$  such that  $\gamma(r) \in \R+i\e$  when $|r|>M$ for some $M>0$ and the following have the following limits 
 $$\lim_{r \to \pm \infty} \gamma(r) = \pm \infty.$$
 \end{defn}
 
 \begin{defn}[Cylindrical extension]
 Let $\gamma:(-t,t)$ be an admissible path. From the definition of admissible path, we know there is an $\e>0$ such that $\gamma(x)=x$ for $x\in(-t,t)\setminus (-\e,\e)$. The cylindrical extension $\gamma^e$ is defined as 
 \begin{align*}
     \gamma^e(x) = x \text{ for x } \in \R\setminus(-\e,\e) \\
     \gamma^e(x) = \gamma(x) \text{ otherwise.}
 \end{align*}
 \end{defn}

We conclude this subsection with the definition of a cylindrical version of  torus segments, that shows up in the later sections as a result of neck stretching. We define a Lagrangian $T_\gamma$ over a cylindrical path $\gamma$ similarly as we define a torus-segment in Definition \ref{defn:torusseg}. We can similarly define $T_\gamma$ for an asymptotically cylindrical path $\gamma$.

Since the path $\gamma$ lies on the real axis near its ends, the Lagrangian $T_\gamma$ has cylindrical ends in the sense that for large $R>0$, $T_\gamma \setminus B_R(0) =  [M,\infty) \times \Lambda $ for some $M>0$ in cylindrical coordinates $\R \times S^{2n-1} $ of $\C^n\setminus \{0\}$  where $\Lambda$ is the disjoint union of two Legendrian tori $T_+,T_-$ in the unit sphere $S^{2n-1}$ given by the equations 
\begin{equation} \label{def:tplusminus}
    T_\pm = \bigg \{ (z_1,\dots,z_n) \in S^{2n-1} \bigg | |z_1| = |z_2| = ..=|z_n| , z_1.z_2.\dots.z_n \in  \R^\pm  \bigg \} .
\end{equation}

\subsection{Mutation preserves monotonicity}
In this subsection, we will describe which Lagrangians can be locally mutated and explain the construction of local mutations and show that monotonicity of Lagrangians is preserved under the operation of mutation.

\begin{lemma}\label{mutmono}

If L is a locally mutable monotone Lagrangian, then its mutation $L_\mu$ is also monotone.
\end{lemma}

\begin{proof}
To fix notation, we use $B$ to denote the mutation neighborhood of $L$, $\gamma$ to denote the curve that determines $T_\gamma$ as in  Definition \ref{defn:locallymutable} and $\gamma'$ to denote the path for the mutation. For any disc $u^\mu $ with boundary on $L_\mu$, we will construct a disc with boundary on $L$, $u$, with the same area and Maslov number as $u^\mu$, which would prove that $L_\mu$ is monotone. Since $L=L_\mu$ outside the mutation neighborhood, we need to modify $u^\mu$ only where it enters the mutation neighborhood. Let $\del u^\mu$ denote the restriction of the disc $u^\mu$ to the boundary circle  and assume that on the domain interval $[p,q]$, $\del u^\mu$ is a map to the mutation neighborhood i.e. $\del u^\mu : [p,q] \to B$. We can choose points $p,q$ on the circle so that the $[p,q]$ is a maximal subset of the boundary circle with this property. There can be two cases,
\begin{enumerate}

    \item $\del u^\mu (p)$ and $\del u^\mu(q)$ lie on the same end of $T_\gamma$, i.e.\ $\pi_m(\del u^\mu(p)) = \pi_m(\del u^\mu(q))=:z$,
    \item $\del u^\mu (p)$ and $\del u^\mu(q)$ lie on the opposite ends of $T_\gamma$, i.e.\ $\pi_m(\del u^\mu(p)) \neq \pi_m(\del u^\mu(q))$.
\end{enumerate}
 
 For case (1), we can homotope the path $\del u^\mu$ over $[p,q]$ through paths in $L_\mu$  while fixing the end points to a path that lies completely on a single torus fiber $T_z \sub T_\gamma$. We can use this homotopy to extend $u^\mu$ so that the boundary on $[p,q]$ lies on $T_z$, call the extended map $u_{pre}$. Note that this doesn't change the area since the homotopy of paths lied in a Lagrangian, also, it doesn't change the Maslov number of the boundary. Since $\del u_{pre}$ maps the boundary $[p,q]$ to $T_z$, using the fact that $z\in \gamma$ we get that $\del u_{pre} ([p,q]) \sub L$ since $L = T_\gamma$ in the mutation neighborhood.

 For case (2), as a preliminary step, we homotope $\del u^\mu$ to a path on $T_{\gamma'}$ such that on $(p+e,q-e)$ for small $e$, $\del u^\mu$ is a horizontal section of $\gamma'$. We then attach horizontal sections over the discs bounding $\gamma \cup \gamma'$ and not containing 0 to extend $u^\mu$ to $u_{prepre}$ so that the boundary lies on $T_\gamma$. We need to deal with the case of discs containing 0 carefully, note that since the horizontal section is given by taking the $n\th$ root, it can't be defined in a neighborhood of 0, but only on a slit neighborhood, i.e., after removing a half-line starting at 0. We can remove a ray $\R^+$ from the neighborhood and then attach two copies of a line to the neighborhood to get a domain topologically looking like a radial sector removed from a disc. Call the two newly attached lines as $l_1, l_2$. Refer to figure (\ref{fig:monotonicity_sections}). We can continuously define an $n\th$ root map over this domain (colored yellow in the figure), and then attach a sector ( colored green in the figure) of 0 symplectic area that lies on $T_{\R^+}$ and is a homotopy between the two branches of the horizontal section over $l_1,l_2$. Use this homotopy to extend $u_{prepre}$ to $u_{pre}$. The boundary of $u_{pre}$ on $[p,q]$ lies on $L$. Note that, as before, the symplectic area of $u_{pre}$ is the same as that of $u^\mu$ from the fact that the area contributions from the horizontal sections cancel out from condition (4) in Definition \ref{defn:mutation}. Also, the Maslov number doesn't change since the generators of $H_1(T_z)$ have Maslov number 0.
\begin{figure}[ht]
    \centering

        \def\svgscale{1.5}
    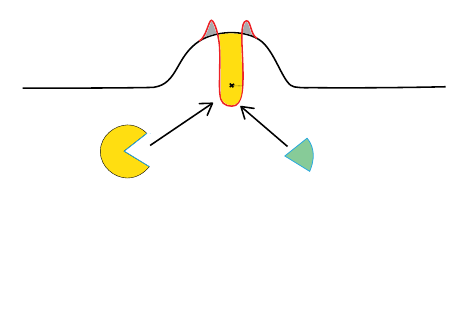

    \caption{Disc used in creating homotopy from $u_{prepre}$ to $u_{pre}$ }
    \label{fig:monotonicity_sections}
\end{figure}

After extending $u^\mu$ as shown in the last two paragraphs for every maximal domain $[p,q]$ on the boundary where $\del u^\mu$ is a map to the mutation neighborhood, we can get a map $u$ that has the following properties, 
\begin{itemize}
    \item the boundary of the map lies on $L$,
    \item the Maslov number is $I_m(u^\mu)$,
    \item  the symplectic area is $\omega(u^\mu)$.
\end{itemize}
Since $L$ is monotone, we get that $I_m(u^\mu)  = \lambda \omega(u^\mu)$. The constant $\lambda$ doesn't depend on $u^\mu$ we began with. Thus, by repeating this construction for each $u^\mu \in \pi_2(M,L_\mu)$ we see that $L_\mu$ is monotone.

\end{proof}

\begin{defn}\label{defn:monotup}(Monotone tuple of Lagrangians, Definition 4.1.2 in \cite{quilt})
A pair of Lagrangians $(L_1,L_2)$ is called a monotone pair if there is a  there exists a $\tau$ for which we have the following an action-index relation for maps with boundary on $L_i$. Let $u : \Sigma\to M$ where $\del \Sigma = \{C_i\}$ and $u(C_i) \sub L_i$, then we have $$2\int u^*\omega = \tau \mu(u).$$
\end{defn}

\noindent We can use the method of Lemma \ref{mutmono} to prove that monotonicity of a Lagrangian pair is preserved when one of the Lagrangians is mutated.

\begin{lemma}\label{lemm:mutmonopair}
If $(L,K)$ is a monotone pair of Lagrangians and $L$ is locally mutable, then the pair $(L_\mu,K)$ is also a monotone pair.
\end{lemma}

%% file: pics/lagoverpath.pdf_tex
\begingroup%
  \makeatletter%
  \providecommand\color[2][]{%
    \errmessage{(Inkscape) Color is used for the text in Inkscape, but the package 'color.sty' is not loaded}%
    \renewcommand\color[2][]{}%
  }%
  \providecommand\transparent[1]{%
    \errmessage{(Inkscape) Transparency is used (non-zero) for the text in Inkscape, but the package 'transparent.sty' is not loaded}%
    \renewcommand\transparent[1]{}%
  }%
  \providecommand\rotatebox[2]{#2}%
  \newcommand*\fsize{\dimexpr\f@size pt\relax}%
  \newcommand*\lineheight[1]{\fontsize{\fsize}{#1\fsize}\selectfont}%
  \ifx\svgwidth\undefined%
    \setlength{\unitlength}{225bp}%
    \ifx\svgscale\undefined%
      \relax%
    \else%
      \setlength{\unitlength}{\unitlength * \real{\svgscale}}%
    \fi%
  \else%
    \setlength{\unitlength}{\svgwidth}%
  \fi%
  \global\let\svgwidth\undefined%
  \global\let\svgscale\undefined%
  \makeatother%
  \begin{picture}(1,0.66666667)%
    \lineheight{1}%
    \setlength\tabcolsep{0pt}%
    \put(0,0){\includegraphics[width=\unitlength,page=1]{lagoverpath.pdf}}%
    \put(0.194814,0.04518633){\makebox(0,0)[lt]{\lineheight{1.25}\smash{\begin{tabular}[t]{l}$\gamma$\end{tabular}}}}%
    \put(0,0){\includegraphics[width=\unitlength,page=2]{lagoverpath.pdf}}%
    \put(0.39980794,0.36693667){\makebox(0,0)[lt]{\lineheight{1.25}\smash{\begin{tabular}[t]{l}$T_\gamma$\end{tabular}}}}%
    \put(0,0){\includegraphics[width=\unitlength,page=3]{lagoverpath.pdf}}%
    \put(0.81722333,0.29277667){\makebox(0,0)[lt]{\lineheight{1.25}\smash{\begin{tabular}[t]{l}$\C^n$\end{tabular}}}}%
    \put(0.73166667,0.03222333){\makebox(0,0)[lt]{\lineheight{1.25}\smash{\begin{tabular}[t]{l}$\C$\end{tabular}}}}%
    \put(0.45416,0.06141266){\makebox(0,0)[lt]{\lineheight{1.25}\smash{\begin{tabular}[t]{l}$0$\end{tabular}}}}%
    \put(0,0){\includegraphics[width=\unitlength,page=4]{lagoverpath.pdf}}%
    \put(0.52276952,0.23149012){\makebox(0,0)[lt]{\lineheight{1.25}\smash{\begin{tabular}[t]{l}$\pi_m$\end{tabular}}}}%
  \end{picture}%
\endgroup%

%% file: pics/gamma.pdf_tex
\begingroup%
  \makeatletter%
  \providecommand\color[2][]{%
    \errmessage{(Inkscape) Color is used for the text in Inkscape, but the package 'color.sty' is not loaded}%
    \renewcommand\color[2][]{}%
  }%
  \providecommand\transparent[1]{%
    \errmessage{(Inkscape) Transparency is used (non-zero) for the text in Inkscape, but the package 'transparent.sty' is not loaded}%
    \renewcommand\transparent[1]{}%
  }%
  \providecommand\rotatebox[2]{#2}%
  \newcommand*\fsize{\dimexpr\f@size pt\relax}%
  \newcommand*\lineheight[1]{\fontsize{\fsize}{#1\fsize}\selectfont}%
  \ifx\svgwidth\undefined%
    \setlength{\unitlength}{225bp}%
    \ifx\svgscale\undefined%
      \relax%
    \else%
      \setlength{\unitlength}{\unitlength * \real{\svgscale}}%
    \fi%
  \else%
    \setlength{\unitlength}{\svgwidth}%
  \fi%
  \global\let\svgwidth\undefined%
  \global\let\svgscale\undefined%
  \makeatother%
  \begin{picture}(1,0.16666667)%
    \lineheight{1}%
    \setlength\tabcolsep{0pt}%
    \put(0,0){\includegraphics[width=\unitlength,page=1]{gamma.pdf}}%
    \put(0.85518565,0.03259234){\makebox(0,0)[lt]{\lineheight{1.25}\smash{\begin{tabular}[t]{l}$\R$\end{tabular}}}}%
    \put(0.16592566,0.02703835){\makebox(0,0)[lt]{\lineheight{1.25}\smash{\begin{tabular}[t]{l}$\gamma$\end{tabular}}}}%
    \put(0.44555634,0.03481366){\makebox(0,0)[lt]{\lineheight{1.25}\smash{\begin{tabular}[t]{l}$0$\end{tabular}}}}%
    \put(0,0){\includegraphics[width=\unitlength,page=2]{gamma.pdf}}%
  \end{picture}%
\endgroup%

%% file: pics/locpic.pdf_tex
\begingroup%
  \makeatletter%
  \providecommand\color[2][]{%
    \errmessage{(Inkscape) Color is used for the text in Inkscape, but the package 'color.sty' is not loaded}%
    \renewcommand\color[2][]{}%
  }%
  \providecommand\transparent[1]{%
    \errmessage{(Inkscape) Transparency is used (non-zero) for the text in Inkscape, but the package 'transparent.sty' is not loaded}%
    \renewcommand\transparent[1]{}%
  }%
  \providecommand\rotatebox[2]{#2}%
  \newcommand*\fsize{\dimexpr\f@size pt\relax}%
  \newcommand*\lineheight[1]{\fontsize{\fsize}{#1\fsize}\selectfont}%
  \ifx\svgwidth\undefined%
    \setlength{\unitlength}{300bp}%
    \ifx\svgscale\undefined%
      \relax%
    \else%
      \setlength{\unitlength}{\unitlength * \real{\svgscale}}%
    \fi%
  \else%
    \setlength{\unitlength}{\svgwidth}%
  \fi%
  \global\let\svgwidth\undefined%
  \global\let\svgscale\undefined%
  \makeatother%
  \begin{picture}(1,0.75)%
    \lineheight{1}%
    \setlength\tabcolsep{0pt}%
    \put(0,0){\includegraphics[width=\unitlength,page=1]{locpic.pdf}}%
    \put(0.27126656,0.28119091){\makebox(0,0)[lt]{\lineheight{1.25}\smash{\begin{tabular}[t]{l}$M$\end{tabular}}}}%
    \put(0.24344391,0.37790009){\makebox(0,0)[lt]{\lineheight{1.25}\smash{\begin{tabular}[t]{l}$B$\end{tabular}}}}%
    \put(0.19816366,0.48099197){\makebox(0,0)[lt]{\lineheight{1.25}\smash{\begin{tabular}[t]{l}$L$\end{tabular}}}}%
    \put(0.45368622,0.46172022){\makebox(0,0)[lt]{\lineheight{1.25}\smash{\begin{tabular}[t]{l}$\phi$\end{tabular}}}}%
    \put(0.74716446,0.51063326){\makebox(0,0)[lt]{\lineheight{1.25}\smash{\begin{tabular}[t]{l}$T_\gamma$\end{tabular}}}}%
    \put(0.82041588,0.43123819){\makebox(0,0)[lt]{\lineheight{1.25}\smash{\begin{tabular}[t]{l}$\C^n$\end{tabular}}}}%
    \put(0.71644615,0.32372398){\makebox(0,0)[lt]{\lineheight{1.25}\smash{\begin{tabular}[t]{l}$\pi$\end{tabular}}}}%
    \put(0.70434571,0.14329392){\makebox(0,0)[lt]{\lineheight{1.25}\smash{\begin{tabular}[t]{l}\small {$\gamma$}\end{tabular}}}}%
    \put(0.50756138,0.13090736){\makebox(0,0)[lt]{\lineheight{1.25}\smash{\begin{tabular}[t]{l}$\R$\end{tabular}}}}%
    \put(0.71502838,0.22826084){\makebox(0,0)[lt]{\lineheight{1.25}\smash{\begin{tabular}[t]{l}$i\R$\end{tabular}}}}%
    \put(0.74385635,0.06001892){\makebox(0,0)[lt]{\lineheight{1.25}\smash{\begin{tabular}[t]{l}$\C$\end{tabular}}}}%
  \end{picture}%
\endgroup%

%% file: pics/gamm.pdf_tex
\begingroup%
  \makeatletter%
  \providecommand\color[2][]{%
    \errmessage{(Inkscape) Color is used for the text in Inkscape, but the package 'color.sty' is not loaded}%
    \renewcommand\color[2][]{}%
  }%
  \providecommand\transparent[1]{%
    \errmessage{(Inkscape) Transparency is used (non-zero) for the text in Inkscape, but the package 'transparent.sty' is not loaded}%
    \renewcommand\transparent[1]{}%
  }%
  \providecommand\rotatebox[2]{#2}%
  \newcommand*\fsize{\dimexpr\f@size pt\relax}%
  \newcommand*\lineheight[1]{\fontsize{\fsize}{#1\fsize}\selectfont}%
  \ifx\svgwidth\undefined%
    \setlength{\unitlength}{854.91868387bp}%
    \ifx\svgscale\undefined%
      \relax%
    \else%
      \setlength{\unitlength}{\unitlength * \real{\svgscale}}%
    \fi%
  \else%
    \setlength{\unitlength}{\svgwidth}%
  \fi%
  \global\let\svgwidth\undefined%
  \global\let\svgscale\undefined%
  \makeatother%
  \begin{picture}(1,0.09985081)%
    \lineheight{1}%
    \setlength\tabcolsep{0pt}%
    \put(0,0){\includegraphics[width=\unitlength,page=1]{gamm.pdf}}%
    \put(0.18669118,0.03605985){\makebox(0,0)[lt]{\lineheight{1.25}\smash{\begin{tabular}[t]{l}$\gamma$\end{tabular}}}}%
    \put(0.61082607,0.02759965){\makebox(0,0)[lt]{\lineheight{1.25}\smash{\begin{tabular}[t]{l}$\gamma'$\end{tabular}}}}%
    \put(0,0){\includegraphics[width=\unitlength,page=2]{gamm.pdf}}%
    \put(0.45289936,0.02540159){\makebox(0,0)[lt]{\lineheight{1.25}\smash{\begin{tabular}[t]{l}$\text{Mutation}$\end{tabular}}}}%
  \end{picture}%
\endgroup%

%% file: pics/monosec.pdf_tex
\begingroup%
  \makeatletter%
  \providecommand\color[2][]{%
    \errmessage{(Inkscape) Color is used for the text in Inkscape, but the package 'color.sty' is not loaded}%
    \renewcommand\color[2][]{}%
  }%
  \providecommand\transparent[1]{%
    \errmessage{(Inkscape) Transparency is used (non-zero) for the text in Inkscape, but the package 'transparent.sty' is not loaded}%
    \renewcommand\transparent[1]{}%
  }%
  \providecommand\rotatebox[2]{#2}%
  \newcommand*\fsize{\dimexpr\f@size pt\relax}%
  \newcommand*\lineheight[1]{\fontsize{\fsize}{#1\fsize}\selectfont}%
  \ifx\svgwidth\undefined%
    \setlength{\unitlength}{225bp}%
    \ifx\svgscale\undefined%
      \relax%
    \else%
      \setlength{\unitlength}{\unitlength * \real{\svgscale}}%
    \fi%
  \else%
    \setlength{\unitlength}{\svgwidth}%
  \fi%
  \global\let\svgwidth\undefined%
  \global\let\svgscale\undefined%
  \makeatother%
  \begin{picture}(1,0.66666667)%
    \lineheight{1}%
    \setlength\tabcolsep{0pt}%
    \put(0,0){\includegraphics[width=\unitlength,page=1]{monosec.pdf}}%
    \put(0.02037033,0.23296367){\makebox(0,0)[lt]{\lineheight{1.25}\smash{\begin{tabular}[t]{l}$\text{Disc with removed sector}$\end{tabular}}}}%
    \put(0.537985,0.264171){\makebox(0,0)[lt]{\lineheight{1.25}\smash{\begin{tabular}[t]{l}$\text{Sector}$\end{tabular}}}}%
    \put(0,0){\includegraphics[width=\unitlength,page=2]{monosec.pdf}}%
  \end{picture}%
\endgroup%

%% file: Sections/SFT.tex
\section{SFT Compactness} \label{sec:SFT}

In this section, we will explain a version of the neck-stretching procedure where the Lagrangian passes through the neck region and prove a compactness result with assumptions on the Legendrian that models the Lagrangian in the neck region.

Bourgeois, Eliashberg, Hofer, Wysocki and Zehnder \cite{bhewz} proves compactness result for closed holomorphic curves under neck-stretching for Morse-Bott Reeb dynamics. Abbas \cite{Abbasbook} proves compactness for open holomorphic curves with boundary on cylindrical Lagrangians  for the case of non-degenerate Reeb dynamics. We prove a version of the compactness result for neck-stretching along a contact type hypersurface $Z$ with a  free $S^1$ action that generates the Reeb field  for discs with boundary on a cylindrical Lagrangian that intersects $Z$ at the Legendrian $\Lambda$ in the neck region such that $\Lambda$ is a finite cover of an embedded Lagrangian in $Z/S^1$. Notice that the case of stretching along the boundary of a Darboux ball where we perform a local mutation falls in this case since the boundary sphere $S^{2n-1}$ has a natural $S^1$ action that generates the Reeb field, the Lagrangian intersects the sphere in two disjoint tori that descends to the Clifford torus in $\P^{n-1}$.

\subsection{Setup}\label{subsec:sftsetup}
We begin by setting the notations and conventions for this section. Let $(W,\omega,J)$ be a closed symplectic manifold and $Z$ be a compact codimension 1 submanifold of contact-type, i.e. in a tubular neighbourhood  $N \cong (-\epsilon,\epsilon) \times Z$, there is a contact form $\lambda$ on $Z$ such that  $\omega = d(e^{kt}\lambda)$ where $t\in(-\e,\e)$ and $k$ is a positive integer. In our case, we assume $Z$ separates $W$ into two connected components for simplifying the combinatorial structure, but this is not a necessary assumption. Let $L$ be a Lagrangian such that in the tubular neighbourhood of the hypersurface $Z$, $L$ is translation invariant over a Legendrian, i.e. there is a Legendrian $\Lambda \sub Z$ such that $N \cap L = (-\e,\e) \times \Lambda$, we call such Lagrangians to be cylindrical near $Z$. Assume $J$ is a compatible complex structure on $(W,\omega)$ such that it is translation invariant on $N$ and is adjusted to the form $\omega$ in the SFT sense, i.e., $J\frac{\del}{\del t} = R$ where $R$ is the Reeb field of $\lambda$ . Call the metric induced by the compatible pair $(\omega,J)$ on $W$ to be $g_W.$ Assume there is a free $S^1$ action on $Z$ that generates the Reeb field, the quotient $Y=Z/S^1$ has a natural symplectic form induced from $d\lambda$. Assume $\Lambda$ is closed under multiplying by $-1$ from the $S^1$ action (i.e. $-\Lambda = \Lambda$), projects to an embedded Lagrangian $L_Y$ in $Y$ and the projection map is a finite cover.

\begin{figure}[ht]
    \centering

        \def\svgscale{1.5}
    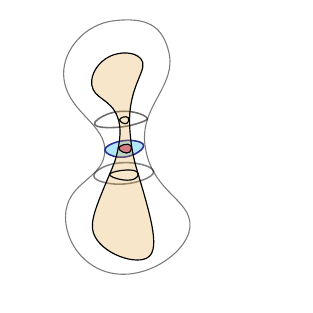

    \caption{Symplectic manifold $W$ and Lagrangian submanifold $L$ with neck pieces of length $2\e$}
    \label{fig:sftsetup}
\end{figure}

\begin{example}
A sphere $S^{2n-1}_r$ of radius $r$ centred at $0$ inside $\C^n$ has a tubular neighbourhood $\Phi : (-\e,\e) \times S^{2n-1}_r \to \C^n $ given by $\Phi(s,p)= e^sp$. Clearly $\Phi^*J_0$ is a translation invariant almost complex structure as scaling is a biholomorphism in $\C^n$. Moreover, $\Phi^*\omega_0 =  d(e^{2s}\lambda)$ where $\lambda = \frac{\lambda_0}{\sum |x_i|^2 + |y_i|^2}$. Here $\omega_0,\lambda_0$ are the standard symplectic and contact forms respectively. 

Note that $T_\gamma$ for an admissible $\gamma$ is cylindrical over disjoint union of two Legendrian tori. This follows from the fact that $f_a(T_z)) = T_{a^nz}$ where $f_a$ is the scaling map, multiplication by $a$.
\end{example}

\subsubsection{Neck stretching}
We now explain the process of \textit{neck stretching}. The $\tau-$stretched manifold $W^\tau$ is defined by removing the tubular `neck' around $Z$ and replacing with a $2\tau$ sized neck as defined below - 
$$W^\tau =  \underbrace{(W \setminus N ) \cup   (-\tau,\tau)\times Z} _{\text{ gluing $Z$ at the ends}}. $$
Since $J$ is translation invariant on $ N$, $W^\tau$ can be given a complex structure $J^\tau$ adjusted to $\omega$   by requiring $J^\tau$ on the extending the translation invariant complex structure $J|_N$ . We define the $\tau-$ stretched Lagrangian $L^\tau$ similarly as $W^\tau$.
$$L^\tau =  \underbrace{(L \setminus N ) \cup   (-\tau,\tau)\times \Lambda} _{\text{ gluing $\Lambda$ at the ends}} .$$

\noindent From our assumption, the manifold with the hypersurface removed $W\setminus Z$ has two connected components. Thus if we remove the neck $N$ from the manifold $W$, there are two connected components of $W \setminus N$,   call the connected components  $W_{in}$, $W_{out}$. Clearly, all the $\tau$ stretched manifolds are diffeomorphic to $W$ by choosing a diffeomorphism between the intervals $(-\e,\e)$ and $(-\tau,\tau)$. The metric on the stretched manifolds is defined to be $g_W$ away from the necks and cylindrical extension of the compatible metric on $Z$ in the neck piece.

\subsubsection{ Energy}
For establishing compactness of a family of objects, it is essential for the family to be bounded in some sense.  The notion of \textit{energy} has been used to bound families of $J-$holomorphic curves, for SFT, there is a similar notion of \textit{energy} which we describe now. 

There are more than one way of defining energy for the neck-stretched manifolds. Cielibak-Mohnke \cite{ciel} defines it to be the symplectic area induced from a particular choice of diffeomorphism between $W$ and $W^\tau.$  
Bourgeois et. al. \cite{bhewz} defines energy  analogous to Hofer energy. We will use the approach in op. cit. for dealing with energy , see section 9.2 in \cite{bhewz}.

We recall the definition of Hofer energy. Given a symplectic manifold $M$ with a cylindrical end $\R_+ \times Z$ where $Z$ is a contact type hypersurface i.e. there is a contact form $\lambda$ on $Z$ such that $\omega=d(e^t\lambda)$.  Denote the compact complement of the cylindrical end by $\overline{M}$. The Hofer energy is of a curve $u:\Sigma \to M$  is defined as follows : 

$$E_H(u) = \int_{u\inv (\bar{M})} u^* \omega   + \text{sup}_{\phi \in \mathcal{T}} \int_{u\inv (\R_+ \times Z)} u^* d(e^{\phi(t)} \lambda),  $$
where $\mathcal{T}$ is the family of increasing diffeomorphism from $\R$ to $(0,1)$ .

The horizontal energy $E_{Hor}$ refers to the energy of the projection to $Z$ in the cylindrical end along with the energy away from the cylindrical end. The formula for  the horizontal energy is  $$E_{Hor}(u) = \int_{u\inv (\bar{M})} u^* \omega   +  \int_{u\inv (\R_+ \times Z)} u^* p_Z^* d \lambda .$$ Here $p_Z$ is the projection on to $Z$. It is clear that  finiteness of Hofer energy implies finiteness of the horizontal energy. In SFT literature this is commonly called the  $E_\omega$ energy, but we avoid this notation since for us $\omega$ is a symplectic form and not a part of stable Hamiltonian structure. The energy for the stretched manifold is defined similarly, let $u:\Sigma \to W^\tau$ be a curve, the energy is defined by $$E(u) =  \int_{u\inv (\bar{W})} u^* \omega   + \text{sup}_{\phi \in \mathcal{T}} \int_{u\inv ((-\tau,\tau) \times Z)} u^* d(e^{\phi(t)} \lambda) .$$
where $\bar{W}$ is the compact manifold obtained by removing the tubular neighborhood $N$ from $W$ and  $\mathcal{T}$ is the family of increasing diffeomorphism from $(-\tau,\tau)$ to $(0,1)$.

\begin{remark} \label{rem:horizontalmajorizes}
    Lemma 9.2 in \cite{bhewz} proves that the horizontal energy $E_{Hor}$ can be used to bound the total energy of holomorphic maps to $( W^\tau, J^\tau )$. More precisely, there exists a $C>0$ depending only on $W,Z,\lambda$ and not on $\tau$ such that $E(u) < CE_{Hor}(u)$ for all holomorphic $u : S \to (W^\tau,J^\tau)$. Later, this lets us use $E_{Hor}$ bounds to compactify moduli of maps to neck-stretched manifold.  
\end{remark}

\begin{figure}[ht]
    \centering

        \def\svgscale{1}
    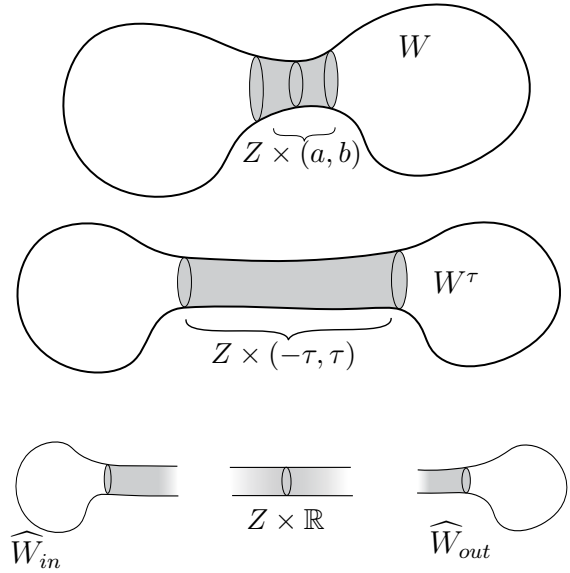

    \caption{Breaking of $W$}
    \label{fig:breaking}
\end{figure}

\subsubsection{Broken manifold and broken discs}\label{subsub:brokenstuff}

The geometric objects that are obtained as a result of neck stretching can be intuitively imagined to be broken manifolds, which we explain in detail below. As before, assume $(W,\omega,J,L)$ be a symplectic manifold with a contact type hypersurface, $Z$ such that $J,L$ are cylindrical near the hypersurface $Z$. Conceptually as $\tau \to \infty$, the neck of the manifold gets longer and longer and one can view the `limit' to be given by adding positive and negative symplectization of $Z$ to the two boundaries in $W\setminus Z$. If we further assume that the hypersurface $Z$ divides the manifold $W$ into two connected components $\bar{W}_{in}$, $\bar{W}_{out}$ with contact type boundaries, a  broken manifold $\mathbb{W}[k]$ with $k$ necks is given by the collection  $$\widehat{W_{in}}, \underbrace{ \R\times Z,\R\times Z,\R\times Z,\dots  }_{\text {k times}} , \widehat{W_{out}},$$
where $\widehat{X}$ denotes symplectization of the contact type boundary of the manifold $X$ and $W_i$ refers to either the $i\th$ neck piece or the corresponding symplectizations. When we want to refer to a broken manifold without specifying the number of necks, we use the notation $\mathbb{W}$.

Similar to a broken symplectic manifold, one has a notion of broken Lagrangian submanifold in a broken symplectic manifold. We can construct $\mathbb{L}$ from a Lagrangian $L$ if we  assume that $L$ is cylindrical near the hypersurface $Z$. We define the broken Lagrangian $\mathbb{L}[k]$ with $k-$ necks to be the collection

$$\widehat{L_{in}}, \underbrace{ \R\times \Lambda,\R\times \Lambda,\R\times \Lambda,\dots  }_{\text {$k$ times}} , \widehat{L_{out}},$$
where $L_{in}$ and $L_{out}$ denote the intersection of $L$ with $W_{in},W_{out}$ respectively and the hat $\widehat\_$ denotes the extension to the symplectization by the cylinder of the Legendrian $\Lambda$ and
$L_i$ refers to the $i\th$ neck piece.

 \begin{defn}[Broken disc]\label{def:brokdisc}
 
 A broken disc with $k-$necks $(\mathfrak{u,D})$ is a map from a  nodal disk with marked points $\mathfrak{D}$ to a broken manifold with $k$ necks with the following properties. 
 \begin{itemize}
     \item $\mathfrak{D}$ is a nodal disc with labels $l\in \{\text{in}=0,1,2,\dots,k,\text{out}=k+1\}$ such that neighboring discs have labels differing at most by 1. We designate some nodes to be ``puncture nodes" that are the nodes in the disc  between two components with different labels.We write the nodal disks as  $ \mathfrak{D} = (\D_{in},\D_1,\dots,\D_{out})$ where $\D$ refers to a collection of discs.
     \item  $ \mathfrak{u} = (u_{in},u_1,\dots,u_{out})$ is a map from $\mathfrak{D}$ to $\mathbb{W}[k]$ such that  $u_i$ is a map from $\D_i$ to $W_i$ with boundary on $L_i$ and near the puncture nodes of $\mathfrak{D}$ between components of differing label $u_i$ and $u_{i+1}$ asymptote to the same Reeb chord or Reeb orbit.
     \item (\textit{Stability}) for all $i \in [1,\dots,k]$, at least one of the component in each neck $\D_i$ has to be a non-trivial strip, or have enough markings to make it stable.
 \end{itemize}
  \end{defn}
\begin{figure}
    \centering

        \def\svgscale{1.5}
    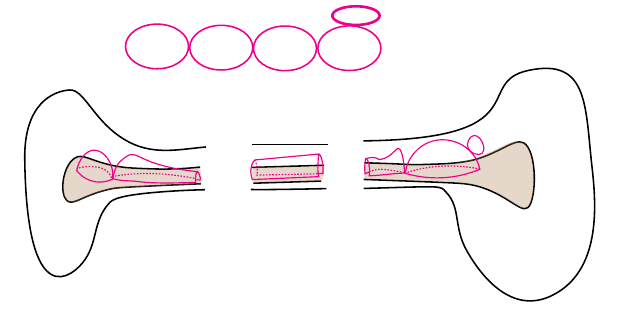

    \caption{Broken disc in $\mathbb{W}[1]$}
    \label{fig:my_label}
\end{figure} 

\begin{remark}
The labels in and out are chosen to refer to the ``inside" and ``outside" where inside means where we do a local mutation and outside means where we don't make any changes.
\end{remark}

Now that we have defined broken discs, we will proceed to define the notion of convergence of broken discs.
We say a sequence of discs $u_n : \D^2 \to (W^n,J^n,L^n)$ converge to a broken disc $(\mathfrak{u,D})$ with $k-$necks if the following holds: \begin{itemize}
    \item there is a sequence of maps $\phi_n:\mathfrak{D_n} \to \mathfrak{D} $ such that they are biholomorphic diffeomorphisms except away from a certain circle or line that might get mapped to a node in $\mathfrak{D}$
    \item there is a sequence $t^{n}_i \in \R$ such that $u^n_i \of \phi_n\inv+t^n_i$ converges in $C^{\infty}$ to $u_i$ on every compact subset. Here $+t$ refers to translation in the $\R$ direction by $t$ 
\end{itemize}

\noindent Similarly, we say a sequence of broken disks  $(\mathfrak{u^n,D_n})$ with $k$ necks converges to a broken disk $(\mathfrak{u,D})$ with $k'>k$ necks if the following hold:
\begin{itemize}
    \item There is a sequence of maps $\phi_n:\mathfrak{D_n} \to \mathfrak{D} $ such that they are biholomorphic diffeomorphisms, except away from certain circles or lines that might get mapped to a node in $\mathfrak{D}$.
    \item There is a sequence $t^{n}_i \in \R$ such that $u^n_i \of \phi_n\inv+t^n_i$ converges in $C^{\infty}$ to $u_i$ on every compact subset. Here $+t$ refers to translation in the $\R$ direction by $t$ .
  
\end{itemize}

\subsection{Asymptotic convergence to Reeb chords}

In this subsection, we show that we can obtain asymptotic data near the puncture nodes if the broken disc has finite Hofer energy. To that end, we prove a version of asymptotic convergence of pseudoholomorphic strips in $\R \times Z$ with finite Hofer energy  with boundary on $\Lambda$.

We start with recalling the assumptions on the contact hyper surface and the complex structures. $(Z,\lambda)$ is a closed contact manifold  with contact form $\lambda$ such that there is a free $S^1$ action on $Z$ that generates the Reeb field\footnote{In literature such contact spaces are also called pre-quantum bundles} i.e.\ the quotient $Y=Z/S^1$ has an natural symplectic form induced from $d\lambda$. Let $\Lambda$ be a Legendrian in the contact manifold $Z$ that is closed under action of  $-1$ from the $S^1$ action (i.e. $-\Lambda = \Lambda$). Also assume that the Legendrian $\Lambda$ projects to an embedded Lagrangian $L_Y$ in $Y$ and the projection map is a finite cover. Clearly, the Reeb dynamics of the contact form $\lambda$ is Morse-Bott in the sense of \cite{bourg}.

We make the following assumptions on the complex structure.  $J$ is an almost complex structure on $\R \times Z$ such that $J$ is translation invariant and is adjusted to the contact form $\lambda$. Thus, the almost complex structure $J$ descends to a compatible almost complex structure on $(Y,d\lambda)$ that makes the projection $ \pi_Y : R\times Z \to Y$ a holomorphic map. The complex structure $J$ induces a cylindrical metric on $\R \times Z$, call that metric $g_{cyl}$.

Now, we explain our goal of this subsection in details. We begin by explaining the strips that model the asymptotics of a holomorphic strip with finite Hofer energy. We define the map $u_\gamma$ on a  strip as follows 
$$u_\gamma(s,t) = (ms,\gamma(mt)),$$ where $ \gamma $ is a Reeb chord with end points on $\Lambda$ and $m>0$ is the length of the Reeb chord $\gamma$. From the definition, we see that $u_\gamma$ is invariant of the $\R$ action on the target $\R \times Z$. We call such translation invariant strips over Reeb chords \textit{trivial strips}.

We can use classical removal of singularity and exponential convergence results of $J$ holomorphic strips to conclude exponential convergence of strips in the compact manifold  $Y$.  Since $L_Y$ is a compact embedded Lagrangian in  the symplectic manifold $Y$, we have that there is a number $\delta_Y>0$ such that any finite-energy pseudoholomorphic half-strip $ v : \halfstrip \to Y $ with boundary on $L_Y$, the strip $v$ converges to a point in the Lagrangian $L_Y$ exponentially fast with speed $O(e^{-\delta_Y s})$ where $(s,t)$ is the standard coordinate on $\halfstrip$. Rigorously, there is a $C>0$ and $p \in L_Y$ such that $d(v(s,t),p) < Ce^{-\delta_Y s}$.  Let $i_\Lambda = \inf  \,(\, \{ \delta_Y\} \cup \{ \theta \in (0,2\pi)| \, \, \exists p\in \Lambda \, \, \text {such that }  e^{i\theta}.p \in \Lambda \})$. Clearly $i_\Lambda > 0$ since $\Lambda$ is a finite cover of $L_Y$ and $L_Y$ is compact.

\begin{Theorem}[Exponential Convergence] \label{thm:exponential conv}
For any $\delta \in (0,i_\Lambda) $ and holomorphic strip $u:[0,\infty) \times \interval \to \R \times Z$ , that satisfies the following,
\begin{itemize}
    \item $E_H(u)<\infty$ ,
    \item $u$ satisfies the boundary condition  $u\, ( \,[0,\infty) \times \{0,1\}) \sub\R \times \Lambda$,
    \item the projection on $\R$ , $\pi_\R (u)$ is unbounded;
\end{itemize}  there is an $s_0 > 0$, a trivial strip $u_\gamma$ over a Reeb chord $\gamma$ and a $C>0$ such that

\begin{equation}
\label{expConv}
     d_{cyl}(u(s,t),u_\gamma(s,t)) < Ce^{-\delta s} \quad \forall s\geq s_0.
\end{equation} 

\end{Theorem}

\noindent We will break the proof of the theorem in the following steps.
\begin{itemize}
    \item Step 1 : Compactify $\R \times Z$ to a symplectic manifold $N$ with compatible almost complex structure such that $\R \times \Lambda$ compactifies to a cleanly self intersecting Lagrangian $L$ in $N$.
    \item Step 2 : Show that $\pi_\R (u)$ diverges to $\pm \infty$ where $\pi_\R$ is the projection of $\R \times Z$ to $\R$.
    \item Step 3 : Extend $u$ to a map  $\bar{u}$ to $N$ with boundary on $L$.
    \item Step 4 : Use analysis done in \cite{schm} to conclude convergence in $N$.
    \item Step 5 : Get convergence in cylindrical metric from convergence in $N$.
\end{itemize}

\begin{proof}

\begin{figure}[ht]
    \centering

        \def\svgscale{1.2}
    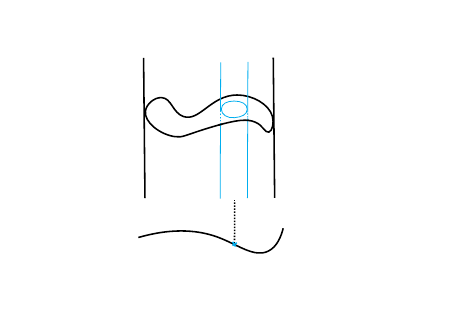

    \caption{$\R \times Z$ as a $\C^*$ fibration}
    \label{fig:cylfib}
\end{figure}

\textbf{Step 1 : Compactification of target} 

We will compactify our target manifold by adding two copies of the symplectic manifold $Y$. Since $Z$ is a principal circle bundle over $Y$, the cylinder $\R \times Z$ is a $\R \times S^1 $ bundle on $Y$. The transition functions  between trivializations of $\R \times Z$ over the base $Y$ is given by rotation action of $S^1$, thus the cylinder $\R \times Z$ can be viewed as a $\C^*$ bundle over $Y$. Moreover, we can extend this bundle to a $\C\P^1$ bundle over $Y$ by compactifying each fiber to a $\C\P^1$ and using the induced transition functions, call this compactification $N$. Thus, $N$ is a $\C\P^1$ bundle over the base $Y$. Clearly $N$ is compact as $Y$ is compact. Denote the zero section of $Y$  by $Y_0$ and the infinity section ($y\mapsto (y,\infty)$ in each trivialization) by $Y_\infty$.  The complex structure $J$ extends naturally to a complex structure $J_N$ on $N$ that restricts to the standard complex structure on the $\CP^1$ fibers. We can define a symplectic form $\omega_N$ that is compatible with $J_N$ by defining it in each trivialization $\C\P^1 \times U$ for $U$ open set in $Y$ as following, 
$$\omega_N =  \left( \begin{matrix} \omega_{FS} && 0 \\ 0 && d\lambda  \end{matrix} \right) ,$$
where $\omega_{FS}$ is the Fubini-Study form. The form $\omega_N$ is a well-defined symplectic form on $N$ as the transition functions of $N$ seen as a  $\CP^1$ fibration are given by the rotation action of the $S^1$ on $\C\P^1$. One can check from the construction of $\omega_N$ and $J_N$ that they are compatible from a direct computation.

\begin{figure}[ht]
    \centering

        \def\svgscale{1.2}
    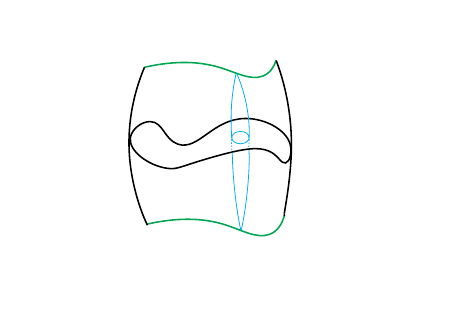

    \caption{Compactification of $\R \times Z$ to $N$}
    \label{fig:compcyl}
\end{figure}

We describe the closure of the cylindrical Lagrangian in the compactified symplectic manifold. Denote the closure of the cylindrical Lagrangian $  \bar {\R \times Z}$ in $N$ by $L$. We can check that $L$ is obtained by adding $L_Y$ in the 0 and infinity divisor of the $\CP^1$ fibration to $\R \times \Lambda$. Indeed, a trivialization $\CP^1 \times  U$ of $N$ over some open $U$ in $Y$ such that $U\cap L_Y$ is nonempty, we see that $$L \cap ( \CP^1 \times  U) = \{0,\infty\}\times (U\cap L_Y)  \cup ( ( \CP^1 \times U) \cap ({\R \times \Lambda})).  $$  
This directly shows adding two copies of the Lagrangian $L_Y$ to the cylinder $\R\times Z$ produces the closure $L$. Since the Legendrian $\Lambda$ is closed under the action of $-1$, we see that $L$ cleanly self intersects at $L_Y$ in the 0 and infinity divisor. Call the intersection loci $L_{Y_0}$, $L_{Y_\infty}$ respectively.

\textbf{Step 2 : Divergence of $\pi_\R (u)$}

We now show that near the ends of the strip, the image escapes to infinity.  We will prove that in the limit as $s \to \infty$ ,  $u$ diverges to either $\infty$ or $-\infty$. The proof employs a version of monotonicity lemma (Proposition 2.69 \cite{Abbasbook}) and gradient-boundedness of finite Hofer energy maps. 

We state the monotonicity lemma without proof.
\begin{lemma}[Proposition 2.69 \cite{Abbasbook}]\label{lemm:gromMon}
    Let $(W,\omega)$ be a compact symplectic manifold, $E_0 >0$, let $J$ be an $\omega$-compatible almost complex structure and let $L \sub W$ be a Lagrangian submanifold. Let $\mathbb D = \big\{ z \in \C \big| |z| < 1 \big\}, N \in \N $ and $I_k = \{ e^{it} | \frac{2\pi}{N}(k-1) \leq t \leq \frac{2\pi}{N}k \} \sub \del \mathbb D $ for $k= 1,\dots N.$ Then there are constants $\delta_0,c_0 > 0 $ such that the following holds: 
    
    Let $u:\ol \D \to W$ be a $J$-holomorphic disk with $E(u)\leq E_0 $ so that $u(I_k)\sub L$ for $k$ even. Let $z_0 \in \D \setminus \cup_{k} I_{2k+1}$ and $0< r < \delta_0$ so that $u(\cup_k I_{2k+1}) \sub W\setminus \ol{B_r(u(z_0))}.$ Then $$\int_{u(\ol D) \cap B_r(u(z_0))} \omega \geq c_0 r^2.$$
\end{lemma}

Another ingredient in the proof of divergence of $\pi_\R$ is a slight modification of Proposition 2 in \cite{erkaoAsymp}.  Define the norm $\| \|_{L^\infty}$ for a map $u:\mathbb H \to \R\times Z$ as follows $$\| \nabla u \|_{L^\infty} := \sup_{(s,t)\in \mathbb H} | \nabla u (s,t)|,$$ where the norm $|\cdot| $ is taken with respect to the standard Euclidean metric on $\mathbb H$ and cylindrical metric on $\R \times Z$. 

\begin{lemma}\label{lemm:prop1erk}
    If $v: \mathbb H \to \R \times Z$  is a holomorphic map with boundary on $\R \times \Lambda$ with zero horizontal energy, $E_\omega (v) = 0$, and bounded derivative, $\| v\|_{L^{\infty}} < \infty$, then $v$ is a constant map.
\end{lemma}

\begin{proof}
    Since $E_\omega(v)=0$, we have that $\pi_Y (v)$ is a constant map. Thus, the map $v$ lies in a fiber $\R \times S^1$. The map $v$ maps the boundary $\R \sub \mathbb H$ to $\R \times \{ * \}$ for some $* \in S^1$ . Since $\mathbb H$ is simply connected, $v$ factors as follows $v = e^{h}$ where $h: \mathbb H \to \C$ is a holomorphic map. The map $h$ sends the boundary $\R \sub \mathbb H $ to $\R + iy$ for some $y\in \R$.  We can use a doubling trick using Schwarz reflection to extend $h$ to $\C$, thus we can extend $v$ to a map from $\C$. The rest of the proof follows from Proposition 1 in \cite{erkaoAsymp}.
\end{proof}

\begin{lemma}[Proposition 2 \cite{erkaoAsymp}] \label{lemm:gradbound}
    Let $u : \R\times \interval \to (\R_+ \times Z, \R_+ \times \Lambda)$ be a holomorphic map with finite Hofer energy, $E_H(u) < \infty$, then we have $\|\nabla u \|_{L^\infty} < \infty $. 
\end{lemma}
\begin{proof}
    The proof is the same as the proof of Proposition 2 in \cite{erkaoAsymp}. The only change one needs to do is to use Lemma \ref{lemm:prop1erk} in place of Proposition 1 when the rescaled maps converge to a map on $\mathbb H$. 
\end{proof}



The idea behind the proof of  divergence is that if a semi-infinite strip $u$ does not diverge to infinity, it needs to keep coming back to some compact subset of $\R \times Z$. The monotonicity lemma states that every time the strip comes back into such a compact set, there is a `cost of energy', hence a strip with finite energy can't keep coming back to any compact subset. We now state and prove the required divergence.

\begin{lemma}
\label{lemm:div2inf}
$\lim_{s\to \infty} |\pi_\R (u(s,t)) | = \infty$
\end{lemma}

\begin{proof}

Assume the contrary, that there exists $M > 0$ such that there is a sequence $(s_n,t_n)$ that satisfies $|\pi_\R (u(s_n,t_n)) | < M$ and $s_n \to \infty$. Since $\pi_\R$ is unbounded, there is also a sequence $(s'_n,t'_n)$ such that $s'_n \to \infty$ and $\pi_\R \of u (s'_n,t'_n) \to \infty$. 
In addition, we can assume $s_{n+1} > s'_n > s_n$ by taking subsequences.

From Lemma \ref{lemm:gradbound} we know that $|\nabla u | < C$ for some $C > 0$.
Set $W=[-M-C,M+C] \times Z$ and select a symplectic structure $\omega_f$ on $W$  of the form $d(e^{f(r)}\lambda)$ where the function $f:\R \to (0,1)$ is a diffeomorphism such that $f$ is linear on $[-M-1,M+1]$. Now use  monotonicity lemma \ref{lemm:gromMon} on the symplectic manifold with boundary $(W,\omega_f)$ with $N=4$ to get $c,R$. Choose $r_0 < \min \{ R,\frac{1}{2}\}$. Let $\Sigma_n  $ be sub-strips of the format $[\alpha_n,\beta_n] \times \interval$ where $s'_{n-1}<\alpha_n < s_n < \beta_n < s'_n$ such that $u(\alpha_n, \cdot)$ and $u(\beta_n, \cdot)$ lies in $[M+ 10C , \infty) \times Z$.  We can now apply the monotonicity lemma to each $\Sigma_n$. We have, for all $n$, \begin{equation}\label{eq:bla}
    \int_{\Sigma_n} u^*\omega_\phi \geq {c}{r_0^2}.
\end{equation} 
Then from summing over $n$ in equation (\ref{eq:bla}), we get $E_{\omega_{\phi}}(u) > \sum_1^\infty {c}{r_0^2}=\infty$ which contradicts that Hofer energy is finite since $E_H \geq E_{\omega_{\phi}}$. Thus, we have $\lim_{s\to \infty} |\pi_\R (u(s,t)) | = \infty$.

\end{proof}

\noindent Without loss of generality, we will assume from now on wards that the $\pi_\R(u)$ goes to $\infty$ as $s$ goes to $\infty$, ie. $\pi_\R (u) \to \infty $ as $s \to \infty$. The case of $\pi_\R(u) \to -\infty$ would follow by similar arguments.

\textbf{Step 3 : Extending $u$ in $N$} \label{comptarg} 

We will now show that the map $u$, when considered as a holomorphic map to $N$, can be extended over the punctures. Under the biholomorphism $z \mapsto e^{-\pi z}$, the half strip $[0,\infty) \times \interval  $ is mapped to the punctured half disc $D^+_0$ where half disc $D^+$ denotes $D\cap \mathbb{H}$ for a disc $D $ centered at $0$.  From this step on wards we view the domain as a punctured half disc using the identification $z \mapsto e^{-\pi z}$. We will show that the map $u$ can be extended by removing the boundary singularity at $0$. As of now, $u$ can be viewed as a map of the half strip to $N$ with boundary on the cleanly self-intersecting Lagrangian $L$. Note that $\pi_Y(u)$ is a holomorphic map from the half strip, such that the boundary of the strip lies on the Lagrangian $L_Y$. As the Hofer energy $E_H(u)$ is finite, we have that the horizontal energy $E_{Hor}(u)= E_{d\lambda} (\pi_Y (u))$ is finite, thus from classical removal of singularities with Lagrangian boundary in compact symplectic manifold, we have $\pi_Y(u)$ can be extended over the singularity at $0$, denote the point $\pi_Y(u(0)) \in L_Y$ as $q$. Assuming that $\pi_\R (u) \to \infty$ as $s\to \infty$, we can extend $u$ to $\D^+$ by defining $u(0) = q_\infty$ where $q_\infty$ is the image of $q$ in $Y_\infty$. This can be checked by taking a trivialization of $N$ over some neighborhood $U$ of $q$, then $u|_{[R,\infty)}$ for large $R>0$ is a map to $\CP^1 \times U$ as $\pi_Y(u)|_{[R,\infty)}$ lies in $U$. Since $\pi_\R (u(z)) \to \infty$ as $z\to 0$, we have that $u(z) \to (\infty,q)$ as $z \to 0$. Thus we have shown that $u$ can be extended. This gives us the following lemma.

\begin{lemma}
If $u:[0,\infty) \times [0,1] \to \R \times Z$ with boundary on $\R \times \Lambda$, has finite Hofer energy $E_H(u) < \infty$, then $u$ has a limit as $s \to \infty$  in the compactification $N$ of $\R \times Z$ and also finite energy $E_{\omega_N} (u) < \infty$ with respect to the symplectic form $\omega_N$.
\end{lemma}

\textbf{Step 4 : Exponential convergence to a point in $N$}

We will use results from \cite{schm} to show the exponential convergence of $u$ to Reeb chords. In \textbf{Theorem} \textbf{3.1 } and \textbf{3.2} the exponential convergence of finite energy strips with boundary on clean intersecting Lagrangians is proved. 

To use Schm\"aschke's result, we first need to ensure that the map $u$ on the boundary of the strip lies on at most two Lagrangians that cleanly intersect. Since $\pi_Y \of u $ converges to $q$, there is a small neighborhood $V$ around $q_\infty$ in $N$ such that $L \cap V$ is diffeomorphic to cleanly intersecting Lagrangians and $u$ has boundary lying on at-most two of them, denote these Lagrangians as $L_0,L_1$ that cleanly intersect at $L_{Y_\infty} \cap V$.

Now, Schm\"aschke's\textbf{ Theorem 3.2} gives us an exponential convergence in terms of eigenvalue of an asymptotic operator. We have that $u(s,t) = \exp _{q_\infty}(e^{-\alpha s}\zeta(t) + w(s,t)) $ where  $\zeta (t)$ is an eigenfunction with eigenvalue $\alpha$ of the operator 
$$A_{q_\infty}:T_{q_\infty}\mathcal{P}(L_0,L_1) \to L^2([0,1],T_{q_\infty} N), \text{      } \xi \mapsto J_N\del_t\xi,$$
where $T_{q_\infty}\mathcal{P}(L_0,L_1) = W^{1,2}(T_{q_\infty}N,T_{q_\infty}L_0,T_{q_\infty}L_1)$, the Sobolev space of paths from $\interval$ with boundary on $T_{q_\infty}L_0,T_{q_\infty}L_1$. As $\Lambda$ is a finite cover of $L_Y$, we can arrange a trivialization $U\times \CP^1$ of $N$ over a neighborhood $U$ of $q$ such that $L_0,L_1$ are determined by two sections of $Z$ over $L_Y$. Call those sections $s_1,s_2 : L_Y \cap U \to S^1$. From definition of $i_\Lambda$ we see that $d(s_1,s_2) > i_\Lambda$ under the standard $S^1$ metric induced from $[0,2\pi]$. As  $J_N$ splits into a direct sum of $J_Y$ and $J_0$ in each trivialization, we obtain the eigenvalues and eigenvectors separately from the $Y$ and $\CP^1$ component in the trivialization $U \times \CP^1$. In the $\CP^1$ component we see that the eigenvalues of $A_{q_\infty}$ is determined by the difference in angle between $s_1(q)$ and $s_2(q)$ that is at most $i_\Lambda$. Denote the angle difference $s_1(q)/s_2(q)$ by $\theta_q$, then the eigenvalues are given by $2k\pi + \theta_q $ for all integers $k$. Since $J_N$ restricts  to the standard complex structure $J_0$ on $\CP^1$ a direct computation shows that the eigenfunctions are $e^{it(2k\pi+\theta_q)}v$ for $v \in T_{q_\infty}L_0$.

\begin{figure}[ht]
    \centering

        \def\svgscale{1.8}
    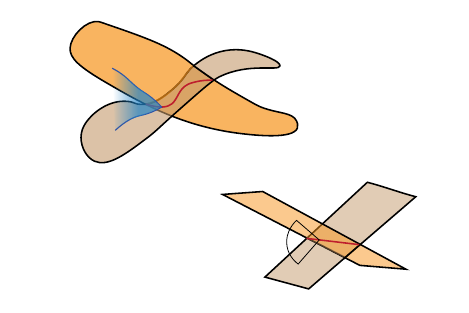

    \caption{Angles between the Lagrangians $L_0$,$L_1$ at $q_\infty.$}
    \label{fig:lagangle}
\end{figure}

We can now collect the results of the previous paragraph to obtain the required exponential convergence. When $\theta_q \neq \pi$, $v$ lies in the vertical component $T_{\infty}\CP^1$ of $T_{q_\infty}N = T_{\infty}\CP^1 \oplus T_{q}U$ in the trivialization $\CP^1\times U $ as the only eigenvalues in the horizontal part are $k\pi$. When $\theta_q = \pi$, the boundary of $u$ lies on a single Lagrangian branch. Thus, we can do a doubling trick to obtain a holomorphic cylinder in a neighborhood of $q_\infty$, then the exponential convergence  follows from \cite{bourg} since the Reeb dynamics in our case is Morse-Bott. From now we assume $\theta_q \neq \pi$, thus from Theorem 3.2 in \cite{schm} we get, for some vertical $v\in T_{\infty}\CP^1 \hookrightarrow T_{q_\infty}N  $ , $s_0>0$  and any $\delta \in (0,i_\Lambda)$,

 \begin{equation}
\label{Nconv}     
u(s,t) = \exp _{q_\infty}(e^{-(2k\pi + \theta_q) (s+it)}v + w(s,t)) , \,\|w\|_{C^l} = O (e^{- (\delta + 2k\pi + \theta_q)s})  \, \forall l\in \mathbb{N}, s\geq s_0.
 \end{equation}
 \blfootnote{Here  $ O(e^{- (\delta + 2k\pi + \theta_q)s})$ means that there is a $c_l >0$ such that $\| w\|_{C^l([s,\infty))} < C_le^{- (\delta + 2k\pi + \theta_q)s}$. We keep using the big-$0$ notation from now on to avoid notational clutter.}

\textbf{Step 5: Convergence to a strip in cylindrical metric}

So far, we have proved exponential convergence of a $J-$holomorphic strip in a compact manifold. In this step, we will prove exponential convergence in the cylindrical manifold $\R \times Z $ with respect to a cylindrical metric. We will use the convergence in the compact manifold $N$ to obtain exponential convergence in the cylindrical manifold $\R \times Z$. 

The choice of metric  in $N$ in (\ref{Nconv}) to define the exponentiation is not vital, since exponentiation is an analytic local diffeomorphism and derivative at $0$ is identity. From a Taylor-series computation, we see that changing the metric doesn't change the leading order decaying term in the equation (\ref{Nconv}). Let $\wt {\exp}_{q_\infty}$ denote exponentiation at $q_\infty$ with a different metric. The Taylor-series of $\wt \exp_{q_\infty} \inv \of  \exp_{q_\infty} $ is as follows,

\begin{align}\label{eqn:exptaylor}
    \wt \exp_{q_\infty} \inv \of  \exp_{q_\infty} (\wt v) &= d(\wt \exp_{q_\infty} \inv \of  \exp_{q_\infty})_0 (\wt v) + O(|\wt {v}| ^2) \nonumber \\
    &= {d{{\wt{\exp}_{q_\infty}\inv}}}_0 \of d{\exp_{q_\infty}}_0 (\wt v) + O(|\wt {v}| ^2) \nonumber \\
    &= \operatorname{Id} \of \operatorname{Id} (\wt v) +  O(|\wt {v}| ^2)
\end{align}
Thus, from equations (\ref{Nconv}) and (\ref{eqn:exptaylor}) we get that $u$ converges exponentially to the strip $\wt{\exp} _{q_\infty}(e^{-(2k\pi + \theta_q) (s+it)}v)$.

We choose a metric on the compactified manifold $N$ that lets us easily obtain trivial strips from geodesic exponentiation. Choose a metric in $g$ for $N$ such that in the trivialization $U\times \CP^1$ in a neighborhood $U$ of $q$, $g = g_1 \oplus g_2$ where $g_1$ is the restriction of some metric on $Y$ and $g_2$ is a metric on $\CP^1$ for which an open disc $B_\infty$ centered at $\infty$ is isometric to a disc $B_0$ in $\C$ with standard Euclidean metric. Since $v$ in (\ref{Nconv}) lies in the vertical component, we have that $$\exp_{q_\infty}((e^{-(2k\pi + \theta_q) (s+it)}v)$$ is a trivial strip $u_\gamma$ over a Reeb chord $\gamma$ given by $$\gamma(t)= e^{it\theta_q}.s_1(q).$$ Thus, from the exponential convergence in (\ref{Nconv}) we get the following estimate with respect to the metric $g$:
\begin{equation}
    d_g (u(s,t) , u_\gamma(s,t)) < C e^{- (\delta + 2k\pi + \theta_q)s}\quad \forall s>s_0.
\end{equation}

We now show that the exponential convergence with respect to $g$ in the compact manifold $N$ would imply our required exponential convergence in cylindrical metric, as proposed in (\ref{expConv}). Choose a trivialization $U\times \CP^1$ of $N$.  We claim that we can arrange $g$, such that it agrees with $g_{cyl}$ on the contact distribution  $\ker \lambda$. Indeed, since $\pi_Y$ induces an isomorphism of $T_{\pi_Y(p)}Y$ with the hyperplane $(\ker \lambda)|_p\hookrightarrow T_pZ$, we can assume $g$ agrees with $g_{cyl}$. From our choices of $g,g_{cyl}$ we have that in the trivialization $U \times \CP^1$,  $\pi_U \of u$ converges exponentially fast to $\pi_U \of u_\gamma$ in $O(e^{-\delta s})$. It remains to check is the convergence of $\pi_{\CP^1} (u)$ to  $\pi_{\CP^1}(u_\gamma)$ in $O(e^{-\delta s})$.

We fix parameterization, which we use in the upcoming estimates. We use the radial parametrization $[1,\infty) \times S^1 \to \C^* $ , $(r,\theta) \mapsto (re^{i\theta})$. In this parametrization, the restrictions of the metric to $[R,\infty) \times S^1$ for  large $R>0$ direction are given as follows : \begin{align*}
g_{cyl}|_{\C^*} &=  \frac{1}{r^2} dr^2 + d\theta^2\\g|_{\CP^1} &=  \frac{1}{r^4} dr^2 +\frac{1}{r^2} d\theta^2.    
\end{align*} The $\frac{1}{r^2}$ coefficient shows up in the formula of $g_{cyl}$ since the radial to cylindrical parameterization is given by $(r,\theta) \mapsto (\log r , \theta)$. Also note that the metric $g|_{\CP^1}$ is isomorphic to the Euclidean metric on the punctured  disc via the map $(r,\theta) \mapsto (\frac{1}{r},\theta)$. Denote the metric induced from these by $d_g$,$d_{cyl}$.

\begin{lemma} There exists $s_0,C >0$ such that
$$d_{cyl} (\pi_{\CP^1} (u(s,t)),\pi_{\CP^1}(u_\gamma(s,t))) < Ce^{-\delta s} \quad  \forall s>s_0.$$
\end{lemma}

\begin{proof}
In the proof, we drop the arguments $(s,t)$ from the maps to avoid notational clutter. Denote $\pi_{\CP^1} \of u$  and $\pi_{\CP^1} \of u_\gamma$ in radial coordinates for $ \C^*$ by  $(u^r,u^\theta)$ , $(u^r_\gamma,u^\theta_\gamma)$ respectively. From $d_{g} (u,u_y)=O(e^{- (\delta + 2k\pi + \theta_q)s})$ we get for $r\in [R,\infty)$, $$d_g((u^r,u^\theta),(u^r_\gamma,u^\theta_\gamma)) = d^{\D^2}_{polar}((\frac{1}{u^r},u^\theta),(\frac{1}{u_\gamma^r},u_\gamma^\theta)) = O(e^{- (\delta + 2k\pi + \theta_q)s}) .$$
In the last equation $d^{\D^2}_{polar}$ is the Euclidean distance between two points in polar coordinates.  We now  prove the estimate in $\R$ direction in $\R \times S^1$.
\begin{align*}
    &|\frac{1}{u^r} - \frac{1}{u^r_\gamma}| = O(e^{- (\delta + 2k\pi + \theta_q)s})\\
     \implies &| \frac{u_\gamma^r}{u^r} - 1|= O(e^{- \delta s})\\
     \implies & |\log u^r  - \log u^r_\gamma |=| \log \frac{u_\gamma^r}{u^r} | = O(e^{- \delta s})
\end{align*}
Here, the first equality follows from $|r_1-r_2| \leq d^{\D^2}_{polar}((r_1,\theta_1),(r_2,\theta_2))$ and the last equality follows from the Taylor expansion of $\log(1+x)$ and the preceding line. 

Now we prove the estimate in the $S^1$ direction. We have that $2r_1r_2(1-\cos (\theta_1 - \theta_2)) \leq (d^{\D^2}_{polar}((r_1,\theta_1),(r_2,\theta_2)))^2$. Thus, we see that, 
\begin{align*}
    & 2\frac{1}{u^ru^r_\gamma}(1-\cos (u^\theta - u^\theta_\gamma)) = O(e^{- 2(\delta + 2k\pi + \theta_q)s})\\
    \implies & 1-\cos (u^\theta - u^\theta_\gamma)= O(e^{- 2\delta s})\\
    \implies & |u^\theta - u^\theta_\gamma|=O(e^{- \delta s}).
\end{align*}
In the above equations, we use that $|\log\frac{u^r}{u^r_\gamma}| =  O(e^{- \delta s})  $ implies that $u^r(s,t)= O(e^{(2k\pi + \theta_q)s})$ and that the Taylor series of $\cos x $.

\noindent Recall that the radial-to-cylindrical parametrization is given by $(r,\theta) \mapsto (\log r , \theta)$. We see that,
$$(|\log u^r -\log u^r_\gamma|^2 + |u^\theta- u^\theta_\gamma)|^2)^{\frac{1}{2}}=  d_{cyl}((u^r,u^\theta),(u^r_\gamma,u^\theta_\gamma)) = O(e^{- \delta s}).$$

\noindent Thus, we have $$d_{cyl} (\pi_{\CP^1} (u(s,t)),\pi_{\CP^1}(u_\gamma(s,t))) < Ce^{-\delta s} \quad  \forall s>s_0.$$

\end{proof}

\noindent Combining the exponential convergence of $\pi_{\CP^1} \of u$ and $\pi_Y \of u$ we get our required exponential convergence, $$d_{cyl}(u(s,t),u_\gamma(s,t)) < Ce^{-\delta s} \quad \forall s\geq s_0.$$

\end{proof}

\begin{remark}
\label{remark:asymptoticLagrangian}
We can prove an exponential convergence for strips with boundary on $L_a$ in $\R \times Z$ where $L_a$ projects to the Lagrangian $L_Y$ in $Y$ and $L_a$ is asymptotically close to the cylindrical Lagrangian $\R \times \Lambda$ so that closure of $L_a$ in $N$ is a Lagrangian with clean self intersection. We encounter such a boundary condition in the case of $T_{\R + i\e}$.  Steps 1,2,3  and 5 follow verbatim for this case. From a direct computation of eigenvalues and eigenfunctions of $A_{q_\infty}$ in this case we get the necessary counterpart of equation (\ref{Nconv}) of Step 4.  
\end{remark}

\subsection{The SFT compactness theorem}

We state the version of SFT compactness for discs in our case of neck-stretching here. We avoid adding marked points and perturbation of domain just for sake of notational clarity, the theorem still holds true with them.

\begin{Theorem} 
\label{theorem:SFTCompactness}
If $u_n : \D^2 \to (W^n,J^n,L^n) $ is a sequence of $J^n$ holomorphic disc to the $n-$stretched manifolds, with boundary on $L^n$, and with bounded energy, i.e. $E(u) < E< \infty$, then there is a subsequence of $u_n$ that converges to a stable broken disc.
\end{Theorem}

This theorem is mentioned in Section 11.3 in \cite{bhewz} for the case of non-degenerate Reeb chords. Our case is a Morse-Bott version for Reeb chords. The proof in this case is the same as the usual SFT compactness theorem apart from some analytic changes in the lemmas in section 5 in \cite{bhewz}. We prove the necessary counterparts of the analysis in section 5 in the following lemmas.

\begin{lemma}[Bubbling]\label{lemma:bubbling}
There exists a $\hbar>0$ such that the following is true. Let $F_n =(a_n,f_n):\D_+^2 \to (\R \times  Z, \R \times \Lambda)$ be a sequence of holomorphic maps from the upper half disc with Lagrangian boundary such that $E(F_n) < C$ for some constant $C$ and $a_n(0) = 0 $. Then there is a sequence of points $y_n \in \D^2_+$ converging to 0 and  sequences of real numbers $R_n,c_n$ such that the rescaled maps $F^0_n(z)=(y_n + \frac{z}{R_n})$ are functions on the half-disk of radius $R_n$ and have a convergent subsequence to $F^0$ in $C^\infty_{loc}$ such that $E_{Hor}>\hbar$
\end{lemma} 

\begin{proof}
If $|\nabla F_n(0)| \to \infty$ as $n \to \infty$.  As usual with bubbling analysis, we apply Hofer's lemma to obtain $y_n \in \D^2_+ $ and $\e_n \in \R$ such that $y_n\to 0$, $\sup |\nabla F_n (z) | \leq 2 | \nabla F_n(y) | $ where the supremum is taken on $|z-y_n| \leq \e_n$ and $|\nabla F_n(y_n)|\e_n \to \infty$. Now consider the rescaled maps $F^0_n(z) = F_n(y_n + \frac{z}{|\nabla F_n(y_n)|})$. This rescaled map is defined for $|z| \leq \e_n|\nabla F_n(y_n)|$ and has upper bound on the gradient by 2, thus by Arzela-Ascoli we get a convergent subsequence that converges to a non-constant function $F^0$ on the upper half plane $\mathbb{H}$ or the complex plane $\C$ depending whether $y_n$ was on the boundary or in the interior. Clearly we have $E(F^0) \leq C$. If the image of $F^0$ is bounded, $F^0$  maps to a holomorphic sphere or a disc with Lagrangian boundary and from standard energy quantization result we have $E_{Hor}(F^0) > \hbar > 0$ where $\hbar$ depends on $Z$. for some fixed  or else it maps asymptotically to a Reeb chord or Reeb orbit. If $F^0$ has unbounded image, from Lemma \ref{lemma:asym},  it has Reeb-orbit/chord asymptotes . Since we are in the contact case, we can directly compute $E_{Hor}$ and see that it is bounded below by the smallest length of the Reeb chord, thus $\pi/n$ for our case. Thus choosing a smaller $\hbar$ we can ensure that $E_{Hor}(F^0) > \hbar$ where $\hbar$ depends only on $Z$.
\end{proof}

\begin{lemma}[Asymptotics of strips]\label{lemma:asym}
Let F=(a,f) : $\R_+ \times \interval \to (\R \times Z, \R \times \Lambda)$ be a $J$ holomorphic map of finite energy $E_H(u)< \infty$, assume that image of F is unbounded, then there is a number $T$ and a Reeb chord $\gamma$  of length $T$ such that $$\lim_{s\to \infty } f(s,t) = \gamma (tT) \quad and \quad \lim_{s\to\infty} \frac{a(s,t)}{s} = T \quad in \quad C^{\infty}\interval.$$
\end{lemma}

\noindent This lemma is a direct consequence of our asymptotic convergence proved in Theorem \ref{thm:exponential conv}.

\begin{lemma}[Strip with long neck and small energy]\label{lemma:energytodistance}
Given $E_0,\e>0$ there exists constants $\sigma,c>0$ such that for every $R>c$ and for every $J$ holomorphic function $F=(a,f) : [-R,R] \times \interval \to (\R \times Z, \R \times \Lambda)$ satisfying the inequalities, 
$$E_{Hor} (F) \leq \sigma \quad and \quad E(F)\leq E_0,$$ 
we have 
$$f(s,t) \in B_\e(f(0,t))$$ for all $s\in [-R+c,R-c] , t\in \interval.$

\end{lemma}

\begin{proof}
Assume the contrary, so we have an $\e_0$ and $R_n>c_n$ such that $c_n \to \infty$ and $s_n \in [-k_n,k_n]$ where $k_n=R_n - c_n$ such that $E_{Hor} < \frac{1}{n}$ and $d(f_n(s_n,t),f_n(0)) > \e_0$. If we select $\sigma$ small enough such that any holomorphic disc $v$ in $Y = Z/S^1 $ with boundary in the Lagrangian $L = \Lambda/S^1$  has $E_{Y}(v) > \sigma$ then  for $n$ large, we see that $ E_{Hor} (F_n) = E_{Y} (\pi_{Y} \of F_n) \leq  \frac{1}{n}$ .  From  the bubbling quantization lemma \ref{lemma:bubbling} we know that there must be a uniform gradient bound or else $E_{Hor}$ has to be at least $\hbar$. Thus, by Arzela-Ascoli we can take a convergent subsequence $F_n \to F^0$ in $C^\infty_{loc} (\R\times \interval)$. Clearly from considering the horizontal energy we see that $F^0$ is a strip whose image lies in a trivial cylinder over a point $p \in L$ with the same Reeb chord as its positive and negative asymptote. Thus, we can consider $F^0$ to be a map to $\C^*$ with boundary on $\R^+$ and $e^{k\alpha}\R^+$. From exponential convergence to Reeb orbits at the asymptotes we see that $F^0$ has to be a trivial Reeb strip i.e. $f^0(s,t) = \gamma (t)$ where $\gamma$ is the Reeb chord at the asymptote. From convergence on compact sets, we have that $f_n(0,t) \to \gamma(t)$. Now consider the sequence of translated maps $\bar{F}_n(s,t) = F_n(s-s_n,t_n)$, possibly taking a further subsequence we get a $C^\infty_{loc}$ convergent subsequence to $\bar{F}^0$. Potentially $\bar{F}^0$ can be a trivial strip over a different Reeb strip. It is classically known (see Lemma A.6 in \cite{FraunThes} and Lemma 3.8 (ii) in \cite{schm}) that long strips of small energy have diameter bounded by the energy. Thus, from the horizontal energy bound  we get,    $\pi_{Y} \of F_n \to p$ as $n \to \infty$ since $\pi_{Y} \of F_n(0,t) \to p$, thus $\bar{F}^0 = F^0$, thus $\bar{f}_n(0,t)=f_n(s_n,t) \to \gamma (t)$. One way to see this is to apply Lemma 3.8 (ii), since $|s_n| < k_n$, $|s_n| <  R_n - \frac{c_n}{2}$. As $c_n \to \infty$, we have for large n, $d(\pi_{Y} \of F_n(s_n,t), \pi_{Y} \of F_n(0,t)) < ce^{\mu \frac{c_n}{2}}$.   Thus $d(f_n(s_n,t),f_n(0,t))$ tends  to $ 0 $ in  limit, which contradicts our assumption. 

\end{proof}

Once one has the last three Lemmas we proved above, the proof for SFT compactness from \cite{bhewz} works in our setting. We give a sketch of the proof in \cite{bhewz} for continuity.

\begin{sketch}{\ref{theorem:SFTCompactness}}
The proof is broken into four steps.
\begin{itemize}
    \item Step 1: Gradient Bound
    \item Step 2: Convergence from gradient bound away from nodes
    \item Step 3: Convergence near the nodes aka `thin-part'
    \item Step 4: Labelling of the components in the nodal disk
\end{itemize}

\vskip 3pt
\textbf{Step 1: Gradient Bound}

The first step in showing SFT Compactness is to pick a subsequence so that the domain converges, for sake of clarity lets assume that our domain is constant for all $n$ and have no marked points. Then Lemma 10.7 in \cite{bhewz} and Proposition 3.7 in \cite{Abbasbook} prove that after removing a finite number of points from, $\D^2$ we get a uniform gradient bound for all $n$. 

\begin{lemma}[Gradient lemma]
There is a natural number $K(E)$ and a sequence $M_n$ of $2K(E)$ marked points on the disc such that after removing these marked points from $\D^2$, we have a uniform bound as follows $$|\nabla u_{i_n}(z)| \leq \frac{C}{\rho(z)} \quad \forall z\in \D^2_n,$$where $\D^2_n = \D^2 \setminus M_n$ and $\rho$ is the radius of injectivity with respect to the standard hyperbolic metric on $\D^2_n$ and $u_{i_n}$ is a subsequence of $u_n$. 
\end{lemma}

\begin{prooflemma*}
Assume $x_n \in \D^2$ such that $|\rho(x_n) \nabla u_n | \to \infty$. We can assume after taking a subsequence (and still indexing it by n) that either of the following conditions hold :
\begin{enumerate}
   
    \item $u_n(x_n) \in W^n \setminus ( (-n,n)\times Z) $,
    \item $u_n(x_n) \in  (-n/2,n/2) \times Z$.
\end{enumerate}

For case 1, after rescaling our domains and using usual the bubbling lemma we can add two marked points $x_n,y_n$, $d(x_n,y_n) \to 0$ such that the bubble captures horizontal energy at least $\hbar$.  For case 2, after translating to ensure $u_n(x_n) \in \{0\}\times Z$ and rescaling the domain and lemma \ref{lemma:bubbling}, we see that there is a bubbling of a possibly punctured holomorphic disk or sphere which captures at least $\hbar$ horizontal energy.  Thus, in either can choose  $x_n,y_n$ , $d(x_n,y_n) \to 0$ such that the $E_{Hor}(u_n) < E - \hbar$ on $D^2\setminus\{x_n,y_n\}$. From positivity of horizontal energy, we see that if we inductively keep adding marked points, our process will end in at most $\frac{E}{\hbar}$ steps. 

\end{prooflemma*}

\vskip 3pt
\textbf{Step 2: Convergence from gradient bound away from nodes}

After putting in the marked points to get gradient bounds, take a further subsequence of $u_n$ such that the domain with the marked points $M_n$ converge to a nodal disk with marked points. From the bubbling lemma (Lemma \ref{lemma:bubbling}), we have that $u_n$ converges on the components containing the markings $M_n$.  The gradient bound on the components not containing $M_n$ lets us use Arzela-Ascoli to get a convergent subsequence of $u_n$ that converges on every piece of the nodal disc.
\vskip 3pt
\textbf{Step 3: Convergence near the nodes aka `thin-part'}

Near the nodes in components with marked points $M_n$, we know that the maps $\mathfrak{u}$ are asymptotic to Reeb chords or points. It might happen that the asymptotic limit at the nodes might not match. Let's fix a node and pick a sequence of  parameterization on a neighborhood of  $\D^2$ which degenerate to that node, call the neighborhoods $V_n$. We can arrange the parameterization such that 

\begin{itemize}
    \item $V_n$ becomes a strip or a cylinder as $n \to \infty$ ,
    \item $u_n|_{V_n}$ has bounded gradient, thus bubbling can't occur ,
    \item $u_n (s,t)  \to $ Reeb chords/strips $\gamma^{\pm}$ as $s\to \pm\infty$  .
\end{itemize}

There are two cases, either $E_{hor}(u_n|_{V_n}) \to 0$ or $E_{hor}(u_n|_{V_n}) > \delta > 0$ for a fixed lower bound $\delta$. For the first case,  from Lemma \ref{lemma:energytodistance} we obtain that if there is no concentration of energy near the node, the asymptotes match i.e. $\gamma^+ = \gamma^-$. Thus, the only case left is if there is some energy concentrating near the nodes. in such a case one adds more marked points which results in formation of a new spherical/disc component which captures at least $\hbar$ horizontal energy, and the new nodes which are formed are treated inductively. The new nodes will eventually have no energy concentration. Thus, finally, we get a nodal disc with matching conditions on the nodes. The limit we obtain finally is denoted by $(\mathfrak{u,D})$

\vskip 3pt
\textbf{Step 4: Labelling of the components in the nodal disk}
Pick points $x_n^C$ for every component $C$ of $\mathfrak{ D}$. Label the components $C$ such that if  $u_n(x_n)\in W^n\setminus(-n,n)\times Z $ then label it $in$ or $out$ based on which component it lies in. If $u_n(x_n) \in (-n,n)\times Z $, order the components such that $C_i \leq C_j$ if $\text{limsup } (pr_\R (u_n(x^j_n)) - pr_\R (u_n(x^i_n))) \to \infty $. Set $C_i \sim C_j$ if $C_i\leq C_j$ and $C_j\leq C_i$.  Label the collection of the minimal components as $\D_1$, then remove them and label the collection of the new minimal components as $\D_2$ and keep going. In the case of such labeling, it might happen that across a node, the label just by more than 1, in such a case add trivial sphere/discs between and label them such that the difference across nodes is at most 1.

\end{sketch}

\subsection{ Compactness for strips bounding pair of Lagrangians under stretching}
    
We now return to our specific situation where we want to apply the neck stretching construction. Since the differential in the Floer complex $CF(L,K)$ is based on the count of  holomorphic strips, we will need a version  SFT compactness of strips under neck stretching.  Recall that a holomorphic strip in $M$ bounding $L,K$ from $x_+$ to $x_-$ is a holomorphic map $u : \R \times \interval \to M$ with boundary conditions $u(\R \times \{1\}) \sub L$ and $u(\R \times \{ 0\}) \sub K$ and asymptotic conditions $\lim_{s\to \pm \infty} u(s,t) = x_\pm $ where $x_\pm \in L \cap K$.  If $L$ is a locally mutable Lagrangian and $K$ is another Lagrangian that intersects $L$ transversally and lies outside the mutation neighborhood $B$  ie. $K \sub M \setminus B$, we can do a neck stretch of $M$ along the contact hypersurface $Z= \del B$ since near a tubular neighborhood of $Z$, the Lagrangian $L$ is of the form $(-\e ,\e) \times T_\pm$ (see Equation ( \ref{def:tplusminus}) for definition of $T_\pm$). Since $K$ is outside the mutation neighborhood, after neck stretching, $K$ naturally lies in the outside piece $\widehat{M}_{out}$ of the broken manifold $\mathbb{M}$. The analog of a broken disc in the setting of a strip is called a \textit{broken strip}. 

\begin{defn}[Broken strip]
A broken strip $\mathfrak({u,D})$ from $x_+$ to $x_-$ is a broken disc where one of the components of $\mathbb{D}_{out}$ is a strip $\R \times \interval$ with puncture nodes $P$ on the boundary $\R \times \{1\}$ and $u_{out}$ satisfies the boundary condition $u_{out}(\R \times \{1\} \setminus P) \sub \widehat{L}_{out}$ and $u_{out}(\R \times \{0\}) \sub K$ and asymptotic conditions $\lim_{s\to \pm \infty}u_{out}(s,t) = x_\pm$.
\end{defn}

We can use an exactly similar technique as in the proof of Theorem \ref{theorem:SFTCompactness} to get the following compactness result about holomorphic strips. 
\begin{Theorem}[SFT compactness for strips]
If $u_n : \R \times \interval \to (M^n,J^n,L^n,K) $ is a sequence of $J^n$ holomorphic strips from $x_+$ to $x_-$ in the $n-$stretched manifolds, with boundary condition on $L^n$ and $K$, and with bounded energy, i.e. $E(u) < E< \infty$, then there is a subsequence of $u_n$ that converges to a broken strip.
\end{Theorem}

%% file: pics/SFTsetup.pdf_tex
\begingroup%
  \makeatletter%
  \providecommand\color[2][]{%
    \errmessage{(Inkscape) Color is used for the text in Inkscape, but the package 'color.sty' is not loaded}%
    \renewcommand\color[2][]{}%
  }%
  \providecommand\transparent[1]{%
    \errmessage{(Inkscape) Transparency is used (non-zero) for the text in Inkscape, but the package 'transparent.sty' is not loaded}%
    \renewcommand\transparent[1]{}%
  }%
  \providecommand\rotatebox[2]{#2}%
  \newcommand*\fsize{\dimexpr\f@size pt\relax}%
  \newcommand*\lineheight[1]{\fontsize{\fsize}{#1\fsize}\selectfont}%
  \ifx\svgwidth\undefined%
    \setlength{\unitlength}{150bp}%
    \ifx\svgscale\undefined%
      \relax%
    \else%
      \setlength{\unitlength}{\unitlength * \real{\svgscale}}%
    \fi%
  \else%
    \setlength{\unitlength}{\svgwidth}%
  \fi%
  \global\let\svgwidth\undefined%
  \global\let\svgscale\undefined%
  \makeatother%
  \begin{picture}(1,1)%
    \lineheight{1}%
    \setlength\tabcolsep{0pt}%
    \put(0,0){\includegraphics[width=\unitlength,page=1]{SFTsetup.pdf}}%
    \put(0.33039032,0.86686804){\makebox(0,0)[lt]{\lineheight{1.25}\smash{\begin{tabular}[t]{l}$W$\end{tabular}}}}%
    \put(0.31505575,0.74000931){\makebox(0,0)[lt]{\lineheight{1.25}\smash{\begin{tabular}[t]{l}$L$\end{tabular}}}}%
    \put(0.60432156,0.5228857){\makebox(0,0)[lt]{\lineheight{1.25}\smash{\begin{tabular}[t]{l}$(-\e,\e) \times Z$\end{tabular}}}}%
    \put(-0.03101767,0.52009758){\makebox(0,0)[lt]{\lineheight{1.25}\smash{\begin{tabular}[t]{l}$(-\e,\e) \times \Lambda$\end{tabular}}}}%
    \put(0,0){\includegraphics[width=\unitlength,page=2]{SFTsetup.pdf}}%
  \end{picture}%
\endgroup%

%% file: pics/breaking.pdf_tex
\begingroup%
  \makeatletter%
  \providecommand\color[2][]{%
    \errmessage{(Inkscape) Color is used for the text in Inkscape, but the package 'color.sty' is not loaded}%
    \renewcommand\color[2][]{}%
  }%
  \providecommand\transparent[1]{%
    \errmessage{(Inkscape) Transparency is used (non-zero) for the text in Inkscape, but the package 'transparent.sty' is not loaded}%
    \renewcommand\transparent[1]{}%
  }%
  \providecommand\rotatebox[2]{#2}%
  \newcommand*\fsize{\dimexpr\f@size pt\relax}%
  \newcommand*\lineheight[1]{\fontsize{\fsize}{#1\fsize}\selectfont}%
  \ifx\svgwidth\undefined%
    \setlength{\unitlength}{225bp}%
    \ifx\svgscale\undefined%
      \relax%
    \else%
      \setlength{\unitlength}{\unitlength * \real{\svgscale}}%
    \fi%
  \else%
    \setlength{\unitlength}{\svgwidth}%
  \fi%
  \global\let\svgwidth\undefined%
  \global\let\svgscale\undefined%
  \makeatother%
  \begin{picture}(1,1.16666667)%
    \lineheight{1}%
    \setlength\tabcolsep{0pt}%
    \put(0,0){\includegraphics[width=\unitlength,page=1]{breaking.pdf}}%
    \put(0.74761358,0.59087438){\makebox(0,0)[lt]{\lineheight{1.25}\smash{\begin{tabular}[t]{l}\large{$W^\tau$}\end{tabular}}}}%
    \put(0,0){\includegraphics[width=\unitlength,page=2]{breaking.pdf}}%
    \put(0.42764414,0.81195109){\makebox(0,0)[lt]{\lineheight{1.25}\smash{\begin{tabular}[t]{l}$Z \times (a,b)$\end{tabular}}}}%
    \put(0.68957616,0.98644521){\makebox(0,0)[lt]{\lineheight{1.25}\smash{\begin{tabular}[t]{l}\large{$W$}\end{tabular}}}}%
    \put(0,0){\includegraphics[width=\unitlength,page=3]{breaking.pdf}}%
    \put(0.37600586,0.47385061){\makebox(0,0)[lt]{\lineheight{1.25}\smash{\begin{tabular}[t]{l}$Z \times (-\tau,\tau)$\end{tabular}}}}%
    \put(0.03955824,0.13931366){\makebox(0,0)[lt]{\lineheight{1.25}\smash{\begin{tabular}[t]{l}\large{$\widehat W_{in}$}\end{tabular}}}}%
    \put(0.74199671,0.15609151){\makebox(0,0)[lt]{\lineheight{1.25}\smash{\begin{tabular}[t]{l}\large{$\widehat W _{out}$}\end{tabular}}}}%
    \put(0.43549962,0.19652374){\makebox(0,0)[lt]{\lineheight{1.25}\smash{\begin{tabular}[t]{l}$Z \times \R$\end{tabular}}}}%
  \end{picture}%
\endgroup%

%% file: pics/brokendisc.pdf_tex
\begingroup%
  \makeatletter%
  \providecommand\color[2][]{%
    \errmessage{(Inkscape) Color is used for the text in Inkscape, but the package 'color.sty' is not loaded}%
    \renewcommand\color[2][]{}%
  }%
  \providecommand\transparent[1]{%
    \errmessage{(Inkscape) Transparency is used (non-zero) for the text in Inkscape, but the package 'transparent.sty' is not loaded}%
    \renewcommand\transparent[1]{}%
  }%
  \providecommand\rotatebox[2]{#2}%
  \newcommand*\fsize{\dimexpr\f@size pt\relax}%
  \newcommand*\lineheight[1]{\fontsize{\fsize}{#1\fsize}\selectfont}%
  \ifx\svgwidth\undefined%
    \setlength{\unitlength}{300bp}%
    \ifx\svgscale\undefined%
      \relax%
    \else%
      \setlength{\unitlength}{\unitlength * \real{\svgscale}}%
    \fi%
  \else%
    \setlength{\unitlength}{\svgwidth}%
  \fi%
  \global\let\svgwidth\undefined%
  \global\let\svgscale\undefined%
  \makeatother%
  \begin{picture}(1,0.5)%
    \lineheight{1}%
    \setlength\tabcolsep{0pt}%
    \put(0,0){\includegraphics[width=\unitlength,page=1]{brokendisc.pdf}}%
    \put(0.22211895,0.41728624){\makebox(0,0)[lt]{\lineheight{1.25}\smash{\begin{tabular}[t]{l}$in$\end{tabular}}}}%
    \put(0.33643124,0.41542751){\makebox(0,0)[lt]{\lineheight{1.25}\smash{\begin{tabular}[t]{l}$in$\end{tabular}}}}%
    \put(0.42750931,0.41356877){\makebox(0,0)[lt]{\lineheight{1.25}\smash{\begin{tabular}[t]{l}$1$\end{tabular}}}}%
    \put(0.52509293,0.41449814){\makebox(0,0)[lt]{\lineheight{1.25}\smash{\begin{tabular}[t]{l}$out$\end{tabular}}}}%
    \put(0.55030205,0.46805298){\makebox(0,0)[lt]{\lineheight{1.25}\smash{\begin{tabular}[t]{l}$out$\end{tabular}}}}%
    \put(0.34107807,0.08643124){\makebox(0,0)[lt]{\lineheight{1.25}\smash{\begin{tabular}[t]{l}$\mathbb{W}[1]$\end{tabular}}}}%
  \end{picture}%
\endgroup%

%% file: pics/cylfib.pdf_tex
\begingroup%
  \makeatletter%
  \providecommand\color[2][]{%
    \errmessage{(Inkscape) Color is used for the text in Inkscape, but the package 'color.sty' is not loaded}%
    \renewcommand\color[2][]{}%
  }%
  \providecommand\transparent[1]{%
    \errmessage{(Inkscape) Transparency is used (non-zero) for the text in Inkscape, but the package 'transparent.sty' is not loaded}%
    \renewcommand\transparent[1]{}%
  }%
  \providecommand\rotatebox[2]{#2}%
  \newcommand*\fsize{\dimexpr\f@size pt\relax}%
  \newcommand*\lineheight[1]{\fontsize{\fsize}{#1\fsize}\selectfont}%
  \ifx\svgwidth\undefined%
    \setlength{\unitlength}{225bp}%
    \ifx\svgscale\undefined%
      \relax%
    \else%
      \setlength{\unitlength}{\unitlength * \real{\svgscale}}%
    \fi%
  \else%
    \setlength{\unitlength}{\svgwidth}%
  \fi%
  \global\let\svgwidth\undefined%
  \global\let\svgscale\undefined%
  \makeatother%
  \begin{picture}(1,0.66666667)%
    \lineheight{1}%
    \setlength\tabcolsep{0pt}%
    \put(0,0){\includegraphics[width=\unitlength,page=1]{cylfib.pdf}}%
    \put(0.61547399,0.2779585){\makebox(0,0)[lt]{\lineheight{1.25}\smash{\begin{tabular}[t]{l}$\R \times Z$\end{tabular}}}}%
    \put(0.44885543,0.56972349){\makebox(0,0)[lt]{\lineheight{1.25}\smash{\begin{tabular}[t]{l}$\C^*\text{ fiber}$ \end{tabular}}}}%
    \put(0.10875429,0.42485908){\makebox(0,0)[lt]{\lineheight{1.25}\smash{\begin{tabular}[t]{l}$\{r\}\times Z$\end{tabular}}}}%
    \put(0.41589221,0.113228){\makebox(0,0)[lt]{\lineheight{1.25}\smash{\begin{tabular}[t]{l}$Y$\end{tabular}}}}%
  \end{picture}%
\endgroup%

%% file: pics/compcyl.pdf_tex
\begingroup%
  \makeatletter%
  \providecommand\color[2][]{%
    \errmessage{(Inkscape) Color is used for the text in Inkscape, but the package 'color.sty' is not loaded}%
    \renewcommand\color[2][]{}%
  }%
  \providecommand\transparent[1]{%
    \errmessage{(Inkscape) Transparency is used (non-zero) for the text in Inkscape, but the package 'transparent.sty' is not loaded}%
    \renewcommand\transparent[1]{}%
  }%
  \providecommand\rotatebox[2]{#2}%
  \newcommand*\fsize{\dimexpr\f@size pt\relax}%
  \newcommand*\lineheight[1]{\fontsize{\fsize}{#1\fsize}\selectfont}%
  \ifx\svgwidth\undefined%
    \setlength{\unitlength}{225bp}%
    \ifx\svgscale\undefined%
      \relax%
    \else%
      \setlength{\unitlength}{\unitlength * \real{\svgscale}}%
    \fi%
  \else%
    \setlength{\unitlength}{\svgwidth}%
  \fi%
  \global\let\svgwidth\undefined%
  \global\let\svgscale\undefined%
  \makeatother%
  \begin{picture}(1,0.66666667)%
    \lineheight{1}%
    \setlength\tabcolsep{0pt}%
    \put(0,0){\includegraphics[width=\unitlength,page=1]{compcyl.pdf}}%
    \put(0.1167234,0.32379317){\makebox(0,0)[lt]{\lineheight{1.25}\smash{\begin{tabular}[t]{l}$\R \times Z$\end{tabular}}}}%
    \put(0.29832715,0.5579306){\makebox(0,0)[lt]{\lineheight{1.25}\smash{\begin{tabular}[t]{l}$Y_\infty$\end{tabular}}}}%
    \put(0.57527883,0.13506816){\makebox(0,0)[lt]{\lineheight{1.25}\smash{\begin{tabular}[t]{l}$Y_0$\end{tabular}}}}%
    \put(0.35966543,0.04306071){\makebox(0,0)[lt]{\lineheight{1.25}\smash{\begin{tabular}[t]{l}$N$\end{tabular}}}}%
  \end{picture}%
\endgroup%

%% file: pics/lagangle.pdf_tex
\begingroup%
  \makeatletter%
  \providecommand\color[2][]{%
    \errmessage{(Inkscape) Color is used for the text in Inkscape, but the package 'color.sty' is not loaded}%
    \renewcommand\color[2][]{}%
  }%
  \providecommand\transparent[1]{%
    \errmessage{(Inkscape) Transparency is used (non-zero) for the text in Inkscape, but the package 'transparent.sty' is not loaded}%
    \renewcommand\transparent[1]{}%
  }%
  \providecommand\rotatebox[2]{#2}%
  \newcommand*\fsize{\dimexpr\f@size pt\relax}%
  \newcommand*\lineheight[1]{\fontsize{\fsize}{#1\fsize}\selectfont}%
  \ifx\svgwidth\undefined%
    \setlength{\unitlength}{225bp}%
    \ifx\svgscale\undefined%
      \relax%
    \else%
      \setlength{\unitlength}{\unitlength * \real{\svgscale}}%
    \fi%
  \else%
    \setlength{\unitlength}{\svgwidth}%
  \fi%
  \global\let\svgwidth\undefined%
  \global\let\svgscale\undefined%
  \makeatother%
  \begin{picture}(1,0.66666667)%
    \lineheight{1}%
    \setlength\tabcolsep{0pt}%
    \put(0,0){\includegraphics[width=\unitlength,page=1]{lagangle.pdf}}%
    \put(0.54739777,0.3460347){\makebox(0,0)[lt]{\lineheight{1.25}\smash{\begin{tabular}[t]{l}$L_0$\end{tabular}}}}%
    \put(0.49488846,0.56908302){\makebox(0,0)[lt]{\lineheight{1.25}\smash{\begin{tabular}[t]{l}$L_1$\end{tabular}}}}%
    \put(0.35824689,0.50159164){\makebox(0,0)[lt]{\lineheight{1.25}\smash{\begin{tabular}[t]{l}$L_{Y_\infty}$\end{tabular}}}}%
    \put(0.15009293,0.44687113){\makebox(0,0)[lt]{\lineheight{1.25}\smash{\begin{tabular}[t]{l}$u(s,t)$\end{tabular}}}}%
    \put(0.58039032,0.1425031){\makebox(0,0)[lt]{\lineheight{1.25}\smash{\begin{tabular}[t]{l}$\theta_q$\end{tabular}}}}%
    \put(0.30181229,0.40470105){\makebox(0,0)[lt]{\lineheight{1.25}\smash{\begin{tabular}[t]{l}$q_\infty$\end{tabular}}}}%
    \put(0.38940516,0.20198263){\makebox(0,0)[lt]{\lineheight{1.25}\smash{\begin{tabular}[t]{l}$T_{q_\infty}L_0$\end{tabular}}}}%
    \put(0.7230483,0.04770757){\makebox(0,0)[lt]{\lineheight{1.25}\smash{\begin{tabular}[t]{l}$T_{q_\infty}L_1$\end{tabular}}}}%
  \end{picture}%
\endgroup%

%% file: Sections/modconstruct.tex
\section{Constructing Moduli Space of Broken Maps } \label{sec:modcons}

In this section, we describe the construction of moduli space of broken strips and broken discs as zero sets of a Fredholm map between Banach manifolds. From this section on wards, we will focus on the broken manifold obtained from performing a neck-stretching along the boundary of a mutation neighborhood of a manifold $M$.

We recall the notation from  Subsection \ref{subsec:sftsetup} and the spaces which show up in the broken manifold after stretching along the boundary of a mutation neighborhood. The contact hypersurface $Z$ is $S^{2n-1}$ and the Legendrian $\Lambda$ is diffeomorphic to a disjoint union of two $n-1$ dimensional tori.  The broken manifold (as defined in Subsection \ref{subsub:brokenstuff}) produced by neck-stretching across the contact sphere $S^{2n-1}$ is the collection,
$$\widehat M_{in}=\C^n,\;  \R \times S^{2n-1}  , \dots , \widehat M_{out} $$
The broken mutable Lagrangian is the collection,
$$\widehat L_{in}=T_\gamma,\; \R \times (T_+^{n-1} \sqcup T_-^{n-1})  , \dots , \widehat L_{out} $$
where the curve $\gamma$ is a cylindrical extension of the curve determining $L$ in the mutation neighborhood. The analysis for setting up moduli spaces of broken discs with boundary on Lagrangian that do not meet the neck has already been developed in the literature, see e.g. \cite{floerflip}. In the situation where the Lagrangian does not meet the neck, the puncture nodes in the broken maps occur only in the interior. In this subsection,  we will deal with the case of broken discs and strips where the puncture nodes occur only on the boundary and refer the reader to \cite{floerflip} for the case of the interior puncture nodes.

\subsection{Constructing Map spaces}

We describe a Banach manifold structure on a subspace of continuous maps from punctured discs to $\C^n$ or $\C^n\setminus \{0\}$ with boundary lying on the torus-segment Lagrangian. We aim to construct a Banach manifold such that any punctured holomorphic disc with boundary on the relevant torus segments. We construct the space of maps by starting with continuous maps that are trivial strips over Reeb chords near the punctures. Then we use geodesic flow exponentiation to define open sets around the maps with trivial strips.
\subsubsection{Continuous maps with strip-like ends}
We begin by describing the cylindrical structure on the domain near punctures.  Let $P$ be a finite set of points on the boundary of the unit disc, let $\Sigma_P$ denote the punctured disc with punctures at $P$ i.e. $\Sigma_P = \mathbb{D}\setminus P$. From now on, we assume that the set of punctures $P$ has $l$ elements and is equal to $\{p_i\}_{i=0}^l$. We label each puncture with either $+$ or $-$. We  fix holomorphic charts $\phi^\pm_p$ on $U_p$, a neighborhood of the puncture $p$, the sign $\pm$ is determined by  the label of the puncture $p$, such that for a large $M_p >0$ we have the following biholomorphisms $$\phi_{p,+}:U_p\setminus\{p\}\to [M_p,\infty) \times [0,1]$$
$$\phi_{p,-}:U_p\setminus\{p\}\to (-\infty, -M_p] \times [0,1],$$
and  $\phi_{p,\pm}$ maps the boundary of $\Sigma$ to $\R \times \{0,1\}$.We can compose this chart map with the biholomorphic map $z \mapsto e^{\mp 2\pi z}$ from the strip $\R \times \interval$ to the punctured upper half-space $\mathbb{H}\setminus\{0\}$ to get new charts $e^{\mp2\pi\phi^\pm_p}$, assume that these charts extend to charts from $U_p$ where $p \mapsto 0$. The neighborhoods $U_p$ under the coordinates $\phi^\pm_p$ are hereafter  called strip-like ends of $\Sigma$ .  We fix a volume form $d\Sigma$ on $\Sigma$ that is equal to the translation invariant volume form  $dsdt$ on the strip-like ends. 

We define a space of continuous maps that are equal to trivial strips near boundary punctures.  The set of continuous maps from $\Sigma$ to $\C^n \setminus \{0\} \iso \R \times S^{2n-1}$ that satisfies the Lagrangian boundary condition on $T_{\gamma}$ where $\gamma$ is cylindrical path and restrict to trivial Reeb-strips, $$u_c \left( s,t \right) =\left(k\frac{\pi}{n}s,c_{R}(kt) \right),$$
 on the  strip-like ends for non-zero integer $k$ and Reeb chord $c_{R}$ on $S^{2n-1}$ with boundary on the Legendrian tori is denoted as $\Maps^0(\Sigma_P, \C^n , T_\gamma)$ i.e. 
 
 \begin{align*}
     \Maps^0 (\Sigma_P, \C^n , T_\gamma) := \bigg\{ u \in C^0(\Sigma, \C^n, T_\gamma )\bigg| &\exists  K,m,R_u > 0 , c_R  \text{ a Reeb chord}, \\&\text{ such that } \\& u\of \phi_{p,\pm}\inv(s,t)=Ke^{\pm ms}c_{R}(\pm mt) \\&\text{ for  } |s|>R_u  \bigg\}.
 \end{align*}

We now define the weighted Sobolev space that we will use to define charts for our space of maps. Let $\beta \in C^\infty(\R)$ be a non-decreasing bump function that is  equal to $1$ on $[1,\infty)$ and equal to $0$ on $(-\infty,0]$. First we describe a function $\kappa^\delta_\Sigma$ that we use to define weighted Sobolev norm. 

\begin{align*}
    &\kappa^\delta_\Sigma : \Sigma \to \R \\
    &\kappa^\delta_\Sigma (z) = 0 \text{ for  } z \in \Sigma \setminus \cup_i U_{p_i} \\
    &\kappa^\delta_\Sigma (\phi_p \inv (s,t)) = \delta \beta( | s |- M_p) \text{ for } p  \in  P
\end{align*}

\noindent Let $g_{cyl}$ be the cylindrical metric  on $\C^n$ obtained from renormalizing  the standard Riemannian metric $g_0$ i.e. $g_{cyl} = g_0/r^2$. Let $u \in \Maps^0_{strip}(\Sigma_P, \C^n , T_\gamma)$, we define the $\delta$ weighted  Sobolev space of sections of $T\C^n$ over $u$  with boundary on $TT_\gamma$  and $k,p$ regularity with $kp>2$ using the cylindrical metric $g_{cyl}$ as follows,

\begin{align*}
    W^{k,p,\delta} (u^*T\C^n,\del  u^* TT_\gamma) = \bigg \{ \eta \in W_{loc}^{k,p} (u^*T\C^n)   \bigg |  \int (&|\eta|^p + |\nabla \eta |^p  + \\&\dots |\nabla^k \eta|^p)e^{p\kappa^\delta_\Sigma}  ) d\Sigma<\infty,\quad\\ & \eta(z) \in TT_\gamma \text{ for } z\in \del \Sigma  \bigg\}.
\end{align*}
where   $W^{k,p}_{loc}$ denotes the set of sections of  the bundle $u^*T\C^n$ that under any trivialization of the bundle $u^*T\C^n$ is in the Sobolev space $W^{k,p}_{loc}$ and  $\nabla$ refers to the weak derivative taken with respect to the standard connection on $T\C^n$. We define the $(k,p,\delta)$ norm on sections as follows $$\| \eta \|_{k,p,\delta}^p = \int (|\eta|^p + |\nabla \eta |^p  + \dots |\nabla^k \eta|^p)e^{p\kappa^\delta_\Sigma}  ) d\Sigma.$$

\subsubsection{Choosing domain dependent metric}
We will use geodesic exponentiation to obtain a sub-basis for the topology on the subspace of maps. We now choose a metric to exponentiate along to construct the sub-basis. See section 5.3 \cite{Mauthesis} for an expansive treatment of domain-dependent metric in the context of quilts. Let $RM(\C^n)$ be the space of metrics on $\C^n$. For a map $u \in \Maps^0_{strip}(\Sigma_P, \C^n , T_\gamma)$  choose $R>0$ such that $u(s,t)$ is cylindrical on all the strip-like ends and the image of $[R,\infty) \times \interval$ lies the cylindrical end $T_\gamma \setminus B_R(0)$ of the Lagrangian $T_\gamma$. Choose a metric $g_u$ on $\C^n$ that makes the Lagrangian $T_\gamma \cap B_{R+1}(0)$ totally geodesic. We can choose such a metric from corollary $3.2$ in \cite{abouzaidExotic} or by using Frauenfelder's metric, see Lemma 4.3.4 \cite{jholobook}. Let $g^u_\Sigma$ be a domain dependent choice of metrics that interpolates between the metric $g_u$   and standard cylindrical metric $g_{cyl} = g_0/r^2$ where $r$ is the radial coordinate as follows 

\begin{align*}
    &g^u_\Sigma : \Sigma \to RM(\C^n) \\
    &g^u_\Sigma(z) = g_u \text{ for  } z \in \Sigma \setminus \cup_i U_{p_i} \\
    &g^u_\Sigma(\phi_p \inv (s,t)) = (1-\beta(|s| -M_p))g_u +  \beta(|s| -M_p)g_{cyl} \text{ for  } p \in P
\end{align*}

\noindent For $\eta \in W^{k,p,\delta} (u^*T\C^n,\del  u^* TT_\gamma)$  we can define the map   $\exp _u (\eta)$ where we do point-wise geodesic exponentiation with respect to the domain dependent metric $g^u_\Sigma$  $$\exp _u (\eta) (z) = \exp_{u(z),g(z)}(\eta (z))$$
\noindent Notice that our choice of domain-dependent metric ensures that $\exp_u \eta $ is a map that preserves the Lagrangian boundary condition $ \exp_u \eta (\del \Sigma) \in T_\gamma$. This follows from the fact that on the cylindrical end of $T_\gamma$, the Legendrian tori $T_+,T_-$ are totally geodesic with the round metric on the sphere.

\subsubsection{Constructing the set of $W^{k,p,\delta}$ maps }
We now describe the set of maps on which we will give a Banach manifold structure. We define the map space $\Maps^{k,p}_\delta(\Sigma_P, \C^n , T_\gamma)$ from which we cut out the moduli space of parameterized holomorphic maps as follows

\begin{align*}
    \Maps^{k,p}_\delta (\Sigma_P, \C^n, T_\gamma) := 
\bigg\{ \exp_{u_0} (X) \bigg|\,u_0 \in &\Maps^0_{strip}(\Sigma_P, \C^n , T_\gamma),\,\,\\& X\in  W^{k,p,\delta}({u_0}^* T\C^n,\partial {u_0}^* TT_\gamma) \bigg\}.
\end{align*}

We show that any finite Hofer energy punctured holomorphic disc with correct boundary conditions lie in the space of maps we defined above. Let $u$ be a holomorphic map  on the domain $\Sigma$ to $\C^n$ with boundary lying on the torus segment $T_\gamma$ which has finite Hofer energy $E_H(u) < \infty$. We have exponential convergence (of decay rate at least $e^{-\pi s/n}$) to trivial Reeb strips on its strip-like ends from Theorem \ref{thm:exponential conv}. Thus, for some $u_0 \in \Maps^0_{strip}(\Sigma_P, \C^n , T_\gamma)$, we have $u=\exp_{u_0}(\eta)$ for some $\eta \in W^{k,p,\delta} (u^*T\C^n,\del  u^* TT_\gamma)$ where the geodesic exponentiation is done based on the domain dependent metric $g^{u_0}_\Sigma$  and $ \delta  < \frac{\pi}{n} $. This proves that any finite Hofer energy punctured disc $u$ can be written as a geodesic exponentiation of a $W^{k,p,\delta}$ vector field on a map $u^0 \in \Maps^0_{strip}(\Sigma_P, \C^n , T_\gamma) $.

\subsubsection{ Banach structure on $\Maps^{k,p}_\delta(\Sigma_P, \C^n , T_\gamma)$}

Now we will describe a Banach manifold structure on the space of maps by prescribing a generating set for an atlas. Let $V$ be the space generated by the  $n-1$ independent vector fields $X_i$ on $\C^n$  where  
$$X_k(z_1,\dots,z_n) = (0,\dots,0, iz_k,0,\dots,-i z_n).$$
Let $H$ be the one dimensional space generated by the radial vector field $X_{rad}$ on $\C^n$ where $\beta$ is the bump function we chose before and $R>0$ such that $T_\gamma \setminus B_R(0)$ is cylindrical.
$$X_{rad}(z)= \beta(|z|/R)z $$
\noindent We can check that all the sections of  $V \oplus H$ restrict to vector fields on $T_\gamma$ under the restriction of the sections to $T_\gamma$. Given a puncture $p$ and a map $u\in \Maps^{k,p}_\delta(\Sigma_P, \C^n , T_\gamma)$, for any $v\in V $, $h\in H$ we can define $v_{u,p},h_{u,p}$ to be sections on $u^*T\C^n$ which are equal to $v,h$ in the strip-like end at $p$ and extended to $\Sigma$  by a bump function so that $v_{u,p},h_{u,p} = 0$ on $\Sigma \setminus U_p$

We describe charts on the space of maps now. The Banach charts on $\Maps^{k,p}_\delta(\Sigma_P, \C^n , T_\gamma)$ at $u_0 \in \Maps^0_{strip}(\Sigma_P, \C^n , T_\gamma)$  are given by geodesic exponentiation 

\begin{equation}\label{eq:chartbanach}
    \exp_{u_0} :  W^{k,p,\delta}({u_0}^* T\C^n,\partial {u_0}^* TT_\gamma) \oplus V^{|P|} \oplus H^{|P|}  \to \Maps^{k,p}_\delta(\Sigma_P, \C^n , T_\gamma)
\end{equation}

\begin{equation*}
    \exp_{u_0} (w,v_1,\dots,v_l,h_1,\dots,h_l) =  \exp_{u_0}(w + \sum_{i=1}^l(v_{i,u_0,p_i} + h_{i,u_0,p_i}))
\end{equation*}

\noindent Since the choice of domain dependent metric differ only on compact sets, the transition functions
$\exp_{u_0} \of \exp_{u_1} \inv$ are maps between the corresponding Banach tangent spaces at $u_0,u_1 \in \Maps^0_{strip}(\Sigma_P, \C^n , T_\gamma)$. Following the traditional notation in treatments of $J-$holomorphic maps, we denote $\Maps^{k,p}_\delta(\Sigma_P, \C^n , T_\gamma)$ as $\mathcal{B}$ when we want to view it as a Banach manifold.

We describe a Banach bundle over the space of maps, and show that the space of holomorphic maps can be viewed as the zero locus of a certain section of the bundle. Define a Banach vector bundle  $\mathcal{E} $ over $\mathcal{B}$ as follows $$\mathcal{E} = \bigcup _{u\in \mathcal{B}} \,\; W^{k-1,p,\delta}(\overline{\operatorname{Hom}}(T\Sigma,u^*T\C^n))    $$
where $W^{k-1,p,\delta}(\overline{\operatorname{Hom}}(T\Sigma,u^*T\C^n))$ is the space of $(0,1)$ forms on $\Sigma$ with finite weighted  norm as before, 
$$\| \eta \|_{k-1,p,\delta}^p = \int (|\eta|^p + |\nabla \eta |^p  + \dots |\nabla^{k-1} \eta|^p)e^{p\kappa^\delta_\Sigma}  ) d\Sigma.$$

\noindent The almost complex structure $J$ is chosen such that in the cylindrical ends (ie. away from a compact domain) $J = J_0$ where $J_0$ is the standard complex structure on $\C^n$. The moduli space of $J$-holomorphic punctured discs is then given by the zero set of the Cauchy-Riemann $\mathcal{F}$ section. $$\mathcal{F}: \mathcal{B} \to \mathcal{E},\; \,\; u\mapsto \delbar_J(u)$$ 

\noindent The moduli space of parameterized punctured discs  $\widetilde{\mathcal{M}}_\Sigma$ can be described as $\mathcal{F}\inv (0)$ and the unparameterized moduli space $\mathcal{M}_\Sigma$ is the quotient $\widetilde{\mathcal{M}}_\Sigma/\text{Aut}(\Sigma)$. Note that the automorphism group of $\Sigma$ can be atmost 2 dimensional and the infinitesimal action of it on $\widetilde{\mathcal{M}_\Sigma}$ lies in the Banach tangent space.

\subsubsection{Space of Reeb chords and evaluation at puncture}

From Section \ref{sec:SFT} we know the asymptotic behavior near a puncture of a holomorphic map with finite Hofer energy is modelled on a Reeb chords or orbit. Thus, one can define an evaluation map at the punctures, taking values in the space of Reeb chords or orbits. We denote the space of Reeb chords in $Z$ with boundary on $\Lambda$ by $\mathcal{RC }(Z,\Lambda)$. In our situation where the Reeb field is generated by a free $S^1$ action on $Z$, the space $\mathcal{RC}$ can be easily described in terms of the length of the Reeb chord and starting point of the Reeb chord. Note that the length of the Reeb chord has to be an integer multiple of $\pi/n$ and the starting point has to lie on either of the two $T^{n-1}$ in $\Lambda$, thus $$\mathcal{RC}=\Bigg \{ \bigg(\frac{l\pi}{n},w \bigg) \bigg|\quad l\in \mathbb{N}, w\in \Lambda \Bigg \}.$$

\noindent We can define evaluation at puncture $p \in P$ from the moduli space of parameterized holomorphic discs to $S^{2n-1}$ that maps the holomorphic map $u$ which is asymptotic to the Reeb chord $\zeta$ at the puncture to the starting point $\zeta(0)$ of the Reeb chord.

\begin{equation*}
    ev_p : \widetilde{\mathcal{M}}_\Sigma \to \mathcal{RC},\qquad ev_p(u) = (k\pi/n,\zeta_p(0))
\end{equation*}
Here $\zeta_p$ is the Reeb chord of length $k\pi/n$ such that $u$ is asymptotic to the trivial strip over $\zeta_p$ at the strip-like end at $p$. We can also define an evaluation map on the charts of the Banach space of maps $\mathcal{B}$ as given in (\ref{eq:chartbanach}) as follows

\begin{align*}
    &EV_p : W^{k,p,\delta}({u_0}^* T\C^n,\partial {u_0}^* TT_\gamma) \oplus V^{|P|} \oplus H^{|P|} \to V\\
    &EV_p(w,v_1,\dots,v_p,\dots,v_l,h_1,\dots,h_p,\dots,h_l)=(v_p).
    \end{align*}

\begin{remark}
We can similarly define a space $\Maps^{k,p}_{\delta}(\Sigma_P, \R \times Z , \R \times \Lambda) $ for cylindrical manifolds and $\Maps^{k,p}_{\delta}(\Sigma_P, \widehat W , \widehat L)$ for manifolds with cylindrical ends over contact hypersurface. 
\end{remark}

\subsection{Linearization of Cauchy-Riemann map and Fredholmness}

In this subsection we show that the Cauchy-Riemann section  is a Fredholm map near zeroes of the section . More precise, we show that the vertical derivative at the points where the section $\mathcal{F}$ vanishes is a Fredholm map of Banach spaces. One can linearize $\mathcal{F}$ at a holomorphic map $u$ as usual in the classical case by choosing a connection on $\mathcal{E}$ induced from a connection  ( cylindrical in our case ) on $\C^n$. We assume that the punctures are non-removable ie. the punctures go to Reeb chords. Since we have already established exponential convergence to Reeb-strips, we first obtain an expression for linearization for trivial Reeb strip, $u_\gamma$ where $\gamma$ is a Reeb chord in $S^{2n-1}$. This will be an analog of  Proposition 6.24 in \cite{wendlSFT}.

Before we begin stating the result, we will define some terms so that it is consistent with Wendl's proposition. Wendl deals with a  cylindrical manifold $\R \times Z$ with framed Hamiltonian structure $(\omega,\lambda)$. In our case, the cylindrical manifold is obtained by setting $M=S^{2n-1}$ and the Hamiltonian structure is given the standard contact form $\lambda = \lambda_0$ on the sphere, where $\xi$ is the contact structure. The Reeb field on $M$ is denoted by $R$. $T(\R \times Z) $ is split into direct sum of two bundles, $\epsilon \oplus \xi$  where $\epsilon$ is generated by $R$ and $T\R$. $\Lambda$ denotes the two disjoint Legendrian tori, which gives the boundary condition. One can define an asymptotic operator associated to Reeb chords with boundary on $\Lambda$  similar to asymptotics associated to Reeb orbits. For Reeb orbits of period $T$, the asymptotic operator $A_\gamma$ was given by the hessian of the action functional $A_\omega$ on the space of loops in  $M$, see section 6.1.1 page 149 in \cite{wendlSFT} . In the case of Reeb chords, $T$ denotes the `length' of the Reeb chords, i.e. $\int \gamma^* \lambda$ . The connection $\nabla$ is chosen to be the Levi-Civita connection induced from the round metric on the sphere. We define  $A_\gamma  : \Gamma (\gamma^* \xi, T_{\gamma(0)}\Lambda, T_{\gamma(1)}\Lambda) \to \Gamma (\gamma^* \xi)  $  as follows $$A_\gamma(\eta)=-J(\nabla_t \eta - T\nabla_\eta R).$$ 
We can ignore the projection $\pi_\xi$ as we are working with a contact hypersurface instead of stable Hamiltonian structure.

\begin{lemma}[Analog of Wendl's 6.24 \cite{wendlSFT}] \label{analog}
The $J$-holomorphic trivial strip $u_\gamma : \R\times \interval \to \R \times Z$ for a Reeb chord $\gamma$ of length $T$ with boundary on $\R \times \Lambda$ has linearized Cauchy-Riemann operator $D_{u_\gamma}:\Gamma(u_\gamma^* \epsilon \oplus u_\gamma^* \xi  , \del u_\gamma^* T\R \oplus \del u_\gamma^* T \Lambda ) \to \Omega^{0,1}(\R \times \interval,u_\gamma^* \epsilon \oplus u_\gamma^* \xi )$ given by   $$(D_{u_\gamma} \eta) \del _s =\del_s \eta + 
\begin{bmatrix} i\del_t  & 0\\0&-A_\gamma \end{bmatrix} \eta.$$

\end{lemma}

\begin{proof}
Wendl's proof by splitting $D_{u_\gamma}$ into blocks given by the splitting $\epsilon \oplus \xi$ works the same  almost verbatim. The main difference is the boundary condition. Define the contact action functional for chords with boundary on  Legendrian $\Lambda$ by $\mathcal{A}_\lambda (\gamma) = \int \gamma^* \lambda $. We can check that the Hessian of this action functional, computed at a Reeb chord $\gamma$ with respect to $\nabla $, is  $-J(\nabla_t \eta - T\nabla_\eta R)$ where $T$ is the length of $\gamma$. The only change from Wendl's proof is that we are dealing with chords instead of orbits, so anytime an integration by parts is used one might encounter boundary terms $\lambda(\eta(1)) - \lambda(\eta(0))$ but they vanish since $\eta(0),\eta(1)$ lie on $ T\Lambda \hookrightarrow \ker \lambda $.
\end{proof}

\noindent Consider the Banach space isomorphisms $E: W^{k,p,\delta} \to W^{k,p}$ given by $f \mapsto e^{\delta s} f$. We have a slight abuse of notation here since $E$ is a map from  $W^{1,p,\delta}(u_\gamma^* T(\R \times Z),\del u_\gamma^* T(\R \times \Lambda))$ and of $\Omega^{0,1}_\delta (u_\gamma ^* T(\R \times Z))$. With this notation, we have 
$$ED_{u_\gamma}E\inv  = \del_s + \begin{bmatrix} i\del_t + \delta  & 0\\0&-A_\gamma + \delta \end{bmatrix}. $$

\begin{lemma} \label{lemma:goodeigen}
For any $   \delta < \frac{\pi}{n}$, we have that $ED_{u_\gamma}E\inv$ is an operator of form $\del_s + A$ on the space of Sobolev sections $W^{1,p}$ where $A$ is a  non-degenerate operator.
\end{lemma}

\begin{proof}

Note that the eigenvalues of $\begin{bmatrix} i\del_t & 0\\0&-A_\gamma \end{bmatrix} $ are given by eigenvalue of $i\del_t$ and $A_\gamma$. The eigenvalue of $A_\gamma$ can be explicitly computed to be $\pi \Z$, notice since we are dealing with  a Morse-Bott situation, 0 is actually an eigenvalue. Note that $A_\gamma = -J \hat {\nabla}_t$ where $\hat{\nabla}$ is the unique symplectic connection for which parallel transport is the linearized Reeb flow , (see \cite{wendlSFT} ), we have that eigenfunction for eigenvalue $\lambda_e$ is given by $ \eta(t) = (\phi_t)_*(e^{i\lambda_e t} \eta(0 ))$ where $\phi$ is the Reeb flow. Thus, to satisfy the boundary conditions we need $\lambda_e \in \pi \Z$. The eigenvalues for $i\del_t$ are a subset of $\frac{\pi}{n}\Z$ since the smallest $\theta$ such that $e^{i\theta} \Lambda \hookrightarrow \Lambda$ is $\frac{\pi}{n}$. Thus, if we choose $\delta$ in the given range, we'd have $ED_{u_\gamma}E\inv = \del_s + A$ where $A$ is non-degenerate.
\end{proof}

\subsubsection{Trivialization }

We explain the choice of trivialization of bundles we make to simplify the linearization. We need a suitable trivialization under which $D_u$ is equal to  $\del_s + J_0 \del_t + A$. We state a version of Lemma 2.13 in \cite{schm} to create a family of trivialization, we leave the proof to the reader since it follows the idea in \cite{schm} very closely.

\begin{lemma} \label{lemma:trvlz}
For any $p\in L_Y$ and $R>0$, there is an open set $V$ in $Y$ such that the tangent bundle of $U = \pi_Y \inv (V)$ has the following trivialization 
$$\Phi : [R,\infty) \times \interval \times U \times \R^{2n} \to T(\R \times Z)|_{U},\qquad (s,t,q,v) \mapsto \Phi_{s,t}(q)v\in T_q(\R \times Z)$$ that satisfies the following
\begin{itemize}
    \item $J_{s,t}(q)\Phi_{s,t}(q) = \Phi_{s,t}J_0$,
    \item $g_{cyl} (\Phi_{s,t}(q)\xi ,\Phi_{s,t}(q)\xi' ) = g_0(\xi,\xi')$,
    \item $T_q (\R \times \Lambda) = \phi_{s,k} (\R ^n \oplus \{0 \})  $ for all $q\in \R \times \Lambda \cap U$ and $k=0,1$,
\end{itemize}

where $J_0,g_0$ are the standard complex structure and metric on $\R^{2n}$.
\end{lemma}

\noindent  Recall when $B$ is a subbundle of $A$ and $D$ is a subbundle of $C$, then by an isomorphism $F$ of the pair of bundles $F: (A,B) \to (C,D)$ refers to an isomorphism $F:A\to C$ such that its restriction to $B$ is an isomorphism to $D$.  Let's assume that $\pi_Y \of u$ converges to the point $p \in L_Y \sub Y$ as $s\to \infty$ and choose a trivialization $\Phi$ as in Lemma \ref{lemma:trvlz}. Using this trivialization, we get the following trivialization of the pair of bundles over the half-strip $[R,\infty)\times \interval$.
\begin{align} \label{eq:pullbacktriv}
    \wt \Phi : (u^*T(\R \times Z), \del u^* T(\R \times \Lambda)|_{[R,\infty)\times \interval}) \to ([R,\infty)\times \interval \times \R^{2n} ,\\ [R,\infty)\times \{0,1\} \times (\R^n \oplus \{0\}) ) \nonumber
\end{align}

\subsubsection{Semi-Fredholmness by local estimates on cylindrical ends}\label{subsec:semifred}

We call a linear map \textit{semi-fredholm} if it has closed image and finite dimensional cokernel. In this subsection, we will show that the linearization of the Cauchy-Riemann section is semi-fredholm by combining local estimates.

We will use the following classic lemma to prove that $D_u$ is Fredholm. This has been proved in numerous textbooks and classical references, e.g., see Lemma A.1.1 in \cite{jholobook}. 
\begin{lemma}\label{semilemma}
If $K$ is a compact operator, $T$ is a bounded operator, and $c>0$ such that $$|x| \leq c(|Tx| + |Kx|)$$ then $T$ is semi-fredholm.
\end{lemma}

\noindent If  the domain surface $\Sigma$  were compact, by standard methods in \cite{jholobook} one has  $$|x|_{W^{1,p}(\Sigma)}\leq c(|D_ux|_{L^p(\Sigma)} + |x|_{L^p(\Sigma)}),$$ and by the Rellich-Kondrachov compact embedding, we have that the inclusion of $W^{1,p} (\Sigma) \hookrightarrow L^p (\Sigma)$ is compact for a compact surface $\Sigma$ with smooth boundary. The argument goes as follows, from Proposition B.4.9 in \cite{jholobook} we have a local estimate, which is then patched together for a finite covering of the compact space. This argument is outlined in Proposition C.2.1 in the book.

 The above approach will not directly work when the domain is non-compact. We will show a stronger inequality in the cylindrical setting to avoid the use of Rellich-Kondrachov theorem. We will use the approach in \cite{wendlSFT} to show a similar inequality for a compact subdomain away from the strips and on the strips we will show there is an inequality  (see  Lemma 4.18 in \cite{wendlSFT})  
            \begin{equation} 
            ||x||_{W^{1,p}(\R \times [0,1])}\leq c||D_ux||_{L^p(\R \times [0,1])}. \tag{$*$}\label{eq:star}
            \end{equation}

We will use the above estimate while patching local estimates on the whole domain. From exponential convergence we have that $u$ converges to $u_\gamma$ on a strip-like end of $\Sigma$ for some Reeb chord $\gamma$, thus it is enough to prove \eqref{eq:star} for $u_\gamma$. To verify this, notice trivialization of $u_{\gamma} ^* T(\R \times Z)$ and $u ^* T(\R \times Z)$ as given in equation (\ref{eq:pullbacktriv}) is exponentially close to each other that contributes to an exponentially decaying $0^{th}$ order term. Now we use the fact that for $\delta >0$, the inclusion $W^{k,p,\delta} \to  W^{k-1,p}$ is a compact operator and thus the difference is a compact operator.

\begin{lemma}\label{operatoriso}
Let $D=\del_s+  A(t)$ be an operator of functions over the strip such that $A$ is an unbounded self-adjoint and non-degenerate operator on $W^{1,p} ([0,1],\C^n,\R^{n-1}\times\lambda_e \R_+, \R^{n-1}\times \R_+)$, the Sobolev space of functions from $\interval$ to $\C^n$ with given boundary conditions.Then $$D:W^{1,p} (\R \times \interval,\C^n ,\R^{n-1}\times \lambda_e \R_+ ,\R^{n-1}\times \R_+  ) \to L^{p} (\R \times \interval,\C^n)$$ is an isomorphism. Here $\lambda_e = e^{k\pi/n}$ for some k. 
\end{lemma}

\begin{proof}
This lemma is a version of Proposition 8.7.3 in \cite{AudDam}, Lemma 2.4 in \cite{salamonNotes}. The proof runs the same, except the fact in the special case of $p=2$  we'd need to use density of $W^{1,p} (\C^n ,\R^{n-1}\times \lambda \R_+ ,\R^{n-1}\times \R_+  ) $ inside $ L^{2} (\C^n)$. To verify the density, one can check that $C^3$ maps from $\interval$ to $\C^n$ with the boundary conditions on the given subspaces is dense in $C^3(\interval,\C^n)$ with $L^2$ metric. Clearly we have that $C^3(\interval,\C^n)$ is a subspace of $W^{1,p}$, thus density of smooth functions in $L^2$ gives us that $W^{1,p} (\C^n ,\R^{n-1}\times \lambda \R_+ ,\R^{n-1}\times \R_+  ) $ is dense inside $ L^{2} (\C^n)$. The other steps in the proof of the lemma run exactly the same as that given in the references.
\end{proof}

\begin{cor}
\label{corr:localcylestimate}
If $u$ is a holomorphic strip $\R \times \interval \to \R \times Z$ that exponentially converges to a trivial Reeb strip $u_\gamma$  (i.e., satisfies equation (\ref{expConv})), we have that for all $ \delta< \frac{\pi}{n}$ there exists a $c>0$ such that $$||x||_{W^{1,p,\delta}(\R \times [0,1])}\leq c||D_ux||_{L^{p,\delta}(\R \times [0,1])}.$$
\end{cor}

\begin{proof}
Note that since we have exponential convergence, it is enough to check for $D_{u_\gamma}$. We know $u_\gamma(s,t) = e^{\frac{k\pi}{n}(s+it)}(x_1,\dots,x_n)$ for some $x_i \in \C$ such that $x_1.x_2.\dots x_n \in \R_{\neq 0}$. Thus, we see that we can choose a trivialization of $u_\gamma ^* T\C^n$ such that we get the boundary conditions as in Lemma \ref{operatoriso}. Now letting $D'=ED_{u_\gamma}E\inv$, from Lemma \ref{lemma:goodeigen} we have that $D'$ fits the hypothesis in Lemma \ref{operatoriso}. The fact that $D'$ is an isomorphism gives us the inequality 
$$||\eta ||_{W^{1,p}} \leq c || D' \eta ||_{L^p}. $$
Since $E$ is a Banach isomorphism from the weighted to unweighted Sobolev spaces, we have 
$$ ||\eta ||_{W^{1,p,\delta}} \leq c || D_{u_\gamma} \eta ||_{L^{p,\delta}}, $$ for a possibly different c.
\end{proof}

\begin{Theorem}
$D_{u}: T_u\mathcal{B} \to \mathcal{E}_u $ is semi-fredholm when $ \delta < \frac{\pi}{n}$.
\end{Theorem}

\begin{proof}
Define $E$ similar to as \ref{lemma:goodeigen}, but now on the whole $\Sigma$ by choosing a function $f:\Sigma \to \R$ that is equal to $\delta s$ on the strip-like ends and zero away from them. This defines  Banach isomorphism

\begin{equation}
\label{eqn:soboleviso}
E : W^{k,p,\delta}(u^* T\C^n) \rightarrow W^{k,p}(u^* T\C^n)\quad \eta \mapsto e^f\eta.
\end{equation}

\noindent As $E$ is a Banach isomorphism, it is enough to show $D' := ED_uE\inv$ is semi-fredholm.

\noindent We will use Lemma \ref{semilemma} to claim semi-fredholmness. We can cover $\Sigma$ by finitely many charts away from the punctures, call the union of the charts $\Sigma_o$ and use Proposition B.4.9 (ii) in \cite{jholobook} to get an estimate 
$$ ||\eta|| \leq c(||\eta||_{W^{0,p}(\Sigma_o)}|| +  ||D_u (\eta) ||_{W^{0,p}(\Sigma_o)}) \qquad  \forall \eta \in W^{1,p}_0((\Sigma_o)). $$
This gives us an estimate for $D'$ as on $\Sigma_o$ for a possibly different $c$, $E$ is just multiplying by a bounded function.\begin{equation}
    ||\eta|| \leq c(||\eta||_{W^{0,p}(\Sigma_o)}|| +  ||D' (\eta) ||_{W^{0,p}(\Sigma_o)}) \qquad  \forall \eta \in W^{1,p}_0((\Sigma_o)) \tag{*}. \label{1}
\end{equation} 

We can assume that the strip-like ends in $\Sigma \setminus \Sigma_o$ under the fixed chart is biholomorphic to $[R+1,\infty) \times \interval$  for some fixed $R>0$. For each puncture $q$, let $Z^q_R$ be the strip neighborhood $[R,\infty) \times \interval$ under the fixed charts. 
From  Corollary \ref{corr:localcylestimate} and Lemma \ref{lemma:goodeigen} we see that 
\begin{equation}
    ||\eta||_{W^{1,p}(\R \times [0,1])}\leq c_q'||D' \eta||_{L^p(\R \times [0,1])} \qquad \forall \eta \in W^{1,p}_0(Z^q_R). \tag{**} \label{2}
\end{equation}

Here we extend $\eta $ by 0 over the whole strip before using Corollary \ref{corr:localcylestimate}. 
Now a standard argument of using bump function $\beta$ that is supported on  $\Sigma_o$, for any $\eta \in W^{1,p}(u^* T(\R \times Z))$, we can use $\beta\eta$ and $(1-\beta)\eta$  in inequalities (\ref{1}) and (\ref{2}) to get 

\begin{equation*}
    || \eta||_ {W^{1,p}(\Sigma)} \leq C (||D'(\eta)||_{W^{0,p}(\Sigma)} + ||\eta ||_{W^{1,p}(\Sigma_o)}).
\end{equation*}
Since  closure of $\Sigma_o$ in $\Sigma$ is compact with smooth boundary, we have that $W^{1,p} (\Sigma_o) \hookrightarrow W^{0,p} (\Sigma_o)$ is a compact inclusion, thus the restriction on $\eta \mapsto \eta|_{\Sigma_o}$ is a compact operator. Thus, from \ref{semilemma} we have that $D'$ is semi-fredholm, thus $D_u$ is semi-fredholm for the given choice of $\delta$

\end{proof}

\subsubsection{Fredholmness}
In this subsection, we will use a similar argument as in subsection \ref{subsec:semifred} to the adjoint of the linearization and show that the linearization is a Fredholm operator. To show that $D_u$ is Fredholm, it is enough to show $D'= ED_uE\inv$ is Fredholm where $E$ is from equation (\ref{eqn:soboleviso}). We use a Hermitian structure $h_{cyl}= g_{cyl}(\_ ,\_) -ig_{cyl}(\_\,,J\_)$ on $T\C^n$ induced from the cylindrical metric $g_{cyl}$ to define a formal adjoint $D'^*$ and show that $D'^*$ is semi-Fredholm. This idea has been widely used in literature to prove Fredholmness of Cauchy-Riemann operators coming out of linearizations, see \cite{wendlSFT}, \cite{jholobook} . We define $C _{0,T_\gamma}( \ol{\operatorname{Hom}}( T\Sigma,u^*T\C^n))$ as follows

\begin{align*}
    C _{0,T_\gamma}( \ol{\operatorname{Hom}}(T\Sigma,u^*T\C^n)) = \bigg\{\alpha \in  \Gamma (\ol{\operatorname{Hom}}(T\Sigma,u^*T\C^n)) \bigg|  \quad \alpha( T \del \Sigma ) \sub \del u ^* T T_\gamma , \\ \alpha \text { has compact support}\bigg\}.
\end{align*}

We use a Hermitian metric to define an adjoint of the linearization. The Hermitian structure $h_{cyl}$ induces an $L^2$ pairing on the space of sections $\Gamma(u^*T\C^n,\del u^*T\C^n)$ and $\Gamma (\ol{\operatorname{Hom}}(T\Sigma,u^*T\C^n))$ given by the  following integral,

\begin{align*}
    \langle \eta,\xi \rangle := \text{Re} \int_\Sigma \langle \eta, \xi \rangle d\Sigma.
\end{align*}

\noindent The formal adjoint $D'^*$ is defined as the unique operator on the Sobolev completion $W^{k,p}_{T_\gamma} ( u^*\ol {\operatorname{Hom}}(T\Sigma,u^*T\C^n))$  that satisfies 

\begin{align*}
    \langle D' \alpha , \beta   \rangle = \langle \alpha , D'^* \beta \rangle \text {   for all  } &\alpha \in \Gamma_0(u^*T\C^n,\del u^*T\C^n), \\ 
    &\beta \in C _{0,T_\gamma}(\ol {\operatorname{Hom}}(T\Sigma,u^*T\C^n)).
\end{align*}

\begin{remark}
 Another approach of defining the formal adjoint $D_u^* $ is to define it on the space $W^{k,p,-\delta}_{T_\gamma} ( \ol{\operatorname{Hom}}(T\Sigma,u^*T\C^n)),$  where we use a negative Sobolev weight. We can define isomorphism $E'$ similar to $E,$ but with negative weight such that $E'D_u^* E'^{-1} = D'^*$.
\end{remark}

\noindent We can follow the same proof as in Proposition 4.32 in \cite{wendlSFT} to get the following lemma about $D_u^*$ which finishes the proof of Fredholmness.

\begin{lemma}
$D'^*$ is semi-Fredholm and we have the following isomorphisms $\coker D' = \ker D'^*$, $\coker D'^* = \ker D' $.
\end{lemma}

\noindent As a corollary of the previous lemma, we have the following theorem

\begin{Theorem} \label{thm:fred}
The linearization $D_u$ of $\delbar$ at any holomorphic map $u$ is a Fredholm operator.
\end{Theorem}

\subsection{Domain dependent almost complex structure on the outside}\label{subsec:domdepacs}

In practice, attaining regularity is not possible without allowing some perturbation of the almost complex structure. We see in Section \ref{sec:reg} that one can achieve regularity for the inside and neck pieces for the standard almost complex structure $J_0$ in $\C^n$, but we need to allow perturbations of the outside piece to achieve regularity of the broken discs and strips. We use the perturbation scheme of domain dependent almost complex structures as initially done in \cite{CMtrans} for the closed case and later generalized to the open case with boundary on Lagrangian in \cite{charwoodstabil},\cite{floerflip}.

\noindent Notice that we can identify $ \widehat M_{out}$ with $M\setminus \{ p \}$ where $p$ is the center of the Darboux ball whose boundary we do the neck-stretching along. In this identification, $\widehat L_{out}$ corresponds to the Lagrangian $T_\gamma$ obtained by the path $\gamma(t) = t$ for $t\in (-\e,\e) \setminus \{ 0 \}$ in the Darboux chart. The closure of $\widehat L_{out}$ is denoted by $L_{sing}$, it is locally a conical manifold, a double cone over the torus $T^{n-1}$ given by $T_{(-\e,\e)}$. We need to stabilize the components in the outside domain that have boundary lying on $L_{out}$ so that we can use domain dependent perturbation on these components. We follow the scheme of obtaining marking by intersection with a divisor similar to \cite{CMtrans}  ,\cite{charwoodstabil},\cite{floerflip}. The treatment of domain dependent perturbations presented here assumes some familiarity with the idea of coherent choice of perturbation datum for Lagrangian intersection Floer theory as done in \cite{charwoodstabil}.

\begin{defn}[Nice divisor adapted to $J$] \label{def:nicediv}
Given a pair of transversely intersecting Lagrangians $L,K$ of $M$ such that $L$ is locally mutable, a submanifold $D$ of codimension two that doesn't intersect $L \cup K \cup L_{sing}$ is called a `` nice divisor" if $M\setminus D$ is exact symplectic manifold and $L_{sing}$ is exact Lagrangian in $M\setminus D$. A nice divisor is called ``adapted" to a complex structure $J$ if the tangent bundle $TD$ is preserved by $J$.
\end{defn}

\begin{remark}
In $\P^n$, if we take $L$ to be the Clifford torus, then there is a nice divisor $D$ adapted to the standard $J_0$. We can take two fibers of the map $\pi_m$ in some standard affine chart of $\P^n$ and perturb them to obtain a smooth divisor in $\P^n$ that makes $L$ exact in the complement.         
\end{remark}

From now on, we assume the existence of a nice divisor. If $u$ is a non-constant holomorphic map in the outside piece with boundary on $L_{out}$, then it must intersect $D$ as otherwise exactness of $L_{sing}$ would imply symplectic area of $u$ is 0. 

\begin{defn}[Stabilizing combinatorial type of broken strip]\label{defn:stabilz}
Given combinatorial type $\Gamma$ of a broken strip where each outside piece stable, we define the stabilized type $\Gamma^{st}$ by forgetting all the inner and neck pieces that are not stable along with the puncture nodes that connected the unstable piece to a stable piece. If removing such a puncture node renders an otherwise stable piece into an unstable piece, we forget that as well in $\Gamma^{st}$. 
\end{defn}

We will now define the notion of domain dependent perturbation of almost complex structure. Let $k\geq 1$ be natural number, we define $\M_k$ to be the moduli space of nodal strips with $k$ interior markings  and $\mathcal{U}_k \xrightarrow{\wt\pi_{k}} \M_{k}$ is the universal moduli space whose fiber over $S\in \M_{k}$ is the set by pairs $(S,z)$ where $z$ is a point on $S$.  Fix an almost complex structure $ J_{fix}$ on $\widehat M_{out}$ that is equal to the cylindrical almost complex structure $J_0$ on the cylindrical end.  We denote $\mathcal{J}_{cyl}(\widehat M_{out},\omega,D)$ to be the set of almost complex structures satisfying the following,
\begin{enumerate}
    \item compatible with the symplectic form $\omega$ on $\widehat M_{out}$,
    \item equal to $J_0$ in the cylindrical end of $\widehat M_{out}$,
    \item $D$ is adapted to it.
\end{enumerate} 

\begin{defn}[Domain dependent almost complex structure for broken strips ]\label{def:domdepstruc}
A domain dependent almost complex structure of a strip with $k$ interior markings is a smooth map $\mathcal{J}_k: \mathcal{U}_k \to \mathcal{J}_{cyl}(\widehat M_{out},\omega,D)$ such that near the puncture nodes, the almost complex structure is equal to $J_{fix}$, i.e., there is an open neighborhood $U$ in $\mathcal{U}_k$ near the locus of puncture nodes such that  $\mathcal{J}_k (w) = J_{fix}$ for $w\in U$ . 
\end{defn}

Now we describe glued domain dependent almost complex structures. Given a domain dependent almost complex structure on a stabilized combinatorial type of broken strip $\Gamma^{st}$ such that near the puncture nodes $J_k$ is equal to $J_{fix}$, we can get a domain dependent almost complex structure on the unbroken  $T$-stretched manifold $M^T$ by gluing $\Gamma^{st}$ at the puncture nodes. A $T - $ gluing across a puncture node, where $T>0$ is a parameter measuring the ``length of the neck", is obtained by identifying truncated neighborhoods (deleting $(T,\infty)$ from positive end and $(-\infty,-T)$ from negative end) of the strip/cylinder coordinate near both sides of the puncture node.

We explain the gluing construction more precisely in the case where there is only one puncture node, the general construction is done similarly. Fix a parameter $T>0$. For now, assume $\Gamma^{st}$ has only one puncture node. Let $\Gamma^{st}_{gl}$ be the combinatorial type obtained from gluing $\Gamma^{st}$ across a puncture node. Thus, the moduli space of strips $\M_{\Gamma^{st}}$ of type $\Gamma^{st}$ is a lower dimensional stratum in the compactified moduli space $\ol \M_{\Gamma^{st}_{gl}}$. We can choose a chart of $\M_{\Gamma^{st}_{gl}}$ near the strata $\M_{\Gamma^{st}}$ by first choosing strip-like parameterization of the surfaces near the puncture node and then gluing the strips.  Let $[\mathfrak{s}]\in \M_{\Gamma^{st}}$, and $C_{[\mathfrak{s}]}$ be the surface corresponding to $[\mathfrak{s}]$, let $p$ be the puncture node lying between the components $C_-,C_+$ of $C_{[\mathfrak{s}]}$. Let $U_\pm$ be  neighborhoods of $p$ in $C_\pm$ such that we can choose charts  $\phi_\pm$ on them as follows,
\begin{align*}
    \phi_-:U_- \to (-\infty ,0] \times \interval \\
    \phi_+:U_+ \to [0,\infty) \times \interval.
\end{align*}
\noindent We obtain the $T$-glued domain $C^T_{[\mathfrak{s}]}$ by removing $\phi_-((-\infty,-T)\times \interval)$  and $\phi_+((T,\infty)\times \interval) $ from the components $C_\pm$ of $C_{[\mathfrak{s}]}$ and gluing the remnant surfaces by the relation $(s,t) \sim (s+T,t)$.

We now explain how to glue domain dependent almost complex structure. Since $\J_k$ is equal to $J_{fix}$ near the locus of the puncture node, we can choose a strip-like neighborhood $U_\pm$ as above, so that $\J_k$ restricted to $U_\pm$ is equal $J_{fix}$. Recall that we do not perturb the almost complex structure on the inside piece. We define a $T$-glued domain dependent almost complex structure $J^T_k$ on the $T$-glued domain $C^T_{[\mathfrak{s}]}$,by defining $J^T_k(z)$ to be equal to the $T$-stretched almost complex structure in the $J^T$ neck region for all $z$ and for $z\in$ the neck region of the domain, define $J_k(z)$ to be equal to $J_{fix}$ on the outside piece of $M^T$.

\begin{align*}
    J^T_k(z) (m) =
    \begin{cases}
     J_{in}  & \text{ if }  m \in M_{in} \\
     J_0        & \text{ if } m \in [-T\times T] \times Z\\
     J_k(z)(m)    & \text { if } m\in M_{out}
     \end{cases}
\end{align*}

\begin{defn}[Coherent families of perturbation data for broken strips]
We call the set $\mathcal{P} = \{ \mathcal{J}_k \}_{k=0}^{\infty}$ where in the $k=0$ case we choose a $t-$dependent complex structure, and for $k\geq 1$,  $\mathcal{J}_k$ is a domain dependent almost complex structure for broken strips, a coherent family of perturbation data if the $T-$glued domain dependent almost complex structure forms a coherent family of perturbation for $M^T$ as in Definition 3.14 in \cite{charwoodstabil} for any $T>0$.

\end{defn}

We explain the notion of being holomorphic with respect to a domain dependent choice of almost complex structure.  A map $u:S_{\textbf{z}} \to \widehat M_{out}$  from a $k$ marked strip $S_{\textbf{z}} = \wt \pi_k \inv ([\textbf{z}]) $ where $[\textbf{z}] \in \M_k$ is called  holomorphic with respect to the domain dependent almost complex structure $\mathcal{J}$ if it satisfies the domain-dependent version of the Cauchy Riemann equation 
\begin{equation*}
   du(z) + J_k(z,(u(z)))\of du \of j_S = 0 .
\end{equation*}

We show that the expected dimension of moduli space of holomorphic maps with respect to domain-dependent choice almost complex structure is equal to domain independent almost complex structure if we add divisor constraints on the marked point.   Recall, for a fixed complex structure, the expected dimension of the moduli space of discs or strips at a point $u$ depended solely on the topological data of $u$. We show that the same formula can be applied to the moduli space of holomorphic broken strips with respect to a domain-dependent perturbation, with the constraint of interior markings mapping to a nice divisor $D$. Denote the expected dimension of the moduli space at a point $u$ that is of homology class $A$ by $\vir(A)$.  Given a (relative) homology class $A$, the expected dimension of $\M^{dd}_{k}(\widehat M_{out})$ at any element $u$ that is of the homology class $A$ is $\vir([A])+\dim(\M_k)$. If we require the interior markings to map to $D$ to create a constrained moduli space $\M_{k}(M_{out},D)$, the expected dimension of $\M_{k}(M_{out},D)$ at $u$ of homology class $A$ is $\vir(A)+\dim(\M_k) - 2k = \vir(A)$.

We introduce some notations regarding domain-dependent almost complex structures. We denote the space of coherent perturbation datum $\mathcal{P}(\widehat M_{out},D)$.  We denote the moduli space of $J_k$ holomorphic maps by $\M^{dd}_k(\widehat M_{out})$ and the moduli space of broken strips or broken disc with point constraint with respect to a coherent perturbation data $\mathcal{P}$ with the constraint of interior markings going to a nice divisor $D$ as $\M_{Br,\mathcal{P}}(\mathbb{M}[l],\mathbb{L}[l],K,D)$ and $\M_{Br,\mathcal{P}}(\mathbb{M}[l],\mathbb{L}[l],p,D)$. In Subsection \ref{subsec:outtran} we prove that there is a choice of a coherent perturbation datum that makes every rigid broken strip and rigid broken strip with point constraint regular.

\subsection{Moduli space of broken strips and discs}
Now that we have dealt with the analysis of the space of maps and Fredholmness of the Cauchy-Riemann section, we will now describe the moduli space of broken discs and strips.  In this subsection, we will explain how the moduli space of broken strips and discs can be cut out as a zero set of a section similar to moduli space of holomorphic discs. Let $\mathfrak{D}=(\mathbb{D}_{in},\mathbb{D}_1,\dots , \mathbb{D}_{out})$ denote a nodal disk or strip with $l$ levels as in Definition \ref{def:brokdisc}. We denote the set of nodes between two components of different label as $P$ and for all $p\in P $, we use $p^\pm$ to view $p$ as an element in the two different components. The moduli space of broken strips with domain $\mathfrak{D}$ to the broken manifold $\mathbb{M}[l]$ with boundary on the broken Lagrangian $\mathbb{L}[l]$ and some boundary components of $\mathbb{D}_{out}$ mapping to the compact Lagrangian $K$ in $\mathbb{M}_{out}$ is denoted by $\mathcal{M}_{Br}(\mathfrak{D},\mathbb{M}[l],\mathbb{L}[l],K)$ and the moduli space of broken discs is denoted by $\mathcal{M}_{Br}(\mathfrak{D},\mathbb{M}[l],\mathbb{L}[l])$.

Assume $\Gamma^{st}$ is the stabilized combinatorial type of the domain, see Definition \ref{defn:stabilz}. Let $C$ denote the nodal surface $\mathfrak{D}$ regarded as a smooth surface with nodes. Let $C_i$ denote the surface in $i\th$ piece of the domain, $\mathbb D_i$, regarded as a smooth manifold with nodes and let 
\begin{equation}\label{domtriv}
    \mathcal{U}_{\Gamma^{st}_i,l} \to \mathcal{M}_{\Gamma^{st}_i,l} \times C_i
\end{equation}
be a collection of trivializations, indexed by $l$, of the universal moduli space of Riemann surfaces of type $\Gamma$. These trivializations induce a family of complex structures on the fibers of the universal space \begin{equation}\label{domchoice}
    \M_{\Gamma^{st}_i,l} \to \mathcal{J}(C_i), \quad m\to j(m).
\end{equation}
\noindent Let $\mathcal{E}_i$ denote the bundle of $(0,1)$ forms over the space $\Maps^{k,p}_{\delta}(\Sigma_i,M_i ,\R \times \Lambda)$, i.e. the fiber $\mathcal{E}_i$ over a map $u_i$ is 
\begin{equation*}
    \mathcal{E}_{i,u}:= W^{k,p}_\delta ( \overline{\operatorname{Hom}}( T\Sigma_i,u_i^*TM_i)).
\end{equation*}
In the equation below, the notation $\M_\Gamma$ stands for a choice of a complex structure on the domain $C$ induced from the choice of trivializations we made in equations (\ref{domtriv}) and (\ref{domchoice}).  We can define the following map by setting it to $\delbar_{j(m),J_i}$ on each component 

\begin{align*}
    \delbar^{strip}_{Br} :\mathcal{M}_{\Gamma^{st}}\times \Maps^{k,p}_{\delta}(\Sigma_{in},\C^n,T_\gamma) \times \Maps^{k,p}_{\delta}(\Sigma_1,\R \times S^{2n-1},\R \times \Lambda)\\ \times \dots \Maps^{k,p}_{\delta}(\Sigma_{out},\widehat M_{out},L_{out},K)\to &\mathcal{E}_{in}\times \dots \mathcal{E}_{out}.
\end{align*}

\noindent Note that the linearization of $\delbar_{Br}$ is a Fredholm operator because we know that linearization of $\delbar$ is Fredholm from Theorem \ref{thm:fred} and $\mathcal{M}_\Gamma$ is finite dimensional. The set $\delbar_{Br} \inv (0)$ consists of holomorphic maps from a broken domain $\mathfrak{D}$  of combinatorial type $\Gamma$ with the correct Lagrangian boundary condition but at the nodes between two levels, the maps might not match up correctly ie. they might be asymptotic to different Reeb chords. Recall that we have an evaluation map at each puncture to the space of Reeb chords $\mathcal{RC}$, we collect such evaluation maps to create a ``big evaluation map"

\begin{equation*}
    ev_{big}: {\delbar^{strip}_{Br}}\inv (0) \to ( \mathcal{RC}_- \times \mathcal{RC}_+ )^{|P|}.
\end{equation*}

\noindent The space $ev_{big}\inv (\Delta^{|P|})$ where $\Delta$ is the diagonal manifold in $\mathcal{RC}_- \times \mathcal{RC}_+$ is precisely the moduli space $\mathcal{M}_{Br}(\Gamma,\mathbb{M}[l],\mathbb{L}[l],K)$. The linearization of $ev_{big}$ is obtained by collecting $EV_{-} \times EV_{+}$ for each node between two levels. 

Similar to broken strips, we can obtain the moduli space of broken discs as a fiber product obtained from a zero section. The only modification we need to do is define the following $\delbar$ section :

\begin{align*}
    \delbar^{disc}_{Br} :\mathcal{M}_{\Gamma}\times \Maps^{k,p}_{\delta}(\Sigma_{in},\C^n,T_\gamma) \times \Maps^{k,p}_{\delta}(\Sigma_1,\R \times S^{2n-1},\R \times \Lambda)\\ \times \dots \Maps^{k,p}_{\delta}(\Sigma_{out},\widehat M_{out},L_{out})  \to &\mathcal{E}_{in}\times \dots \mathcal{E}_{out}
\end{align*}
where $\mathcal{E}_{out}$ is the obvious bundle of $(0,1)$-forms over the space $\Maps^{k,p}_{\delta}(\Sigma_{out},\widehat M_{out},L_{out})$.

We collect some definitions below. In all the following definitions, we work with a fixed combinatorial type $\Gamma$ and the domains are allowed to vary only in the strata of the associahedra of the moduli of domains determined by the type $\Gamma$.

\begin{defn}[Regular broken strip] \label{defn:brokendiscreg}
We call a broken strip $\mathfrak{u} \in \mathcal{M}_{Br}(\Gamma,\mathbb{M}[l],\mathbb{L}[l],K) $ regular if the linearizations of $\delbar_{Br}$ is surjective at $\mathfrak{u}$  and the evaluation map, $ev_{big}$,  is transverse to the diagonal $\Delta^{|P|}$.

\end{defn}

\noindent Assuming transversality of everything, the expected dimension of the moduli spaces $\mathcal{M}_{Br}(\Gamma,\mathbb{M}[l],\mathbb{L}[l],K)$ and $\mathcal{M}_{Br}(\Gamma,\mathbb{M}[l],\mathbb{L}[l])$ at a point $\mathfrak{u}$ is $\Ind(\wt D_\mathfrak{u}) - |P|(n-1)$ where $\wt D_u$ is 
the linearization of $\delbar^{strip}_{Br}$ or $\delbar^{disc}_{Br}$ at $\mathfrak{(D,u)}$. 

\begin{defn}[Virtual dimension]
For a broken strip or a disc $(\mathfrak{u,D})$, the virtual dimension is given by the formula $$\vir (\mathfrak{u}) = \Ind(\wt D_\mathfrak{u}) - |P|(n-1) - \dim \aut(\mathfrak{D}).$$
\end{defn}

\begin{defn}[Rigid broken strip]
We call a broken strip $\mathfrak{u}$ rigid if $\vir(\mathfrak{u}) = 0$. The moduli space of rigid strips is denoted as $\mathcal{M}_{Br}(\Gamma,\mathbb{M}[l],\mathbb{L}[l],K)_0$.

\end{defn}

\begin{defn}[High-index broken strip]
We call a broken strip $\mathfrak{u}$ high-index if $\vir(\mathfrak{u}) > 0$.

\end{defn}

In the unbroken compact case, the moduli space Maslov 2 discs with a boundary constraint of passing through a fixed point is of expected dimension $0$ and contributes to the disc potential. We define a notion of point constrained broken disc so that we can later compare counts of Maslov two discs passing through a point of the unbroken Lagrangian $L$ with the broken Lagrangian $\mathbb{L}$. Given a point $p\in L_{in}$ that is not in the cylindrical part of the Lagrangian $L_{in}$, we define the notion of broken disc with point constraint, the choice of $p$ is not relevant since $L$ is assumed to be connected.

\begin{defn}[Broken disc with point constraint]
Let $q$ be a boundary marking on one of the interior pieces of  $\mathfrak{D}$ which gives us an evaluation map $ev_q$ with values in $L_{in}$. The moduli space of broken discs with point constraint, $\M_{Br}(\mathbb M [l],\mathbb L [l],p)$ is defined as the space $ev_{q}\inv (p)$.

\end{defn}

\noindent The expected dimension of the point constrained moduli space at a broken disc $(\mathfrak{u,D})$ is $\vir(\mathfrak{u}) - n$.

 \begin{defn}
  
 We say a broken disc with point constraint is \textit{rigid} if $\vir(\mathfrak{u}) = n$.

 \end{defn}

\begin{defn}[Compactified energy]
For a $two$-level broken strip $(\mathfrak{u,D})= ((u_{in},u_{out}), \mathfrak D)$. The compactified energy $E_c(\mathfrak{u})$ is defined as 
$$E_c(\mathfrak u) = E(\bar u_{in}) + E(\bar u_{out})$$
where $\bar{u} _{in} , \bar u_{out}$ are compactifications obtained by extending the maps to the compactified targets $\bar M_{in}, \bar M _{out}$ obtained by symplectic cutting.
\end{defn}

\subsection{Gluing broken things }

In this section, we will discuss how to obtain unbroken maps from broken discs and strip. Given a broken strip or disc $(\mathfrak{u,D})$ in a broken manifold $\mathbb{M}$ with boundary on broken (or unbroken) Lagrangians $\mathbb L$ or $K$, we want to investigate if this broken map can be viewed as a limit of a neck-stretching process. More precisely, by ``gluing" a broken disc $\mathfrak{u}$,  we mean obtaining a family of holomorphic discs $u_\tau$ in a $\tau-$stretched manifold $M^\tau$ such that $u_\tau \to \mathfrak{u}$ as $\tau \to \infty$.  We state some gluing results from \cite{palmerwoodwardsurg} without proof, since there is no change in the proof for the case of broken strips. We refer the reader to the op. cit. article for the case of broken discs. The gluing method in \cite{palmerwoodwardsurg} is similar to gluing discs with Lagrangian boundary conditions, as done in \cite{abouzaidExotic}.

We recall the notation relevant to gluing from \cite{palmerwoodwardsurg}.  Given a nodal domain $\mathfrak{S}$ with $k$ puncture nodes, and gluing parameters $\delta_1,\dots \delta_k > 0$, we define the glued domain $S^{\delta_1,\dots,\delta_k}$ by gluing strip-like $(-L,L) \times \R \times \interval$ or cylindrical necks $(-L,L) \times S^1$ of length $2L =|\ln (\delta_i)|$ at the i$^{th}$ puncture node. Given a $l-$broken manifold $\mathbb M$, we can similarly define $M_{\delta_1,\dots,\delta_l}$ by gluing necks of appropriate length. Note that the notational convention we choose is to use subscript to denote parameters of a glued broken manifold while superscript is used to denote a stretched manifold.

\begin{Theorem}[Theorem 6.8 in \cite{palmerwoodwardsurg}]\label{thm:glu} 
Let $(\mathfrak{u,S})$ be a regular broken strip with boundary on broken Lagrangian $L$ and unbroken Lagrangian $K$ with limiting eigenvalues $\lambda_1,\dots\lambda_k$ of Reeb chords or orbits. Then, there exists  $\delta_0$ such that for each $\delta \in (0,\delta_0)$ there exists an unbroken strip $u_\delta: S^{\delta/\lambda_1,\dots,\delta/\lambda_k} \to M_\delta$ such that $u_\delta \to \mathfrak{u}$ as $\delta \to 0$.
\end{Theorem}

The theorem above shows that for any regular broken strip, we can view it as an SFT-limit of unbroken strips. InSection \ref{sec:reg} we show that we can choose a perturbation scheme on the outer pieces of broken strips so that any broken strip is regular. From now on, until the end of this subsection, we assume that all broken strips are either high index or regular.  The next theorem we state from \cite{palmerwoodwardsurg} shows that once we arrange regularity for all the broken strips, the gluing construction in Theorem \ref{thm:glu} would provide bijection of set of rigid broken disc with rigid unbroken disc in a neck-stretched manifold. 
 
In Section \ref{sec:classif}, from the proof of Lemma \ref{lemm:2level} we see that any regular rigid broken strip has at most two levels. Hence, in the following theorem we consider only two level rigid strips. We use $\mathcal{M}^{<E}_{Br}(\Gamma,\mathbb{M}[l],\mathbb{L}[l],K)_0$ to denote the moduli space of rigid broken strips with compactified energy at most $E$ and constraint of interior markings going to the nice divisor $D$  and $\mathcal M ^{E_{Hor}< E}(M_{\delta},L,K,D)_0$ to denote rigid holomorphic strips with respect to the glued coherent perturbation data with boundary on $L,K$ with  at most $E$ horizontal energy and constraint of interior markings going to $D$. 

\begin{Theorem}[Theorem 6.10 in \cite{palmerwoodwardsurg}] \label{thm:bij}
Let $\mathbb M$ be a broken symplectic manifold with 2 levels and $\mathbb L$  be a broken Lagrangian with 2 levels .  Assume perturbations have been chosen so that every strip in $\mathcal{M}^{<E}_{Br,\mathcal{P}}(\mathbb{M}[l],\mathbb{L}[l],K,D)_0$ is regular. There exists $\delta_0$ such that for $ \delta <\delta_0$ the correspondence $\mathfrak{u} \mapsto u_\delta$ from Theorem \ref{thm:glu} defines a bijection between the moduli spaces $\M^{<E}_{Br,\mathcal{P}}(\mathbb{M}[l],\mathbb{L}[l],K,D)_0$ and $\mathcal M ^{E_{Hor}< E}(M_{\delta},L,K,D)_0$.

$$\M^{<E}_{Br,\mathcal{P}}(\mathbb{M}[l],\mathbb{L}[l],K,D)_0 \leftrightarrow \mathcal M ^{E_{Hor}< E}(M_{\delta},L,K,D)_0$$
\end{Theorem}

\begin{remark}
    The proof of \ref{thm:bij} involves an application of the SFT compactness theorem. Note that we apply SFT Compactness to a sequence of maps with bounded horizontal energy due to Lemma $9.2$ in \cite{bhewz}.
\end{remark}

\noindent If $(L,K)$ is a monotone pair, then from the action-index relation we get that there is an energy bound $E_0$ such that any index-$1$ strip has energy less than $E_0$. Thus, as a corollary of  Theorem \ref{thm:bij} we have the following result.

\begin{cor}\label{corr:rigidregularbijection}
If $(L,K)$ is a monotone pair and a perturbation scheme is chosen that makes every broken strip either high-index or regular and all the elements in $\mathcal{M}^{<E}_{Br,\mathcal{P}}(\mathbb{M}[l],\mathbb{L}[l],K,D)_0$ is regular, then there is a $\delta_0$ such that for $0<\delta<\delta_0$, there is a bijection between the rigid broken holomorphic strips and the unbroken holomorphic strips in $M_\delta$ and all the rigid unbroken strips in $M_\delta$ are regular.

$$\mathcal {M} _{Br}(\mathbb{M}[l],\mathbb{L}[l],K)_0 \leftrightarrow \mathcal {M} (M_{\delta},L,K,D)_0$$

\end{cor}

\begin{proof}
By using $E=E_0$ in Theorem \ref{thm:bij} and the fact that any index $1$ strip has energy less than $E_0$ we directly get the bijection. The fact that the unbroken strips are regular follows from the Floer's Picard Lemma used in the gluing construction of Theorem \ref{thm:glu}, see \cite{palmerwoodwardsurg} for details. 
\end{proof}

%% file: Sections/reg.tex
\section{Regularity} \label{sec:reg}

We prove that under a generic perturbation of complex structure, every broken disc or strip with only puncture nodes (the nodes between two levels) are either regular or high-index. As a consequence of this regularity, we will see in Section \ref{sec:classif} that the only rigid broken strip or rigid broken disc with point constraint has two levels, and the inner level consists of simple discs.

From the construction of moduli space of broken discs and strips, we know that broken maps are fiber products of moduli spaces of levels under the evaluation maps at the punctures. Thus, it is enough to show each level is regular, and the evaluation maps are transverse, we will take this approach to prove regularity. We will prove that the standard complex structure $J_0$ of $\C^n$ is enough to ensure the inner and neck pieces to be transversally cut out. Finally, we use a domain dependent choice of almost complex structure on the outside level to cut the last layer of evaluation fiber product transversally.

\noindent Here is an outline of this section :
\begin{itemize}
    \item In the first subsection, we will focus on the interior pieces. A short exact sequence argument shows that when we forget about the matching condition at Reeb chord/orbits, discs in the interior piece are transversally cut out when they don't touch the critical locus (set of critical points) of the fibration $\pi_m$ (see Equation (\ref{eq:pim}) for definition). As for when the discs touch the critical locus, we show that we can obtain a stronger regularity result. We show that a constrained moduli space of discs that touch the critical locus is transversely cut out.
    \item In the second subsection we will focus on the neck piece and show that all of them are regular and the evaluation map at any of the puncture is a surjective map to the Legendrian torus.
    \item In the third subsection, we will use domain dependent perturbation to regularize the outside component. Notice that using domain dependent perturbations, we can cut out any matching condition coming from the necks transversally. Thus, we can choose domain dependent perturbation to make the evaluation map of the outer pieces transverse to the matching conditions coming from the inside and neck pieces.
    \item In the fourth subsection, we explain how to obtain the fiber product under the matching condition of Reeb chords between the levels transversally cut out.
\end{itemize}

\begin{remark}[Transversality of punctured spheres]
    We leave the details of transversality of punctured spheres with cylindrical ends in the broken maps, since it has already been well studied in the literature. See subsection 7.5 in \cite{floerflip} for details, specifically Lemma 7.20 proves all the required transversality result needed for spheres. Thus, from now on we will deal with transversality of pieces with nonempty boundary.  
\end{remark}

\subsection{Interior piece} \label{ssec:interior}

In this subsection, we show that with the standard  almost complex structure, we can obtain regularity of the moduli space of punctured discs in the interior piece $\C^n$ of the broken manifold $\mathbb M$. We prove a stronger result that the constrained moduli space where the points that touch the critical locus are constrained to be on a particular subspace in $crit(\pi_m)$ is regularly cut out. This shows that the moduli space of discs that touch the critical locus is regular, and of high index.

\subsubsection{When discs don't touch the critical locus}
We show regularity of discs that don't touch the critical locus. Let $u$ be a map from a punctured disc $\Sigma$ with $l$ boundary punctures to $\C^n$ with boundary on the Lagrangian torus $T_\gamma$. If we assume that $u$ never touches a point where $d\pi_m$ is 0, we can obtain the following short exact sequence of 2 – chains where $D_i$ refers to the corresponding linearization of $\delbar$ and the boundary conditions are suppressed from the notation for clarity.

\begin{adjustwidth*}{-5em}{-5em} 
\begin{tikzcd}[sep=0.8em] \label{diag:ses} \centering
            & 0 \arrow[d]                                               & 0 \arrow[d]                                & \arrow[d]0 & \\
0\arrow[r]& {W^{1,p,\delta}(u^*T^v\mathbb{C}^n)}\oplus V^l   \arrow[r,"i"] \arrow[d,"D_1"] &{W^{1,p,\delta}(u^*T\mathbb{C}^n)} \oplus V^l \oplus H^l                  \arrow[r,"d\pi_{m}"] \arrow[d,"D_2"]& {W^{1,p,\delta+\frac{\pi}{n}}((\pi_m\circ u)^*T\mathbb{C})}\oplus H^l                      \arrow[r]\arrow[d,"D_3"]&  0\\
 0\arrow[r]& {W^{0,p,\delta}( \bar{Hom}(T\Sigma, u^*T^v\mathbb{C}^n)}   \arrow[r,"i"]\arrow[d] &{W^{0,p,\delta}( \bar{Hom}(T\Sigma, u^*T\mathbb{C}^n)}                      \arrow[r,"d\pi_{m}"]\arrow[d] & {W^{0,p,\delta+\frac{\pi}{n}}( \bar{Hom}(T\Sigma, (\pi_m \circ u)^*T\mathbb{C})}                        \arrow[r]\arrow[d]&  0\\
  & 0                                               & 0                               & 0 & \\
\end{tikzcd}

\end{adjustwidth*}

\noindent Since the above is a short exact sequence, to prove $D_2$ is surjective, it is enough to prove that $D_1 , D_3$ are surjective. We show that $\wt D_2$, the parameterized version of $D_2$ that takes into account the moduli space of domains $\Gamma_{in}$. We need the domain perturbations to prove $\wt D_3$ is surjective. To prove $D_1$ is surjective, one approach is to directly show that index of $D_1$ is $n-1$ and kernel is of dimension $n-1$. We outline another approach to show regularity using a splitting argument similar to Lemma 3.3.2 in \cite{jholobook}.

\subsubsection{Compactifying along the ends and index computations} 

The index problem and regularity of Cauchy Riemann operators on bundles over discs have been well understood in the literature, see Appendix C of \cite{jholobook}.  We aim to extend our bundles to the closed disc so that we can use the results of appendix C in \cite{jholobook}. We can compactify $u^*T^v\C^n$ to a holomorphic vector bundle over the closed disc such that $\del u^* T^v T_\gamma$ compactifies to a real sub-bundle. Notice that near the ends, $u^*T^v\C^n$ is holomorphically trivial such that on the boundary since it is isomorphic to the bundle $T^v ((\pi_m \of u)^*\C^n)$ where $(\pi_m \of u)^*\C^n$ is the pullback of the fibration $\pi_m : \C^n \to \C$. We can take the pullback of the fibration $\pi_m$ since near the ends $u$ doesn't hit the critical locus as the map $u$ is close to a trivial strip over a Reeb chord or a trivial cylinder over a Reeb orbit. We compactify $u^*T^v\C^n$ by gluing a trivial bundle at the ends, we denote this bundle as $(u^*T^v\C^n)_{comp}$. Similarly, $\pi_m \of u$ extends to a holomorphic map $\ol{\pi_m \of u}$ from the closed disc to $\C\P^1$ with boundary on $\ol{\gamma}$. We can compactify $(\pi_m\of u )^* T\C$ by $(\ol{\pi_m \of u})^*T\C\P^1$. We denote the linearization of the Cauchy-Riemann section $\delbar$ at $\ol{\pi_m \of u}$ as  $D^c_i$. We can check that the restriction of $D^c_i$  to $(\pi_m \of u)^*T\C$ equal to $D_i$ by a direct computation.

We can identify the kernel and cokernel of the compactified problem with the cylindrical problem by adding constraints. Since we are using the standard complex structure, the linearizations $D_i$ are just $\delbar$, we can use Taylor series expansion to identify the kernel and cokernel of the uncompactified linearized problems with the kernel and cokernel of the compactified problem with constraints. See Proposition 5.9 in \cite{palmerwoodwardsurg} for a proof of isomorphism between the (co)kernels of the cylindrical and compactified linearized operator. 
\begin{lemma}\label{lemm:identif}
$\ker D_i \cong \ker D^c_i$ and $\coker D_i \cong \coker D^c_i$ where the compactified problem for $D_2$ with constraint of having sections that vanish at points mapping to $\infty$ at the same order as $ \ol {\pi_m\of u}$.
\end{lemma}

We now use Lemma \ref{lemm:identif} to compute indices. From the short exact sequence, we know that $\Ind(D_2) = \Ind(D_1) + \Ind (D_3)$. The index of $D_1$ is $n-1$ by an index gluing argument applied to the vertical index of the fibration $(\pi_m \of u )^* \C^n$ where we pinch off the singular fibers. Note that for index computation, we can perturb any $u$ smoothly over a compact region while preserving the Fredholm index, thus even if $u$ passes through the critical locus, we can assume that after small perturbation over a compact region, $u$ misses the critical locus. Thus we still get a short exact sequence whose $D_1$ has index $n-1$.  The index of $D^c_2$ without any constraints is given directly by the Riemann-Roch formula as $1 + \mu(\ol{\pi_m \of u})$ where $\mu$ is the Maslov index, and with the constraint of vanishing at the same order as $\ol{\pi_m \of u}$, the index is $1 + \mu(\ol{\pi_m \of u}) - \ol{\pi_m \of u }_{weighted} \inv (\infty)$ where $\ol{\pi_m \of u }_{weighted} \inv (\infty)$ is a weighted sum of intersection of $\pi_m \of u$ with $\infty$ where interior intersections are weighted with $2$ and boundary intersections are weighted with $1$. Thus, we have the following index formula :

\begin{lemma}\label{indexLemm}
Let $D_u $ denote the linearization of the $\delbar$ operator at $u$ in the cylindrical setting, then 
$ \Ind(D_u) = n + \mu(\ol{\pi_m \of u}) - \ol{\pi_m \of u }_{weighted} \inv (\infty).$
\end{lemma}
\begin{remark} \label{rem:indexsinglepunct}

As a quick application of the index formula, we  compute the index of a disc with only one boundary puncture. The index of a disc $u : \mathbb{H} \to (\C^n, T_{\gamma)}$ that is asymptotic to a Reeb chord of length $k\pi/n$ is $n+k$. We can also prove the surjectivity of $D_2$ at this stage by using the identification of kernel and cokernel with compactified problem.
\end{remark}
\begin{Theorem}[Horizontal Transversality] \label{thm:hortran}
$\coker \wt D_3 = \{0\}$.
\end{Theorem}

\begin{proof}
From the identification of the cokernels in Lemma \ref{lemm:identif}, it is enough to check that $\coker \wt D^c_3$ is 0 with the appropriate constraints. Since this is a one dimensional problem, a modification of the automatic transversality result in \cite{seidelFuk} proves transversality in this case. The only modification to Lemma 13.1 and 13.2 in \cite{seidelFuk} which we need for our case is working with marked points in the interior. See the proof of horizontal transversality in Theorem \ref{thm:vctran} below for a detailed explanation.
\end{proof}

\subsubsection{Splitting of holomorphic bundles over disc}

We show that a Cauchy Riemann operator on a holomorphic bundle over a closed disc with  totally real boundary conditions can be split using results of Oh in \cite{riemHilb}. Any holomorphic bundle over the closed disc is holomorphically trivial by applying the Oka-Grauert principle to an extension of the bundle to a slightly larger open disc. Thus, a pair ($\wt E,\wt F$) of holomorphic bundle with real sub-bundle on the boundary of a disc is holomorphically isomorphic to the pair of bundles $(\D \times \C^n, F)$ where $F$ is a real sub-bundle of $S^1\times \C^n$. In \cite{riemHilb}, Oh proves that a linearized Cauchy-Riemann problem with boundary for the standard complex structure splits into a direct sum of 1 dimensional problems (see Theorem 1 in \cite{riemHilb} and Corollary 5.19 in
\cite{palmerwoodwardsurg} for adaptation to the cylindrical setting). Applying Oh's splitting result to the vertical part, we have the following: 

\begin{lemma}[Splitting Lemma]\label{lemm:splitt}
The Cauchy-Riemann problem, $$D^c_1: W^{k,p}(u^*T^v\C^n,\del u^* T^vT_\gamma)_{comp} \to W^{k-1,p}(\ol \Hom (T\ol \D,(u^*T^v\C^n)_{comp}),$$ splits into direct sum of $n-1$ one dimensional problem under a change of trivialization, $$D^c_1:\oplus_{i=1}^{n-1}W^{k,p}(E_i,F_i) \to \oplus_{i=1}^{n-1}W^{k-1,p}(\ol\Hom (T\D,E_i))$$ where $\oplus E_i \cong u^*T^v \C^n$ and $\oplus F_i \cong \del u^* T^vT_\gamma$.
\end{lemma}

We can apply Lemma \ref{lemm:splitt} to obtain the following transversality result.

\begin{Theorem}(Vertical Transversality)\label{thm:verttran}
$\coker D_1 = 0$
\end{Theorem}
\begin{proof}
From the splitting Lemma, we have that $\coker D_1 = \oplus_{i=1}^{n-1} \coker D_1|_i$ where $D_1|_i$ is the restriction of $D_1$ to the bundle pair $(E_i,F_i)$. Thus, it is enough to show that each of the $\coker$ are 0. We also know that for 1 dimensional Cauchy-Riemann problem with boundary condition, either the kernel or the cokernel is non-zero. Thus it is enough to show that $\ker D_1|i$ is non-zero for all $i$.

Note that there is a $T^{n-1}$ action on the space of maps given by $$(\theta_1, \dots \theta_{n-1}). (u_1,\dots,u_n) = (e^{i\theta_1}u_1,e^{i\theta_2}u_2\dots,e^{-i\sum\theta_j}u_n),$$ and this action preserves the projection of the map under the map $\pi_m$. This action restricts to an action on the space of holomorphic maps, thus taking the infinitesimal action at $u$, we get a map $\mathfrak{t}^{n-1} \to \ker D_1 \cong D^c_1$. Also note that if we compose $\mathfrak{t}^{n-1} \to \ker D_1 \cong \ker D^c_1$ with the evaluation map $ev_q$ for a point $q$ on $\del \Sigma$, we have an isomorphism of vector spaces $\mathfrak{t}^{n-1} \xrightarrow{\sim} u^*T^v\C^n|_q $. This proves that $\ker D_1|_i$ is non-zero for all $i$ as otherwise we wouldn't get an isomorphism since $\oplus E_i \cong u^*T^v\C^n$. Thus $\coker D_1 = 0$
\end{proof}

\subsubsection{When discs touch the critical locus}

In this section, we show that discs that touch the critical locus are regular. Moreover, we prove that a constrained moduli space is regularly cut out. When a disc $u$ touches the critical locus of $\pi_m$ at the points $\{\zeta_1,\dots \zeta_i\}$, we lose the vertical bundle which helped us prove regularity.  The vertical bundle still exists over $\D \setminus \{\zeta_1,\dots \zeta_i\} $. We will show that the vertical bundle can be holomorphically extended over the points $\zeta_j$ and then use a splitting argument with domain dependent perturbations to prove transversality of a ``vertically constrained"" problem where we put high codimension constraints to the extension of the vertical bundle. This allows us to show that broken strips that touch the critical locus are high index, hence can't be rigid.

The primary obstruction to defining a vertical bundle when a disc touches the critical locus is that the dimension of the kernel jumps up at the critical points. When $p\in \C^n \setminus \text{crit}(\pi_m)$, the vertical direction is $$T^v_p\mathbb C^n = \ker \langle [p_2\dots p_n \quad p_1p_3\dots p_n \quad p_1\dots p_{n-1}], \rule{0.5cm}{0.15mm}\rangle, $$ where $\langle \;, \, \rangle$ is the complex bilinear dot product. Define the rational map $\phi$ as follows

\begin{align*}
    \phi : \mathbb C^n &\dashrightarrow \mathbb P^{n-1}\\
    (z_1,\dots,z_n)&\mapsto [z_2z_3\dots z_n : z_1z_3\dots z_n:\dots : z_1\dots z_{n-1}].
\end{align*}

\noindent The vertical fiber over a point $p$ is equal to $\ker \langle \phi (p), \rule{0.4cm}{0.15mm} \rangle$ since the kernel of the dot product is invariant under scaling by $\C^*$. We will show that this rational map, $\phi$, can be lifted to a smooth map.

\begin{Theorem}
There exists a complex manifold $\wt \C^n$ obtained by iterated blow up of $\C^n$ along complex subspaces and a smooth map $\wt \phi: \wt \C^n \to \P^{n-1}$ such that  the following diagram commutes,
$$\begin{tikzcd}
\centering
        \wt \C^n \arrow{r}{\wt \phi}  \arrow{d}{\wt \pi} &\P^{n-1} \\
        \C^n \arrow[ur,dashed,"\phi"]
\end{tikzcd}
$$
where $\wt \pi$ is a composition of blowdown maps.
\end{Theorem}

 Now, $\phi$ can be extended to a map $\widetilde \phi : \widetilde{\mathbb{C}} \rightarrow \mathbb P^{n-1}$ by doing an inductive iterated blowup. Now the vertical subbundle of $u^*T\mathbb C^n $  over $z $ is defined by $\ker \widetilde\phi \circ \widetilde u (z)$ where $\widetilde u$ is the lift of the disc to the blowup. Because the disc $u$ maps to trivial cylinders (or strips) over Reeb orbits (or  chords) near the punctures, we can check that there is an extension of $\widetilde \phi \circ \widetilde u $ to the whole disc.

\begin{proof}

The iterated blowup is constructed by first blowing up the origin in $\C^n$ to obtain $\wt \C^n_0 \xrightarrow{\pi_0} \C^n$, then blowing up the lift of the lines $\C\langle e_i \rangle$ in $\C^n$ generated  by the standard basis vector $\{e_1,\dots,e_n\}$ of $\C^n$ to obtain $\wt \C^n _1 \xrightarrow{\pi_1} \wt \C_0^n$. We continue this process of blowing up and stop at the $(n-2)^{th}$ step where we have blown up all (lifts of) the standard subspaces of complex dimension $\leq n-2$ in $\C^n$ to obtain a chain of blowups as follows 
\begin{equation}\label{iterblow}
    \wt \C^n _{n-2} \xrightarrow{\pi_{n-2}} \wt \C_{n-3}^n \xrightarrow{\pi_{n-3}} \dots  \wt \C^n_0 \xrightarrow{\pi_0} \C^n
\end{equation}

\noindent We denote the $n-2$ iterated blown up space as $$\wt \C^n \xrightarrow{\pi_{tot}} \C^n.$$ Note that the map $\phi$ factors as $\psi \of pr$ where $pr:\C^n \setminus \{0\} \to \P^{n-1}$ is the standard projection map to the projective space and 

\begin{align*}
    \psi : \P^{n-1} &\dashrightarrow \P^{n-1} \\
    [z_1 : \dots : z_n]&\mapsto [z_2z_3\dots z_n : z_1z_3\dots z_n:\dots : z_1\dots z_{n-1} ].
\end{align*}

Let $\wt \P^{n-1}$ be the space obtained by a similar construction as in (\ref{iterblow}), $(n-2)$ times iterated blowups of $\P^{n-1}$ along the standard coordinate projective subspaces. We can check that $\wt \C^n$ is a $\C^*$ bundle over $\wt \P^{n-1}$. To get a lift of $\phi$, it is enough to show that the rational map $\psi$ lifts to a smooth map $\wt \psi$ on $\wt \P^{n-1}$. Indeed, we can define $$\wt \phi (z)  = \wt \psi (\wt{pr}  (z))$$ where $\wt{pr}$ is the projection map from $\wt \C^n \to \wt \P^{n-1}$.

We use an inductive process to lift the map. When the dimension is two, we can define the map $\wt \psi$ to be just $\psi$. Assume $\wt \psi$ has been defined for dimension $n$ and thus $\wt \phi$ is defined for dimension $n+1$. We can define $\wt \psi$ for dimension $n+1$ by taking  a standard $\C^{n+1}$ chart near the vertex of the moment polytope of $\P^n$ and extending the map $\psi$ to $\wt \C^{n+1}$ by $\wt\phi$ that has already been defined. For example, take the chart $\C^{n+1} \to \P^{n+1}$ given by $z\mapsto [1:z]$, then we define, under the identification of the domain using this chart,  $$\wt\psi ( z) := [0:\wt\phi(z)]$$ on the exceptional divisors of $\wt \C^{n+1}$. We can do this for all the vertices of the moment polytope of $\P^{n+1}$ and check that we get a well-defined map $\wt \psi : \wt \P^{n+1} \to \P^{n+1} $ which is a smooth lift of the rational map $\psi$.

\end{proof}

\begin{cor}
Let $u:\Sigma \to \C^n$ be a holomorphic punctured disc with boundary on a torus segment $T_\gamma$, such that $u\inv(crit(\pi_m))=\{\zeta_1,\dots,\zeta_l\}$. Then the pullback bundle $u^*T\C^n$ splits as 
\begin{equation*}
    u^*T\C^n  = \wt V \oplus \wt H
\end{equation*}
where $\wt V$ is a subbundle that matches with the vertical bundle $ V = \ker d\pi_m\big|_{\Sigma\setminus \{\zeta_i\}}$ on away from the critical points $\Sigma\setminus \{\zeta_1, \dots \zeta_l\}  $.
\end{cor}
\begin{proof}

We use the lifted map to obtain extension of vertical sub-bundles over the critical points in a disc.  From the lifting property of blowup, the holomorphic map $u$ to $\C^n$ extends to $\wt u$, a holomorphic map to the blowup $\wt \C^n$. We have a map $\wt \phi \of \wt u$ to $\P^{n-1}$, thus we can define a holomorphic sub bundle $\wt V$ of $u^*T\C^n$ by defining $\wt V_z = \ker \langle \wt \phi \of \wt u (z),\sbl \rangle $. We can then define a quotient line bundle $u^*T\C^n/\wt V$ over the disc after compactifying the bundles, and since any short exact sequence of holomorphic bundles over a closed disc split, we get a splitting $$ u^*T\C^n = \wt V \oplus \wt H$$ where $\wt H$ denotes the quotient bundle.

\end{proof}

The previous corollary lets us split the Cauchy-Riemann operator. Denote the split operator as follows, $$D_u= D_V \oplus D_H.$$

\noindent We now compute the index of $D_V$. For simplicity of notation, let $u$ be a map from a punctured disc that touch the critical locus $\text{crit}\pi_m$ only at $0$ and assume  $$u(0) = (0,0,\dots,p_{k+1},\dots,p_n),$$where $p_i \neq 0 $. Thus, we have, $$u_i(z)=z^{m_i} + \text{higher order terms}$$ for some positive integers $m_i$, we can assume without loss of generality that $m_i\leq m_{i+1}$. We will prove that the index of $D_V$ is $n-1+ 2(\sum_{i=1}^{k-1} m_i)$. Consider the holomorphic sections, $$X_i(z) = (0,\dots,u_i(z),\dots,-u_k(z),\dots,0),$$ then the sections $X_1,\dots, X_{k-1},X_{k+1},\dots,X_n$ forms a frame of $\wt V$ over $\D^2\setminus\{0\}$ and the boundary condition of the bundle pair $(\wt V,\del u^* T^v T_\gamma)$ is just $i\R^{n-1}$. Note that $X_{k+1}(0),\dots,X_{n}(0)$ are linearly independent and the $X_i(z)$ vanishes at the order $m_i$ for $i<k$. Thus these sections give an isomorphism of $\wt V$ near $0$ with the twisted trivial bundle $\oplus_{i=1}^{k-1} O(m_i)$, thus from a gluing argument we have $$\Ind(D_V) = n-1 + 2c_1(\oplus_{i=1}^{k-1} O(m_i)) =n-1 + 2(\sum_{i=1}^{k-1} m_i). $$

\noindent  Assume the map $u$ that touches the critical locus at $(\zeta_1,\dots,\zeta_i)$ and $k_l$ is the number of 0's in $u(\zeta_l)$ for $1\leq l \leq i$. Let $m_{l_1} \leq m_{l_2} \leq \dots m_{l_{k_l}}$ be the order of the vanishing of the coordinates of $u$ at $\zeta_l$. We define the critical multiplicity $m_c$ of the map $u$ as $$m_c(u)=\sum_{l=1}^i\sum_{q=1}^{k_l-1} m_{l_q}.$$

\begin{lemma}\label{lemm:vertcrit}
If $u$ touches the critical points of $\pi_m$, the index of the vertical linearized operator $D_V$ is  $n-1 + 2m_c(u)$ and the horizontal index is $1 + \mu(\ol{\pi_m \of u}) - \ol{\pi_m \of u}_{weighted}\inv (\infty) - 2m_c(u) $

$$\Ind (D_V) = n-1 + 2m_c(u)$$
$$\Ind (D_H) = 1 + \mu(\ol{\pi_m \of u}) - \ol{\pi_m \of u}_{weighted}\inv (\infty) - 2m_c(u)$$
\end{lemma}
\begin{proof}
    By using an index-gluing argument we can generalize the discussion to obtain the required index formula.
\end{proof}

We will now relate the kernel and cokernel of the linearized operators with a linearized operator of a constrained problem on. 

\begin{lemma}\label{lemm:horiden}
$\ker \wt D_H \cong \ker \wt D^{const}_{\pi_m \of u} $ and 
$\coker \wt D_H \cong \coker \wt D^{const}_{\pi_m \of u} $
\end{lemma}
\begin{proof}
Let $D^{const}_{\pi_m \of u}$ be the linearization of $\delbar$ of the constrained Cauchy Riemann problem at $\pi_m \of u$ with the constraint of vanishing at $\zeta_l$ to the order $\sum_{q=1}^{k_l-1} m_q$ where $m_1 \leq m_2 \leq \dots m_{k_l}$ are the orders of vanishing of the coordinates of $u$ at $\zeta_l$.
The push-forward $d\pi_m$, maps $\ker D_H$ isomorphically to $\ker D^{const}_{\pi_m \of u}$. This follows from the fact that $d\pi_m$ has this order of vanishing at $\zeta_l$. Since $d\pi_m$ is isomorphism on fibers of $\wt H$ except over $\zeta_i$, we see that the map $\ker D_H \xrightarrow{d\pi_m} \ker D^{const}_{\pi_m \of u}$ is injective. By checking the order of vanishing at $\zeta_i$ we see that this map is also surjective. We also know that $\Ind (D_H) = \Ind (D^{const}_{\pi_m \of u})$, thus, their cokernels are also isomorphic. We can use the fact that since $d\pi_m$ is complex linear and the linearizations $D$ we deal with are $\delbar$ under holomorphic frames to see that the kernels of the parameterized linear operators $\wt  D_H$ and $\wt D^{const}_{\pi_m \of u}$ are isomorphic, and the index argument shows that the cokernels are isomorphic. 
\end{proof}

\begin{defn}[Level Structure]

We can define a \textit{level structure} on the set crit$(\pi_m)$ by defining a point $(p_1,\dots,p_n) \in \text{crit}(\pi_m)$ to be of level $k$ if exactly $k$ of the $p_i$'s are 0. Thus, every critical point has a level ranging from $2$ to $n$.
\end{defn}

We will now define the notion of being \textit{vertically constrained} at the critical points and the notion of point constrained.

\begin{defn}[Vertically Constrained moduli space]
Let $u$ be a holomorphic map from a punctured disc to $\C^n$ with boundary on $T_\gamma$ that meets the critical locus $crit(\pi_m)$ at the points $\{ \zeta_1 , \dots, \zeta_i \}$. Assume that the \textit{level} of $u(\zeta_j)$ is $k_j$. The family of punctured holomorphic discs with the constraint of touching a critical point of level $k_i$ at the point $\zeta_i$ is called the vertically constrained moduli space.
\end{defn} 
\noindent We denote the linearization of $\delbar$ at $u$ of the vertically constrained problem by $D_u^{vc}$. By using the splitting of $u^*T\C^n$ into $\wt V \oplus \wt H$ as before, we get a splitting of $D_u^{vc} = D_V^{vc} \oplus D_H^{vc}$, where $D_H^{vc} = D_H$ since the constraints are only on the subbundle $\wt V$. We denote the linearized operator including domain parameterization by $\wt D_u^{vc}$. From Lemma \ref{lemm:vertcrit} we have the following, $$\Ind (D^{vc}_u) = \Ind (D_u) - 2\sum_{j=1}^i{k_j}.$$ If a disc doesn't meet the critical locus, then by definition, being vertically constrained would impose no constraints on the linearized Cauchy-Riemann problem. We denote the moduli space of vertically constrained discs in $\C^n$ with boundary on $T_\gamma$ as $\mathcal{M}^{vc}_{in}$.  

Recall that the disc potential is obtained by count of rigid discs with point constraints. We define the notion of \textit{point constrained} moduli space of broken discs. The notion of point constrained interior moduli space is defined similar to the notion of broken disc with point constraint. Given a point $p$ in $L_{in}$, the moduli space of discs with constraint at $p$ is defined to be the subspace of the moduli of discs that pass through $p$.  The moduli space with point constraint can be realized as $ev_{q}\inv(p)$ where $q$ is a point on the boundary of a disc and $ev_q$ is the evaluation map at $q$.

\begin{Theorem}\label{thm:vctran}
For any disc $u$ touching  crit $(\pi_m)$, the parameterized vertically constrained moduli space of discs with point constraint preserving the level structure of the critical points is regularly cut out, i.e. $\wt D_u^{vc}$ is a surjective map.
\end{Theorem}

\begin{proof}
First, we split the operator $\wt D_u^{vc} = D_V^{vc} \oplus \wt D_H$. Since the action of $T^{n-1}$ as defined in Theorem \ref{thm:verttran} on the space of maps preserves the level structure of the critical points, the action descends to the vertically constrained problem and gives transversality of $D_V^{vc}$. From Lemma \ref{lemm:horiden}, we see that it is enough to prove that linearization of $\delbar$ at $\pi_m \of u $ with the constraint of vanishing at the critical points $\zeta_i$ is surjective. From the proof of Lemma 4.3.2 of \cite{JholoWendl} we see that for any smooth vector field $X$ that vanishes at the points $\zeta_i$s, we have that $d(\pi_m \of u) X$ is a smooth vector field that satisfies the constraint of vanishing at $\zeta_i$ to the correct order since $d(\pi_m \of u)$ has order of vanishing $1$ less than that of $\pi \of u$ and $X$ has an order of vanishing at least $1$ at the points $\zeta_i$. Thus, we can apply the result Lemma 4.3.2 in our case to conclude that image of $ D^{const}_{\pi_m \of u}$ doesn't change if we expand our domain to $W^{k,p} (\ol {\text{End}} (T\Sigma)) \oplus W^{k,p,\delta + \frac{\pi}{n}}_{const}((\pi_m \of u )^*T\C,\del(\pi_m \of u )^*T\gamma)$. Now since $\pi_m \of u$ is a non-constant map to one dimensional manifold, we can use the argument of automatic transversality as done in Lemma 13.2 in \cite{seidelFuk} to get that $\wt D^{const}_H$ is surjective. We finally use the vertical and horizontal transversality to show that $p$ is a regular value of the evaluation map $ev_q$ to prove the regularity of the point constraint.
\end{proof}

\begin{cor} \label{corr:submatch}
The moduli space $\mathcal{M}_{in}$ of vertically constrained discs in $\C^n$ with boundary on $T_\gamma$ is a smooth manifold and the evaluation map $ev_p$ at any boundary puncture $p$ that takes a map $u$ to the starting point of the Reeb chord $u$ is asymptotic to at the puncture $p$, 
$$ev_p : \mathcal{M}^{vc}_{in} \to T_\pm,$$ is a submersion to the Legendrian torus.
\end{cor}

\begin{proof}
Theorem \ref{thm:vctran} shows that $\mathcal{M}^{vc}_{in}$ is a manifold. The submersion to the Legendrian torus follows from the $T^{n-1}$ action. 
\end{proof}

\subsection{Neck Pieces}

In this subsection, we show that regularity of maps in the neck piece can be obtained without perturbing the standard complex structure. The treatment of the neck piece is very similar to the treatment of the interior piece, albeit a bit easier.

We use a similar splitting trick as in Subsection \ref{ssec:interior} to split the Cauchy-Riemann operator. Let $u$ be a disc mapping to the neck piece $\C^n\setminus\{0\}$ with boundary on $T_{\R^*}$. Using \cite{riemHilb}, we get that $(u^*T\C^n,\del u^* TT_{\R^*})$ splits into 1 dimensional bundle pairs $\oplus_{i=1}^n(E_i,F_i)$. We can show regularity of $D_u$ by showing that the kernel of the restriction, $\ker D|_i$ is non-zero for all $i$ since these are one dimensional problems.  
\begin{Theorem}\label{thm:necktranss}
If $u$ is map to the neck-piece with boundary on the cylindrical Lagrangian $T_{\R^*}$, then $\coker D_u = \{0\}$
\end{Theorem}
\begin{proof}
 In the neck piece, there is an action of $\R \times T^{n-1}$ where the $\R$ action is by translation and $T^{n-1}$ action is same as in the interior piece. The infinitesimal action provides a map $\R \times \mathfrak{t}^{n-1} \to \ker D_u$ such that under the evaluation map at a boundary point $q$, we have an isomorphism $\R \times \mathfrak{t}^{n-1} \xrightarrow{\sim} T_{u(q)}T_{\R^*} $. A similar argument as in Theorem \ref{thm:verttran} gives us $\coker D_u=\{0\}$.
\end{proof}

\noindent Let $p$ be a boundary puncture. Let $ev_p$ denote the evaluation map at $p$ that takes values in the space of Reeb chords. Specifically,  $ev_p(u)$ is the starting point of the Reeb chord the map $u$ is asymptotic to at the puncture $p$. As a corollary of Theorem \ref{thm:necktranss}, we have the following result,
\begin{cor}\label{corr:neckmatch}
The moduli space $\mathcal{M}_{neck}$ of discs in the neck piece is a smooth manifold and the evaluation map $ev_p$ \begin{equation*}
ev_p : \mathcal{M}_{neck} \to T_\pm,    
\end{equation*} is a submersion to the Legendrian torus.
\end{cor}
\begin{proof}
Same argument as in proof of Corollary \ref{corr:submatch}.
\end{proof}

\subsection{Outer Piece} \label{subsec:outtran}

We will now show that we can choose a domain dependent perturbation data on broken strips and broken discs to form a coherent perturbation data that attains our goal of showing every broken strip is either regular or of high index. For any combinatorial type $\Gamma^{st}$ of stable connected nodal strip with $k$ markings, the smooth components of $\Gamma^{st}$ can be divided into the following collections 
\begin{itemize}
    \item strip piece,
    \item discs without interior marked points,
    \item discs with interior markings.
\end{itemize}
\noindent  Given a stabilized combinatorial type $\Gamma^{st}$ of broken disc with $k$ intersections with the nice divisor $D$, we will show that we can choose generic domain dependent complex structure to obtain regularity for all outside pieces along with regularity of the matching condition at the contact $S^{2n-1}$ and Legendrian $T_\pm$. The matching conditions between the last neck piece and outside piece is given by countably many images of smooth maps coming the interior pieces given by evaluating at the positive puncture of last level. For a broken strip with $l$ neck pieces, the smooth map it contributes is   $$ev_{l+} : \mathcal{M}^{vc}_{in} \times_{ev_{in},ev_{1^-}} \mathcal{M}_1 \times_{ev_{1^+},ev_{2^-}} \dots \times_{ev_{l-1^+},ev_{l^-}}\mathcal{M}_l.$$

\noindent The choice will be made such that it satisfies the following properties on each type of smooth component.

\begin{itemize}
    \item \textbf{(P1)} On strips, $J_k$ will be only $t$ dependent and equal to $J_{out}$ at the boundaries .
    \item  \textbf{(P2)} On discs without interior marked points, $J_k$ will be the constant almost complex structure $J_{out}$.
    \item \textbf{(P3)} On discs with interior marked points, $J_k$ is a domain dependent almost complex structure that is equal to $J_{out}$ on the boundary .
\end{itemize}

\noindent Once we set such a perturbation data $\mathcal{P}$, given any combinatorial type $\Gamma$ of a broken strip or a broken disc, we define the notion of domain dependent broken holomorphic strip by using the standard complex structure on the inside and neck pieces and the domain dependent complex structure on $\Gamma^{st}$.

\begin{Theorem}[Outside Transversality] \label{thm:outtrans}
There is a coherent family of perturbation data $\mathcal{P}$ satisfying the properties (P1-3) which makes moduli space of strips $\M^{dd}_k(\widehat M_{out},L_{out},K,D)$ or the moduli space of discs with $\M^{dd}_k(\widehat M_{out},L_{out},D)$ regular and the evaluation map at the cylindrical punctures transversely intersecting with the matching condition coming from the inside and neck pieces.
\end{Theorem}

\begin{proof}
For satisfying the coherence conditions we inductively choose the perturbations based on the number of interior marked points. After choosing perturbations for at most $k-1$ markings, we use the perturbations of $\M_{<k}$ to define perturbation on lower dimensional strata of $\M_k$. Since the space of almost complex structures we care about is contractible, we can extend the perturbations from the lower dimensional strata of $\M_k$ to the whole of $\M_k$.

For $k=0$ markings, there is no stable combinatorial type of strip. The only possible broken type $\Gamma$ without any outside interior marking mapping to $D$ is the case when the outside piece is just a strip with punctures going to Reeb chords or orbits. Let $S$ be a strip with punctures and $u:S \to \widehat M_{out}$ with boundary on $L_{sing}$ and $K$ where the $u$ is asymptotically a trivial strip or cylinder over a Reeb chord or orbit, a similar proof as Theorem 4.3 in \cite{FlHStrans}  shows that the set $R(u)$ of regular points that are somewhere injective in the $s-$direction (notation from \cite{FlHStrans} ) is open dense in the strip $S$. Thus, by constructing a universal moduli space and usual generic transversality arguments we can see that for a generic $t$ dependent family of almost complex structure, the moduli space $\M_t(\widehat M_{out})$ of $J_t$ holomorphic punctured strips in $\widehat M_{out}$ with boundary on $L_{out}$ and $K$ is a smooth manifold and the evaluation maps at the punctures are transverse to countably many matching conditions. The important fact of generic transversality we use here is that the set of regular $J_t$ families  for which each evaluation map is transverse forms a comeagre/Baire set hence countable intersection of  such comeagre sets is also comeagre.

For $k=1$, we first define a domain dependent almost complex structure on a punctured disc with 1 marking where the punctures go to Reeb orbits or chords so that the holomorphic punctured discs are regular and the evaluation at the punctures are transverse to the evaluation from the inside and neck pieces and the evaluation at the interior marking is transverse to $D$. This is obtained by the usual technique of creating a universal moduli space of pairs $(J_k,u)$  of domain-dependent almost complex structure and holomorphic map, showing that it is a Banach manifold and then using Sard-Smale. See proof of Theorem 4.1 in \cite{charwoodstabil}  The boundary strata of a strip with $1$ marked point have strip breaking due to the marking escaping to infinity or disc bubbling when the interior marking goes to the boundary. The domain from the strip breaking will not have any holomorphic maps from it since the ends of the strip go to $L_{out}\cap K$ and the interior marking goes to $D$ but $D\cap (L_{out} \cup K) = \emptyset $.  By choosing domain-dependent almost complex structure on the punctured disc with $1$ marking and a $t$ dependent $J_t$ on the strip without any markings, we force the choice of domain dependent almost complex structure at the boundary strata obtained by disc bubbling. 

In the case of $k=2$ we first encounter a situation where the lower dimensional strata will have a disc without any marking. This combinatorial type can only occur as a domain of a holomorphic strip when the disc without marking is in the neck piece or the map is constant on the disc, as $L_{out}$ is exact in $\M_{out} \setminus D$. In either of those cases, we know that such discs are $J_{neck}$ or $J_out$ regular. We continue this method of choosing domain dependent perturbation inductively by extending the pre-determined choice of almost complex structure on lower dimensional strata of $\M_k$ to choose a coherent perturbation data which makes the moduli space of outer piece regular such that the matching condition with the inside pieces are transversely cut out.
\end{proof}

\subsection{Moduli space of broken discs and strips}

At this stage, we know that the moduli space of maps to the blow-up of interior, neck and outer piece are regularly cut out and the evaluation map at the punctures are submersions when the complex structure on the blown-up interior and neck piece is the standard complex structure and the complex structure on the outer piece is a domain-dependent perturbation, $J_{\Gamma_{out}}$, of the complex structure $J_{out}$. We will show that for such a choice of complex structure, any holomorphic broken disc that avoids the critical locus in the interior is regular in the sense of Definition \ref{defn:brokendiscreg}. 

We define the notion of a \textit{vertically constrained} broken strip by constraining the inner piece of the broken strip to be vertically constrained. We denote the moduli space of such broken strips by $$\mathcal{M}^{vc}_{Br}(\Gamma,\mathbb{M}[l],\mathbb{L}[l],K).$$

Similar to strips, we can add vertical constraints to broken discs with point constraints. We denote the moduli space of broken discs with point and vertical constraints by $$\M^{vc}_{Br}(\Gamma,\mathbb M [l], \mathbb{L}[l],p).$$

We will show that for any combinatorial type $\Gamma = (\Gamma_{in},\Gamma_1,\dots,\Gamma_{out})$, any ($J_{std},\dots,J_{\Gamma_{out}}$)-holomorphic broken strip $\mathfrak{(u,D)}$ is regularly cut out as a fiber product. We illustrate an inductive method which ensures that $$\mathcal{M}^{vc}_{in} \times_{ev_{in},ev_1} \mathcal{M}_1 \times_{ev_1,ev_2} \dots \times_{ev_{l-1},ev_l}\mathcal{M}_l $$ is a smooth fiber product.

\begin{lemma}\label{lemm:brokenpretrans}
There is a coherent perturbation datum $\mathcal{P}$ such that every $l-$ broken strip $(\mathfrak{u,D}) \in \mathcal{M}_{Br}(\Gamma,\mathbb{M}[l],\mathbb{L}[l],K) $ or broken disc $(\mathfrak{u,D})\in \M^{vc}_{Br}(\Gamma,\mathbb M [l], \mathbb{L}[l],p)$ where there are no multiply covered pieces and the only nodes are puncture nodes that map asymptotically to Reeb orbit or chord is regular where the standard complex structure of $\C^n$ is chosen in the neck pieces and inner pieces.
\end{lemma}

\begin{proof}
From Corollaries \ref{corr:submatch}, \ref{corr:neckmatch} and Theorem \ref{thm:outtrans}, we get that for any  holomorphic $(u_{in},u_1,\dots,u_{l})$ satisfying the matching condition, it is regular in the moduli space of (vertically constrained) punctured discs to the inside and neck pieces respectively.  The corollaries also show that the evaluation map at a puncture is a submersion at the matching condition, thus by an inductive argument we see that $(u_{in},u_1,\dots, u_l)$  is a regular broken strip with an inner piece and $l$ neck pieces since the transversality of the evaluation maps at the puncture nodes among these levels can be ensured by using the submersion of the evaluation map at a single puncture as illustrated above. Thus $\mathcal{M}^{vc}_{in} \times_{ev_{in},ev_1} \mathcal{M}_1 \times_{ev_1,ev_2} \dots \times_{ev_{l-1},ev_l}\mathcal{M}_l $ is a smooth manifold and we have an evaluation map $ev_{medium}$ from it to $T_\pm \hookrightarrow S^{2n-1}$, we use a domain dependent perturbation on $\D_{out}$ as shown in Theorem \ref{thm:outtrans} to achieve regularity of $\mathcal{M}_{out}$ and to ensure the evaluation map $ev_{out}$ is transversely intersecting with $ev_{medium}$, hence the total fiber product that produces the moduli of broken strips is transversely cut out. 
\end{proof}

\noindent Using Lemma \ref{lemm:brokenpretrans}, we can now prove the main transversality theorem of this section.

\begin{Theorem}\label{thm:ultimatetransversality}
There is a coherent perturbation datum such that every broken strip or a broken disc with point constraint is either regular or high-index.
\end{Theorem}

\begin{proof}
Once we have attained regularity for broken strips with only puncture nodes, we can conclude that the only possible non-regular broken strips are of high index by using monotonicity of $M$ and $L,K$. Given a broken strip, after replacing multiple covered components by the underlying simple and forgetting nodal components to obtain a single nodal strip with only puncture nodes, we get a holomorphic broken strip which is regular by Lemma \ref{lemm:brokenpretrans} hence non-negative index. Removing non-constant nodal components reduces the total energy, by an index gluing argument we see that monotonicity implies that removing nodal components reduces the index. Replacing multiple covers by single covers also cause the index to decrease. Thus, we have attained our goal of showing that either the broken strip is of high index or is regular.
\end{proof}

%% file: Sections/classif.tex
\section{Classifying Rigid Broken Maps}\label{sec:classif}
The aim of this section is to classify combinatorial types of rigid broken strips or rigid broken discs with the choice of perturbation scheme as done in Section \ref{sec:reg}. We will show that any rigid broken strip has only two levels and all the inside pieces are either punctured spheres or  discs with one boundary puncture that is asymptotic to a Reeb chord of length $\pi/n$. 

\begin{lemma}\label{lemm:2level}
If $(\mathfrak{u,D})$ is a regular vertically constrained broken strip  that is also rigid, it has at most two levels and the inside piece does not intersect the critical locus $\text{crit}(\pi_m)$. Similarly, if  $(\mathfrak{u,D})$ is a regular vertically constrained broken disc with point constraint, it has  at most two levels and the inside piece does not intersect the critical locus $\text{crit}(\pi_m)$.
\end{lemma}

\begin{proof}
We have already chosen a perturbation scheme on the outside piece to make the moduli space of broken strips transversally cut out. If $(\mathfrak{u,D})$ is
a rigid broken strip, the only deformations of the broken map $\mathfrak{u}$ are the ones obtained by the action of the automorphism group of the domain, i.e., by reparameterizing the broken strip. If $\mathfrak{u}$ has more than two levels, there will be a non-trivial neck piece $u_i$, (recall that any neck level should have a non-trivial map from the definition of a broken strip). Since $u_i$ is a non-trivial map,  translation of $u_i$ by the $\R $ action on $\R \times S^{2n-1}$ produces neck maps which have the same asymptotic Reeb chords and orbits, so there is a one dimensional deformation of  $\mathfrak{u}$ obtained by replacing $u_i$ by its translates. Notice that this deformation cannot be obtained by a reparameterization, as the map $u_i$ is not a trivial strip or a cylinder. Thus, we get a contradiction to the assumption that $\mathfrak{u}$ is rigid.

\noindent Now assume that $u_{in}$ touches the critical locus, since it is regular as a vertically constrained broken strip, we have that $\Ind(\wt D^{vc}_{\mathfrak{u}}) < \Ind(\wt D_{\mathfrak{u}}) $, we see that $\mathfrak{u}$ is of high index and thus cannot be rigid.
\end{proof}

\subsection{Elementary discs}

We will now classify the moduli space of single punctured discs with boundary on an asymptotically cylindrical Lagrangian manifold. The heuristic of classifying elementary discs by showing they are sections over the fibration is inspired from  \cite{seidelLongEx}. Recall, we can define a torus segment over any simple path in the punctured complex plane.   Thus, we have  a Lagrangian  $T_{\R+i\e} $ for the path $\gamma: \R \to \C^*$ given by $\gamma(s)=s+i\e$.  $$T_{\R+i\e} = \{(z_1,\dots,z_n)| \pi_m(z)\in \R+i\e, |z_1| = \dots  = |z_n|  \}.$$ Note that this Lagrangian is asymptotically close to the cylindrical Lagrangian $T_{\R \setminus \{0\}}$.

\begin{defn}[Elementary disc in $\C^n$ ]
A map $u: \mathbb{H} \to \C^n$ such that the restriction $u|_\R$ lies on the Lagrangian $T_{\R+i\e}$ or $T_\gamma$ for a cylindrical $\gamma$ and $u$ is asymptotic to a trivial strip over a Reeb chord of length $\frac{\pi}{n}$ is called a elementary disc. We denote the set of elementary discs with Reeb chord $\zeta$ at infinity with boundary on $T_{\R+i\e}$ and $T_\gamma$ by $\mathcal{M}_\zeta(\R+i\e)$, $\mathcal M_\zeta{(\gamma)}$ respectively.
\end{defn} 

In this subsection, we will classify single punctured discs with boundary on $T_{\R+i\e}$ that have a minimal length Reeb chord at its puncture. After classifying such discs, we will use a Gromov-compactness and regularity argument to conclude a classification of elementary discs provide us a classification of discs that occur as interior pieces of rigid broken strips and discs.

\begin{lemma}\label{lemm:simpleissection}
 An elementary disc $u$ with boundary on $T_{\R+i\e}$ can be viewed as sections of $\pi_m$ over the upper half-space $\mathbb{H}+i\e$ or the lower half-space   $-\mathbb{H}+i\e$ i.e. $\exists \phi_u$, a biholomorphism from $\mathbb{H} \to \pm\mathbb{H}+i\e$ and a section $s_u:\pm\mathbb{H}+i\e \to \C^ne$ such that $u= s_u \of \phi_u $.

\end{lemma}

\begin{proof}
It is enough to show that $\phi_u = \pi_m \of u$  is a biholomorphism to $\pm\mathbb{H}+i\e$, then $s_u= u \of \phi_u\inv$ will be the required section. As the map $u=(u_1,\dots,u_n)$ is asymptotic to a trivial strip from Lemma \ref{remark:asymptoticLagrangian}, let's denote the trivial strip by  $$u_\zeta (z)= (z^{1/n}l_1,\dots,z^{1/n}l_n),$$ over a Reeb chord $$\zeta(t)=\frac{1}{\sqrt{n}|l_1|}(e^{it} l_1,\dots,e^{it}l_n ),$$ of length $\frac{\pi}{n}$ in $S^{2n-1}$, where we have $|l_i|=|l_j| $ for all $i,j$ since the boundary of $u$ lies on $T_{\R + i\e}$. The idea of the proof is to show that since asymptotically all the coordinates $u_i$ are away from 0, we get that the product $u_1.u_2\dots.u_n$ is  asymptotic to a trivial strip over a Reeb chord of length $\pi$ and thus essentially, a linear map.

Note that $u$ can be extended to a map from the closed disc to $\P^n$  such that $$u(\infty) = [l_1: \dots:l_n:0] \in T_{Cliff} \sub \P^{n-1}_\infty, $$ where $\P^{n-1}_\infty =[z_0:\dots:z_{n-1}:0]$. The fact that $u(\infty)$  is a point in the torus $T_{Cliff}$ follows from the fact that the boundary of $u$ lies on $T_{Cliff}$ from the definition of $T_{\R +i\e}$.

By checking the $S^{2n-1}$ component of the metric, we know that in limit $d_{cyl} (u(z),u_\zeta(z)) \to 0$ as $|z| \to \infty$.  We can arrange for  $\frac{u(z)}{|u(z)|}$ to lie in a small neighborhood  $\mathcal{N}$ of the Reeb chord $\zeta$ in the unit sphere $S^{2n-1}$, such that $\pi_m (\mathcal{N}) $ lies in a small neighborhood $\mathcal{N}_\zeta$ around $\pi_m(\zeta )$ for large $|z|$. We can also assume that the distance from origin to $\mathcal{N}_\zeta$ has a lower bound $\e_0$ for some fixed $\e_0 > 0$ . This implies that we can extend $\pi_m \of u$ to  a map from the closed disc to $\P^1$ by defining $\pi_m \of u (\infty) =  [1:0]$.

Since the boundary of $\pi_m \of u$ lies on $\mathbb{R} + i\e $ and it extends to a map to $\P^1$, after applying the Schwarz reflection principle and classification of maps from  the disc $\D$ to $\P^1$ with boundary on $\R$ we see that $$\pi_m \of u = p_\R (z) + i\e,$$ where $p_\R$ is a polynomial with real coefficients. Now, using the fact that $d_{cyl}(u(z),u_\zeta(z))$ goes to 0, by checking the $\R$ component of the metric, we have

$$\lim_{|z|\to \infty} \frac{|u(z)|}{|z^{1/n}|\|(l_1,\dots,l_n)\|} = 1.$$

Thus, we can conclude  that 
$$\lim_{|z|\to \infty} \frac{|\pi_m \of u(z)|}{|z||l_1.\dots.l_n|} = 1$$
\noindent Thus combining the fact  $\pi_m \of u = p_\R + i\e$ with the limit above, we see that $p_\R$ is a non-constant linear polynomial and thus a biholomorphism. If $p_\R = az +b$ with $a>0$ we get $\pi_m \of u$ is a biholomorphism to the upper half space $\mathbb{H} + i\e$ , or else if $a<0$, it is a biholomorphism to the lower half-space $-\mathbb{H} + i\e$.

\end{proof}

\begin{remark} \label{remark:onlyupperisenough}
    Using Lemma \ref{lemm:simpleissection} we can reduce the problem of classifying elementary discs to the problem of classifying sections. As  we have $\pi_m(z) = -\pi_m (z\xi)$  for an $n\th$ root of unity $\xi$, a section $s_u$ over $-\mathbb{H} + i\e$ is given by a section over $\mathbb{H}-i\e$ and then multiplying it with $\xi$. Thus, it is enough to classify sections over $\mathbb{H}\pm i\e$ for $\e > 0$. Also note that the \textit{phase action} of $T^{n-1}$ on $\C^n$ defined by $(\theta_1,\dots \theta_{n-1}).(z_1,\dots,z_n)= (e^{i\theta_1}z_1,\dots,e^{i\theta_{n-1}}z_{n-1}, e^{-i(\sum \theta_j)}z_n)$ induces a free action on the set of sections of $\pi_m$ over $\mathbb{H}\pm i\e$. 
\end{remark}

\begin{Theorem}[Classification of sections]\label{thm:classif}
Let $\e > 0$.
\begin{enumerate}

    \item There is a unique $T^{n-1}$ orbit of sections with finite Hofer energy over $\mathbb{H}+i\e$ with boundary on $\R+i\e$ given by $s_\theta = \theta.(z^{1/n}, \dots , z^{1/n} )$ for $\theta \in T^{n-1}$.
    \item There are $n$ disjoint $T^{n-1}$ orbits of sections with finite Hofer energy over $\mathbb{H}-i\e$ with boundary on $\R-i\e$ given by $$s^k_\theta = \theta .((z+2i\e)^{1/n} , \dots,  \underbrace{\frac{z}{z+2i\e} (z+2i\e)^{1/n}}_{k^{th}\text{ coordinate}} , \dots , (z+2i\e)^{1/n} )$$
\end{enumerate}

\end{Theorem}

\begin{proof}
(1). Let the section be $\wt{u}= ({\widetilde{u_1},\wt{u_2},\dots ,\wt{u_n}})$. Thus $\wt{u_1}(z)\dots \wt{u_n}(z) = z$ and $|u_1| = \dots |u_n|$ on $\R + i\e$. As the domain of definition is in a complex plane with the negative imaginary axis removed, $\C\setminus[0,-i\infty)$, we can define an $n\th$ root in this domain. Let $u_k= \frac{\wt{u_k}}{z^{1/n}}$. Thus $u_1.\dots.u_n=1$ on $\mathbb{H}+i\e$ and $|u_k| = 1 $ on $\R + i\e$ , so we can extend $u_k$ to $\C$ by Schwarz reflection by defining $u_k(z) = \frac{1}{\bar{u_k}(\bar{z}+2i\e)}$. Note that $E_H(\wt{u}) < \infty $ implies that $u$ is asymptotic to a trivial strip over a Reeb chord. From $|\wt{u_k} (z)| =  |z|^{1/n}$ for $z \in \R + i\e$, we conclude that $u$ is asymptotic to a trivial strip over a Reeb chord of length $\pi/n$. Thus, we have $u_k$ is bounded near $\infty$. Since $u_k$ can be extended to an entire function, we have by Liouville's theorem that $u_k$ is constant. Now, from $\wt{u_1}(z)\dots \wt{u_n}(z) = z$ and $|u_1| = \dots |u_n|$ on $\R + i\e$ we see that $\wt{u}=\theta.(z^{1/n},\dots,z^{1/n})$ for some $\theta \in T^{n-1}$.

(2). Let the section be $\wt{u}= ({\widetilde{u_1},\wt{u_2},\dots ,\wt{u_n}})$. From $\pi_m\of \wt{u} (0) = 0 $ we have that at least one of $\wt{u}_k$ is 0 at 0, without loss of generality assume $\wt{u}_1(0)=0$. Thus, in a neighborhood around 0, we have $\wt{u_1}(z) = zh(z)$ which implies $h\wt{u_2}\dots \wt{u_n}=1$ and $\wt{u_2},\dots,\wt{u_n}$ are non-vanishing on $\mathbb{H}-i\e$.  Define $u_k = \frac{\wt{u_k}}{(z+2i\e)^{1/n}}$, for $k\neq 1$ we have $u_k$ is non-vanishing and $|u_k(z)|= 1$ for $z\in \R-i\e$, thus we can do a Schwarz reflection as above to see that $u_k$ extends to an entire function. From $E_H(\wt{u}) < \infty$ and $|\wt{u_k}| = |z|^{1/n}$ on $\R -i\e$ we have that $u_k$ is bounded near infinity for $k\neq1$. Thus, from Liouville's theorem, we have $u_k$ is constant for $k \neq 1$. Now from $\pi_m \of \wt{u} (z) = z$ and the definition of $u_k$, we have that $\wt{u_1} = \frac{z}{z+2i\e} (z+2i\e)^{1/n}$. Thus, $\wt{u} = \theta.  (\frac{z}{z+2i\e} (z+2i\e)^{1/n},(z+2i\e)^{1/n},\dots,(z+2i\e)^{1/n}) $ for some $\theta \in T^{n-1}$. This gives one of $n$ $T^{n-1}$ orbits. The other orbits are obtained from the cases $u_k(0)=0$ where $k \neq 1$. Using a similar technique as above, by writing $u_k (z) = z h(z)$, one can show that the section would be of the form $s^k_\theta$.

\end{proof}
\begin{cor}\label{corr:disccount}
For a  Reeb chord $\zeta$ of angle $\pi/n$ from $T_+$ to $T_-$ we have that for $\e>0$ $|\mathcal{M}_\zeta (\R + i\e)| =1 $ and $|\mathcal{M}_\zeta (\R - i\e)| = n $
\end{cor}

\begin{Theorem}
\label{thm:regularsimple}
The sections classified in Theorem \ref{thm:classif} are transversally cut out, that is, the linearized Cauchy-Riemann operator $\wt D_u$ is surjective for all the sections $u$ classified in Theorem \ref{thm:classif}.
\end{Theorem}

\begin{proof}
This follows from Theorems \ref{thm:hortran} and \ref{thm:verttran}. 
\end{proof}

\subsection{Flattening the curve:  $\R + i\e$  to a cylindrical $\gamma$}

We will now show that the count of elementary discs with boundary on $T_\gamma$ for an admissible cylindrical $\gamma$ is the same as the disc count we obtained in Corollary \ref{corr:disccount}.  We will assume $\e>0$ hereon. We will compactify the rigid discs to discs in $\C\P^n$ and prove bijection between the compactified discs.  Recall that the closure of $T_{\R +i\e}$ is a self-clean intersecting Lagrangian. A compactified elementary disc will map the  boundary puncture mapping to the self-clean intersection at the infinity divisor $T^{n-1}_{Cliff} \hookrightarrow \C\P^{n-1}_\infty \hookrightarrow \C\P^n$.

Later we will use a Lagrangian isotopy coming from a family of asymptotically cylindrical paths. We explain the process of compactifying elementary discs with boundary on $T_{\gamma_t}$ for an asymptotically cylindrical path $\gamma_t$ . We can compactify elementary discs to discs in $\C\P^n$. Given any elementary disc $u:\D^2\setminus \{1\}$ with boundary on $T_{\R\pm i\e}$ or $T_{\gamma_t}$ we can extend $u$ to $\overline{u}$ that maps $1$ to a point in the Clifford torus $T^{n-1}_{Cliff}$ in the infinity divisor $\C\P_\infty^{n-1} \hookrightarrow \C\P^n$. The Clifford torus is the self-intersection locus of the Lagrangians $\overline T_{\R\pm i\e}$ (and $\overline T_{\gamma_t}$), thus given any holomorphic map $v: \D^2 \to \C\P^n$ with boundary on $\overline T_{\R\pm i\e}$ (or $\overline T_{\gamma_t}$) such that $v\inv (\C\P^{n-1}_\infty) = \{1\}$, we have that $v: \D^2 \setminus \{1\} \to \C^n$ is asymptotic to  a Reeb chord with boundary on $T_\pm$. Similar to Step 5 in the proof of asymptotic convergence in Theorem \ref{thm:exponential conv}, by using Theorem 3.2 of \cite{schm} the length of the Reeb chord is the same as the eigenvalue of the leading order eigenfunction. We define the moduli space of discs in $\C\P^n$ that restrict to elementary discs in $\C^n$ when restricted to $\D\setminus\{1\}$ as elementary discs in $\C^n$.

\begin{defn}[Elementary disc in $\C\P^n$]
A map $u: \mathbb{D} \to \C\P^n$  with boundary on the Lagrangian $\overline T_{\R+i\e}$ or $\overline T_\gamma$ for a cylindrical $\gamma$  is called elementary when $u\inv(\C\P_\infty^{n-1})=\{1\}$ and the eigenvalue of the leading order eigenfunction of $u$ as in Theorem 3.2 of \cite{schm} is $\pi/n$. We denote the set of elementary discs that maps $1 \mapsto p \in T^{n-1}_{Clif}$  and boundary on $\overline T_{\R+i\e}$ or $\overline T_\gamma$ by $\mathcal{M}_p(\R+i\e)$, $\mathcal M_p{(\gamma)}$ respectively.
\end{defn}

 Note that for each $p\in T_{Cliff}^{n-1}$ there are $2n$ Reeb chords $\zeta_i$ in $S^{2n-1}$ of length $\pi/n$ with boundary on $T_\pm$ that descend to $p$. From the definition of elementary discs we have the following bijections 

\begin{equation*}
    \mathcal{M}_p(\R \pm i\e) \leftrightarrow \cup_{i=1}^{2n} \mathcal{M}_{\zeta_i}(\R\pm i\e)
\end{equation*}
\begin{equation*}
    \mathcal{M}_p(\gamma) \leftrightarrow \cup_{i=1}^{2n} \mathcal{M}_{\zeta_i}(\gamma).
\end{equation*}

\noindent We denote the collection of discs in $\mathcal M_p(\R \pm i\e)$ that converges to the Reeb chord $\zeta$ in cylindrical coordinates as $$\mathcal M_p^\zeta (\R \pm i\e),$$ the notation $\mathcal{M}_p^{\zeta}(\gamma)$ is defined similarly. Clearly  we have the following bijections by extending and restricting the maps, $$\mathcal M_p^\zeta(\R \pm i\e) \leftrightarrow \mathcal{M}_{\zeta_i}(\R\pm i\e) $$ $$\mathcal M_p^\zeta(\gamma) \leftrightarrow \mathcal{M}_{\zeta_i}(\gamma) .$$

We will now prove the main theorem of the subsection, which states that there is a bijection between the moduli spaces of rigid discs if we `flatten the curve'. More precisely, if we flatten $\R + i \e$ to an admissible path $\gamma_0$, counts of elementary discs don't change.

\begin{Theorem}
\label{thm:flattenning}
For each $\e>0$ and Reeb chord $\zeta$ of length $\pi/n$ with boundary on $T_\pm$, the moduli space of elementary disc with Reeb asymptote $\zeta$ and  with boundary on $T_\gamma$ for an admissible cylindrical $\gamma$ is in bijection with the moduli space of elementary discs with Reeb asymptote $\zeta$ with bounding $T_{\R + i\e}$
$$\mathcal M_\zeta{(\gamma)} \leftrightarrow \mathcal M_\zeta{(\R + i\e)}$$
For a mutated admissible cylindrical $\gamma'$, we have a similar bijection,
$$\mathcal M_\zeta{(\gamma')} \leftrightarrow \mathcal M_\zeta{(\R - i\e)}$$

\end{Theorem}

We collect some results about regularity and area computation here before proceeding to the proof of Theorem \ref{thm:flattenning}

 \begin{remark}[Regularity for compactified elementary discs bounding $\ol T_{\gamma_t}$] 
We collect some regularity results here that allows us to obtain a cobordism. From Theorem \ref{thm:regularsimple} we have that any elementary disc in $\C^n$ is cut out regularly, i.e., the linearized Cauchy-Riemann operator at any elementary disc $u$ ,  $D_u$ is surjective. Consider the extension to the elementary disc $\overline{u}$ in $\C\P^n$, this is a holomorphic disc that maps a boundary point to the self-clean intersection locus. We can set up the moduli problem of discs with boundary on $\overline T_{\R \pm i\e}$ or $\overline T_\gamma$ similar to the treatment of moduli space of strips with boundary on cleanly intersecting Lagrangians as done by Schm\"aschke in \cite{schm}. The linearization of the $\delbar$ section in this setting is denoted as follows,$$ D_{\ol u}: W^{k,p,\delta }(\overline{u}^*T\C\P^n, \del \bar u ^* T \ol T_{\R+i\e} ) \to  W^{k-1,p,\delta } (\ol {\operatorname{Hom}} (T\mathbb{H},\ol u ^* T\C\P^n ))  $$

The following lemma shows that $\overline{u}$ is regular for any elementary disc $u$.

We begin by comparing the formal adjoints. Note that these formal adjoints $D_{\ol u}^*$  and $D^*_{u}$ are defined through different metrics, the first one is defined through a hermitian metric on $\C\P^n$ and the second one is defined through a cylindrical hermitian metric on $\C^n$. The cylindrical hermitian structure on $\C^n$ at $z$ is given by $\frac{\langle \,,\,\rangle_{0}}{|z|^2}$ where $\langle\,,\, \rangle_0$ is the standard Hermitian metric on $\C^n$  The hermitian metric on $\C\P^n$ is given by $h_{\C\P^n}(\_,\_) = g_{\C\P^n} (\_\_) - g_{\C\P^n} (\_,J\_) $
where $g_{\C\P^n} = dr^2/r^4 + \lambda_0^2/r^2 + \pi^*g_{FS}$ in the cylindrical end $\R \times S^{2n-1}$ of $\C^n$. The subscripts refer to the projection along the  splitting $T(\R \times S^{2n-1}) = T\R \oplus \R \langle R \rangle \oplus \xi,$ where $R$ is the Reeb field and $\xi$ is the contact structure on $S^{2n-1}$ to write sections as a sum of 3 projections. Consider the following rescaling of sections, 

\begin{align*}
    \label{eq:rescale}
    &S_1 : \Gamma (u^*T\C^n) \to \Gamma (u^*T\C^n)\\
    &S_1(\eta_r + \eta_R + \eta_{\xi}) = r^2(\eta_r + \eta_R) + \eta_\xi\\
    &S_2 : \Gamma(\ol{\operatorname{Hom}}(T\Sigma, u^*T\C^n)) \to \Gamma(\ol{\operatorname{Hom}}(T\Sigma, u^*T\C^n)) \\
    &S_2(\alpha_{r} + \alpha_{R} + \alpha_{\xi}) = r^2(\alpha_r + \alpha_R ) + \alpha_{\xi}.
\end{align*}

\noindent For these rescalings we have the following relation between the cylindrical and $\C\P^n$ metric 

\begin{equation}
        \langle \alpha, S_i \beta \rangle_{\C\P^n}=\langle S_i \alpha,\beta \rangle_{\C\P^n}  = \langle \alpha ,\beta \rangle_{cyl} .
\end{equation}

\noindent Using these rescalings we have, 

\begin{align}
    \label{eq:cokernelrel}
    D^*_{\ol u} = S_2 D_u^* S_2\inv.
\end{align}
     
 \end{remark}
\begin{lemma}
\label{lemm:cyltocompreg}
For any elementary disc $\ol u$ in $\C\P^n$ with boundary on  $\ol T_{\gamma_t}$ ,  the linearization $D_{\overline{u}} $  is surjective.
\end{lemma}

\begin{proof}
It will be enough to show that the kernel of the formal adjoint $D^*_{\ol u}$ is trivial.  Let $\eta$ be an element in the kernel of the formal adjoint $D^*_{\overline{u}}$, thus it is an element in $W^{k,p,-\delta} (\ol {\operatorname{Hom}} (T\mathbb{H},\ol u ^* T\C\P^n )) $. From Equation (\ref{eq:cokernelrel}) we have that $S_2\inv \eta$ will be an element in kernel of $D_u^*$. Since we have $\eta \in W^{k,p,-\delta} (\ol {\operatorname{Hom}} (T\mathbb{H},\ol u ^* T\C\P^n )) $, we have $\| \eta(s,t) \|_{\C\P^n} \in O(e^{\delta s}),$ which implies that $\| S_2\inv \eta \|_{cyl} \in O(e^{\delta s}) $ thus $S_2\inv \eta \in W^{k,p,-\delta} (\ol {\operatorname{Hom}} (T\mathbb{H}, u ^* T\C^n ))$. But we have already proved that for any elementary disc $u$, the linearization $D_u$ is surjective in Theorem \ref{thm:regularsimple}, thus $S_2\inv \eta = 0$, hence we have $\eta = 0 $. Thus, the kernel of $D^*_{\ol u}$ is trivial, hence $D_{\ol u}$ is surjective.

\end{proof}
\begin{remark}[Area computations] \label{rem:PrepSurj}
We calculate some symplectic areas that will be required for dealing with bubbling.  Recall that $\C^n$ with the Fubini-Study, $\omega_{FS}$ is symplectomorphic to the unit ball $B_1(0),\omega_0$ by the following rescaling symplectomorphism $\phi_{FS} : B_1(0) \to \C^n$, $$\phi_{FS} ((z_1,z_2,\dots,z_n)) \mapsto \frac{1}{1-\sum_i |z_i|^2} (z_1,\dots,z_n).$$ It is clear that $\phi_{FS} \inv (T_\gamma)=T_{\wt \gamma}$ for some $\wt \gamma : [0,1] \to \C^*$.  Given a elementary disc $\ol{u}_t$ in $\C\P^n$ with boundary on $\ol T_{\gamma_t}$, we can calculate it's symplectic area by using the symplectomorphism $\phi_{FS}$. Note that from convergence to Reeb chords, we can extend the map $\phi_{FS}\inv \of \ol u_n$ by adding a Reeb chord of length $\pi/n$. Recall from Lemma \ref{lemm:Texactlag} and Remark \ref{rem:texactlag} that $T_{\wt \gamma}$ is an exact Lagrangian. Moreover, for the primitive $f_{\wt \gamma}$ of $\lambda_0$,i.e. $\lambda_0|_{T_{\wt \gamma}} = df_{\wt \gamma}$, we have that $f_{\wt \gamma}$ is constant over the fiber $T^{n-1}$, in other words, there is a $f^{\wt \gamma}$ such that $f_{\wt \gamma} = \pi_m \of f^{\wt \gamma} $. Now we can use Stoke's theorem to calculate the symplectic area of $u_n$ as follows,
$ \int_{u_n} \omega_{FS}= \int_{\phi_{FS}\inv \of \ol u_n} \omega_0 = \pi/n \pm ( f^{\wt \gamma}(\wt \gamma (1)) - f^{\wt \gamma}(\wt \gamma (0)))$.

\end{remark}
\noindent Now we can proceed to the proof of Theorem \ref{thm:flattenning}.

\begin{proof}[Proof of \ref{thm:flattenning}]

The proof is based on a cobordism argument. Fix a family $\gamma_t$ of paths in $\C^*$ such that the compactification $\ol \gamma_t$ in $\C\P^1$ is a $C^2$ isotopy from $\gamma$ to $\R + i\e$. This induces a Lagrangian isotopy of $\ol T_{\gamma}$ to $\ol T_{\R + i\e}$. 

Consider the parameterized moduli space $\wt {\M} (\zeta)$ of triples $(t,\Sigma,u)$ where $u$ is a $J_0$ holomorphic map from the one-marked-disc $\Sigma$ to $\C\P^n$ with boundary on $\ol T_{\gamma_t}$ which takes the marked point to the Reeb chord $\zeta$. This is a cobordism between $\M_\zeta (\gamma)$ and $\M_\zeta(\R + i\e)$.

From Remark \ref{rem:PrepSurj} we see that there is no bubbling in this moduli space. Indeed, from exactness of $T_{\gamma_t}$ in $\C^n$, we have that any possible disc or sphere bubble has to touch the infinity divisor $\C\P^{n-1}_\infty$. One can use the area computations from Remark \ref{rem:PrepSurj} to conclude that there can't be a bubble touching the infinity divisor. From Lemma \ref{lemm:cyltocompreg} we get that the parametrized moduli space $\wt {\M}(\zeta)$ is a regularly cut-out manifold with boundary of dimension one. The Lemma \ref{lemm:cyltocompreg} also proves that the natural projection map $\wt {M}(\zeta) \to [0,1]$ is a submersion.

Thus, we get that $\wt \M (\zeta)$ is a trivial  cobordism between the zero-dimensional manifolds $\M_\zeta (\gamma)$ and $\M_\zeta(\R + i\e)$. Hence we get the required bijection $\mathcal M_\zeta{(\gamma)} \leftrightarrow \mathcal M_\zeta{(\R + i\e)}$. A similar cobordism argument works for mutated admissible paths by choosing a family $\gamma'_t$ interpolating between $\gamma'$ and $\R - i\e$.

\end{proof}

\begin{remark}[Theorem \ref{thm:flattenning} for longer Reeb chords]
    \label{rem:genflatt}
The proof technique for Theorem \ref{thm:flattenning} gives us a bijection between single punctured holomorphic discs with higher multiplicity Reeb chords as well. Let $k>0$ be a natural number. We denote the moduli space of single punctured holomorphic discs with boundary on $T_{\R+i\e}$, that is asymptotic to a Reeb chord $\zeta$ of length $k\pi/n$, as $$\M_{\zeta}(\R+i\e).$$ We similarly denote the moduli space of single punctured disc with boundary on $T_{\gamma}$ with Reeb asymptote $\zeta$  as $$\M_{\zeta}(\gamma).$$ A trivial cobordism argument as above shows that we get a bijection (more precisely, a diffeomorphism) $$\M_\zeta(\R + i\e ) \leftrightarrow \M_\zeta (\gamma).$$

\end{remark}

\begin{remark}\label{remark:counting discs}
    The results of this subsection should be interpreted as a $1\longleftrightarrow n$ correspondence between elementary discs with fixed Reeb chords under mutation. Let $\zeta_+$ be a Reeb chord of minimal length starting on $T_+$ and $\zeta_-$ be a Reeb chord of minimal length starting on $T_-.$ From Remark \ref{remark:onlyupperisenough} and Theorems \ref{thm:classif} and \ref{thm:flattenning} we get the following disc counts. For the pre-mutated Lagrangian boundary $T_{\gamma}$, there is exactly one elementary disc with Reeb chord $\zeta_+$ and exactly $n$ elementary discs with Reeb chord $\zeta_-$. For the mutated Lagrangian boundary $T_{\wt \gamma}$, the disc count gets reversed. There are exactly $n$ elementary discs with Reeb chord $\zeta_+$ and exactly one elementary disc with Reeb chord $\zeta_-$.
\end{remark}

\subsection{Classifying inner piece maps with boundary on $T_{\R+i\e}$ }
We will use a similar argument as in Theorem \ref{thm:classif} to classify all possible holomorphic punctured discs. We only need classification of single punctured disks for this article, we add this subsection for the sake of completion.

\subsubsection{Case 1: Disc with only boundary punctures}

Up to an automorphism, a punctured disc is equal to the upper half-plane with punctures on the real line. Assume that we have $k+2$ punctures on the real line,   at the points $0,1,r_1,\dots,r_k$, we denote it by $\mathbb H_{k+2}$. Let $$u:\mathbb H_{k+2} \to  \C^n$$ be a disc with boundary on $T_{\R + i\e}$. Since near the punctures , the disc $u$ is asymptotic to Reeb chords starting on the tori $T_\pm$, we have that $\pi_m \of u$ has poles at the punctures. Thus, $\pi_m \of u $ is a meromorphic map from $\mathbb H$ with boundary lying on $\R+i\e$. The usual Schwarz reflection trick gives us that $$\pi_m \of u = \frac{p}{q} + i\e,$$ where $p,q$ are real polynomials such that degree of $p$ is strictly bigger than that of $q$. Here $q(z)$ is a polynomial of the form $$z^{m_0}(z-1)^{m_1}(z-r_1)^{m_2}\dots (z-r_k)^{m_{k+2}}$$ where $m_0,\dots m_{k+2}$ are the multiplicities of the limiting Reeb chords at the punctures. Assume that $\{u_1,\dots,u_l\}$ are roots of $\pi_m \of u$ in $\mathbb H$, the roots here can have repetitions to accommodate roots with multiplicities. Let $f$ be a rational function defined as follows $$ f(z) = \pi_m \of u(z) \times \Pi_{j=1}^l \frac{z-\ol u_j}{z-u_j}.$$

\noindent From construction of $f$, we see that all the roots of $f$ lie in the lower half plane, and the poles of $f$ are on the real line. Thus, we can take a $n\th$ root of $f$ by taking $n\th$ root of each factor. This construction can be generalized for any rational function on $\mathbb H$, we define it below.

\begin{defn}[Cleaning of rational functions]
    Let $p/q$ be a rational function  on $\mathbb H$,  with zeroes at $\{u_1,\dots, u_l\}$ and poles at $\{\zeta_1,\dots \zeta_o\}$. Here the roots are poles repeat according to its multiplicitiesx. We call $f$ to be a cleaning of the rational function of $p/q$ if $$ f(z) = \frac{p(z)}{q(z)} \times \Pi_{j=1}^l \frac{z-\ol u_j}{z-u_j} \times \Pi_{j=1}^o \frac{z- \zeta_j}{z-\ol \zeta_j}. $$
\end{defn}

\begin{defn}[Weight Matrix]\label{defn:weightmatrix}
    We call a $n \times l $ matrix $B$ a \textit{weight matrix} if each column is a standard basis vector $e_i$ of $\C^n$.  $$e_i = \begin{bmatrix}
0 \\
\vdots \\ 
0 \\
1  \\
0  \\
\vdots \\
0\\ 
\end{bmatrix}$$

\end{defn}

\begin{defn}
    Let  $W$ be the set of all possible $n \times l$ weight matrices. We define an assignment function as follows,
    $$\mathfrak{A} : W \times \C^l \to \C^n$$
    $$\mathfrak{A}(B,(z_1,\dots,z_l)) = ( z_1^{b_{11}} z_2^{b_{21}}\dots z_n^{b_{n1}},   z_1^{b_{12}}z_2^{b_{22}}\dots z_n^{b_{n2}},\dots,  z_1^{b_{1l}}\dots z_n^{b_{nl}} ) $$
\end{defn}

\begin{example}
    If $B =  \begin{bmatrix}
0 & 1 \\
1  & 0 \\
0   & 0\\
\end{bmatrix}$, then $\mathfrak{A}(B,(x,y)) = (y,x,1) .$
\end{example}

\noindent A quick counting argument shows that  there are $n^l$ weight matrices of size $n\times l$. We can now state the classification theorem.

\begin{Theorem}\label{thm:classifboundarypunct}
    Let $u:\mathbb H_{k+2} \to \C^n$ be a holomorphic disc with boundary on $T_{\R + i\e}$. Assume the zeroes of $\pi_m \of u$, counted with multiplicity, are $\{x_1,\dots, x_l\}$. Let $f$ be the cleaning of the rational function $\pi_m \of u$. Then,  $$u(z) = \theta . f^{1/n} (z) \mathfrak{A} \bigg( B, \bigg( \frac{z-x_1}{z-\ol x_1},  \dots \frac{z-x_l}{z-\ol x_l} \bigg ) \bigg) $$ where $B$ is a weight matrix and $\theta$ is an element of $T^{n-1}$ acting on $\C^n$ by the phase action.
\end{Theorem}
\begin{proof}
    Let $u=( u_1, \dots,  u_n )$. We know that $u_1u_2\dots u_n$ has roots at $\{ x_1 , \dots x_l \}$, thus each root is  from some component $u_i$. If $u_i$ has roots $\{y_1,\dots y_h\}$ we define $$\wt u_i = \Pi_{j=1}^h \frac{z-\ol y_j}{z-y_j} \times u_i.$$In case $u_i$ doesn't have any roots, then we set $\wt u_i = u_i$. Denote the product $\wt u_1 \dots \wt u_n$ by $f$. Note that from the construction $f$ is actually the cleaning of the rational function $\pi_m \of u$. Thus, the $n\th$ root $f^{1/n}$ is a well-defined meromorphic function on $\mathbb H$.
    Also, since $|r-z| = |r-\ol z | $ for any real number $r$, we have  $$|\wt u_i(r)| =|\wt u_j(r)| $$ for any $i,j$ and real number $r$. Since near the punctures $u$ is asymptotic to a Reeb chord starting on the Legendrian tori $T_\pm$, the holomorphic function  $\frac{\wt u_i}{f^{1/n}}$ extends over the punctures. Since $|\wt u_i| = |f^{1/n}|$ on $\mathbb R$, from the Schwarz reflection argument in Theorem \ref{thm:classif} we have that $\frac{\wt u_i}{f^{1/n}}$ is a constant function, thus $$\wt u_i = e^{i\theta_i}f^{1/n}.$$
    Thus we have $$u_i = e^{i\theta_i}\Pi_{j=1}^h \frac{z- y_j}{z-\ol y_j} f^{1/n}.$$
    The appearance of the weight matrix in the result comes from distributing $l$ Blaschke factors $\frac{z-x_i}{z-\ol x_i}$ into $n$ coordinates. Depending on the vanishing of each coordinate functions $u_i$, we can choose a weight matrix $B$ such that 
    $$u(z) = \theta . f^{1/n} (z) \mathfrak{A} \bigg( B, \bigg( \frac{z-x_1}{z-\ol x_1},  \dots \frac{z-x_l}{z-\ol x_l} \bigg ) \bigg) .$$
   
\end{proof}

\begin{prop}\label{prop:singlepunct}
    
    Let $c$ be a Reeb chord with boundary on $T_\pm $, then there is a single punctured holomorphic disc with boundary on $T_{\R + i\e}$ with the puncture asymptotic to $c$.
\end{prop}
\begin{proof}
    Recall that we have a transitive $T^{n-1}$ action on the space of Reeb chords. Thus, it is enough to show that there is a single punctured disk asymptotic to a Reeb chord $c_\pm$ where $c_\pm$ is Reeb chord starting at $T_\pm$ of length equal to the length of $c$.

    Let the length of $c$ be an $l$ multiple of the length of the smallest Reeb chord (i.e. length of c is $\frac{l\pi}{n}$). Let $f_\pm$ be a cleaning of the polynomial $\pm z^l + i\e$ and  $\{ x_1,\dots, x_{m} \}$ be the set of zeroes of $$\pm z^l + i\e : \mathbb{H} \to \mathbb{C} .$$ Then $$u_\pm(z) = f^{1/n}_\pm(z) \bigg[ \frac{z-x_1}{z-\ol x_1} \dots \frac{z-x_m}{z-\ol x_m} \bigg ] ^T,$$ is a single punctured disk satisfying the required properties.
\end{proof}

\begin{remark}
    Note that in Proposition \ref{prop:singlepunct} the number of roots in $\mathbb{H}$ is $\lfloor{\frac{l}{2}}\rfloor$ if the Reeb chord starts at $T_+$ or else it is $\lceil{\frac{l}{2}}\rceil$.
\end{remark}

\begin{remark}
    Using Remark  \ref{rem:genflatt} we see that Proposition \ref{prop:singlepunct} holds for the boundary condition $T_{\gamma}$ for a cylindrical path $\gamma$. 
\end{remark}

\subsubsection{Case 2: Disc with interior punctures}
The classification result for discs with interior punctures ( in other words, domains with cylindrical ends) follows a similar proof strategy, except that there is more combinatorial freedom.

We illustrate the combinatorial freedom by considering a neighborhood of an interior puncture. Let $p \in \D^2$ be an interior puncture and $u$ is a holomorphic disc with boundary on the Lagrangian $T_{\R+i\e}$ which is asymptotic to a Reeb orbit near the puncture $p$. We can't read off information about the Reeb orbit by looking at $\pi_m \of u $ near $p$. This stems from the fact that the Reeb orbits can have some zero coordinates, so $\pi_m \of u$ can be any of the following: \begin{enumerate}
    \item finite at $p$.
    \item  vanish at $p$
    \item have a pole at $p$.
\end{enumerate}
Thus to account for all the cases, we need to embed the orders of vanishing and poles of the coordinates of $u$ in the hypothesis of the classification theorem. 

We define the notion of negative degree of vanishing to mean the order of the pole and negative order of pole to mean the degree of vanishing. e.g. the map $z^2$ has a pole of order $-2$ at $0$ and the map $1/z^3$ vanishes to the degree $-3$ at $0$.  We also define the degree of a function at a point to be $0$ if the evaluation of the function at the point is in $\C^*$. 

\begin{Theorem} \label{thm:interiorpunctclass}
        Fix a finite set $P$ of punctures in $D^2$. Let $u=(u_1,u_2,\dots,u_n):\D^2_{P} \to \C^n$ be a holomorphic disc from the punctured disc $\D^2_{P}$ with boundary lying on $T_{\R+i\e}$. Assume that for each puncture $p\in P$, the coordinates $u_i$ have orders of vanishing $p_i$. Let $$\wt f(z)= \pi_m \of u(z) \times \Pi_{p\in P,i=1}^n \bigg(\frac{z-\ol p}{z - p}\bigg)^{p_i}.$$ Assume $\wt f$ has zeroes, counted with multiplicities, at $\{ x_1,\dots ,x_l\}$. Let $f$ be the cleaning of $\wt f$. Then there is a weight matrix $B$ such that $$u(z) = \theta. f^{1/n} \begin{bmatrix}
\Pi_{p\in P} \big(\frac{z-p}{z-\ol p} \big)^{p_1} & 0 & 0 \dots \\
0 & \Pi_{p\in P} \big(\frac{z-p}{z-\ol p} \big)^{p_2}  & 0 \dots \\
0 & \ddots \\
\vdots  \\ 
0  \\
\vdots \\
0 & \dots & \Pi_{p\in P} \big(\frac{z-p}{z-\ol p} \big)^{p_n}  \\ 
\end{bmatrix} 
B \bigg[ \frac{z-x_1}{z-\ol x_1} \dots \frac{z-x_l}{z-\ol x_l} \bigg ] ^T
$$
\end{Theorem}

\begin{proof}
    The idea of the proof is exactly the same as Theorem \ref{thm:classifboundarypunct} after we factor out the zeroes and poles from the interior puncture. Define $\wt u_i = u_i \times  \Pi_{p\in P}^n \bigg(\frac{z-\ol p}{z - p}\bigg)^{p_i} $. Define $\wt u$ as the function $(\wt u_1, \dots , \wt u_n)$ and applying the same argument as in Theorem \ref{thm:classifboundarypunct} gives us the result.
\end{proof}

\subsection{Only elementary discs show up in rigid things}

We now prove the main result of this section : the only possible interior pieces for a rigid broken strip or a disc, is an elementary disc. The idea of the proof is to utilize monotonicity of the Lagrangian pair.

\begin{defn}[Maslov index for 2-level broken maps]
    Let $(\mathfrak{u,D})$ be a 2-level broken disc (or a broken strip) with boundary on $\mathbb L$  (and on $  \mathbb K)$. The Maslov index $\mu(\mathfrak u )$ is defined to be the Maslov index of a T-gluing $u^T$. Since the homotopy type of the glued map does not depend on the gluing parameter $T$, the Maslov index for 2-level broken maps is well-defined.  
    
\end{defn}

\begin{lemma}\label{lemm:maslovbrokenstrip}
    Any regular broken strip $\mathfrak u$ has Maslov index at least 1. A regular rigid broken strip has Maslov index 1.
\end{lemma}

\begin{proof}
    From Theorem \ref{thm:glu}, we know that for a large gluing parameter $T$, there is a regular holomorphic glued strip $u^T$ whose limit as $T\to \infty$ is $\mathfrak u$. As $u^T$ is regular, the Maslov index is at least one. Thus, Maslov index of $\mathfrak u$ is at least one.

    From Theorem \ref{thm:bij}, we know there is a large gluing parameter $T$ such that gluing gives a bijection between rigid broken strips and rigid unbroken strips. Thus, the gluing $u^T$ is rigid, hence Maslov index is one.
\end{proof}

\begin{lemma}\label{lemm:maslovbrokendisc}
    Any regular broken disc $\mathfrak u$ has Maslov index at least 2. A regular rigid broken disc with point constraint has Maslov index 2.
\end{lemma}
\begin{proof}
    
The broken disc analog can be proved by a similar gluing argument as done in the proof of Theorem \ref{lemm:maslovbrokenstrip}.
\end{proof}

\begin{Theorem}\label{lemm:onlysimple}
Each disc component of the inside piece of a rigid broken strip $\mathfrak{(u,D)}$ or a rigid broken disc with point constraint is a elementary disc.
\end{Theorem}

\begin{proof}
Since the Lagrangian pair is monotone, there is an upper limit of energy of rigid broken discs and strips. From Remark \ref{rem:PrepSurj} we see that the symplectic area of an inside disc increases as the number of  punctures increase. Similarly, we see that the symplectic area of an inside disc increases if the Reeb chord length at any of the punctures increases. Assume we are doing a neck-stretching around the boundary sphere of a Darboux ball of radius $R$. Thus, there is an upper limit $k_u>0$ such that any rigid broken strip or disc has Reeb chord asymptote of length at most $\frac{k_u\pi}{n}$.

Since the torus segment $T_\gamma$ is exact, and the punctures map to Reeb orbits or chords, from the Stokes' theorem argument in Remark \ref{rem:PrepSurj} we see that the symplectic area is a sum of contributions from the punctures.  Each boundary puncture contributes at least $c_R(\pi/n - |f^{\gamma}(0)-f^{\gamma}(1)| ) > 0$ area and each interior puncture contributes $c_R\pi > 0$ where $c_R$ depends continuously on the radius $R$.

\begin{figure}[ht]
    \centering

        \def\svgscale{1.6}
    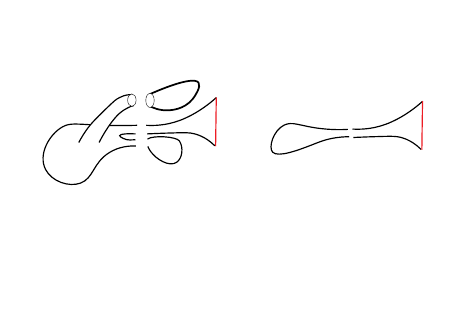

    \caption{Creating $\wt{\mathfrak u}$ by replacing the disc $A$ with $A'$ and removing the outside pieces  $B$ and $C$ }
    \label{fig:smolbuild}
\end{figure}
Assume that the inner piece of a rigid strip $(\mathfrak {u,D})$ has a disc $u$ with more than one punctures. Only one of the boundary punctures of $u$ matches with a boundary puncture of the strip piece lying outside. Denote the boundary puncture on the strip piece by $p$.  We can replace $u$ with a single boundary punctured disc $\wt u$ that
matches with the Reeb chord limit of $u_{out}$ at $p$.  Repeat this for each multiple punctured disc inside, and forget some punctured spheres or discs on the outside to obtain a broken strip $\wt {\mathfrak u}$. The total symplectic area of $\wt  {\mathfrak u} $ is less than that of $\mathfrak u$. Since the Lagrangians $L,K$ form a monotone pair, we know that there is a $c_{x,y}$ for each $x,y\in L \cap K$ such that there is an action-index relation for any strip $u$ from $x$ to $y$,
$$I_\omega (u) = \lambda I_\mu (u) + c_{x,y}.$$
Thus, since $\wt {\mathfrak u}$ has smaller area, it will have smaller Maslov index. We have that $$\mu (\wt{\mathfrak u}) < \mu (\mathfrak{ u}).$$

Since $\mathfrak u$ is a rigid disc, from Lemma \ref{lemm:maslovbrokenstrip} we have that it's Maslov index is one. Thus we have $$\mu (\wt{\mathfrak u}) < 1,$$ which contradicts that $\wt{\mathfrak u}$
is regular. As we have fixed a regularization scheme which makes every broken map of index at most one regular, we have a contradiction.

Thus, we know for any rigid strip $\mathfrak u$, every inner piece disc has exactly one puncture. From index computation in Remark \ref{rem:indexsinglepunct} we see that if the Reeb chord at boundary puncture is not of minimal length ( i.e., $\pi/n$), the dimension $\M_{in}$ is larger than the dimension of the matching condition $\mathcal{RC}$. Thus, for non-minimal Reeb chords, there is a non-trivial kernel of the evaluation map $$ev_{in} : \M_{in} \to \mathcal{RC},$$ hence such configurations cannot arise in rigid broken strips. This ends the proof for rigid broken strips. For rigid broken discs, exactly the same argument would work.
\end{proof}

%% file: pics/smolbuild.pdf_tex
\begingroup%
  \makeatletter%
  \providecommand\color[2][]{%
    \errmessage{(Inkscape) Color is used for the text in Inkscape, but the package 'color.sty' is not loaded}%
    \renewcommand\color[2][]{}%
  }%
  \providecommand\transparent[1]{%
    \errmessage{(Inkscape) Transparency is used (non-zero) for the text in Inkscape, but the package 'transparent.sty' is not loaded}%
    \renewcommand\transparent[1]{}%
  }%
  \providecommand\rotatebox[2]{#2}%
  \newcommand*\fsize{\dimexpr\f@size pt\relax}%
  \newcommand*\lineheight[1]{\fontsize{\fsize}{#1\fsize}\selectfont}%
  \ifx\svgwidth\undefined%
    \setlength{\unitlength}{225bp}%
    \ifx\svgscale\undefined%
      \relax%
    \else%
      \setlength{\unitlength}{\unitlength * \real{\svgscale}}%
    \fi%
  \else%
    \setlength{\unitlength}{\svgwidth}%
  \fi%
  \global\let\svgwidth\undefined%
  \global\let\svgscale\undefined%
  \makeatother%
  \begin{picture}(1,0.66666667)%
    \lineheight{1}%
    \setlength\tabcolsep{0pt}%
    \put(0,0){\includegraphics[width=\unitlength,page=1]{smolbuild.pdf}}%
    \put(0.28337428,0.21621628){\color[rgb]{0,0,0}\makebox(0,0)[lt]{\lineheight{1.25}\smash{\begin{tabular}[t]{l}$\mathfrak u$\end{tabular}}}}%
    \put(0.35708438,0.45167895){\color[rgb]{0,0,0}\makebox(0,0)[lt]{\lineheight{1.25}\smash{\begin{tabular}[t]{l}$C$\end{tabular}}}}%
    \put(0.33865685,0.34111382){\color[rgb]{0,0,0}\makebox(0,0)[lt]{\lineheight{1.25}\smash{\begin{tabular}[t]{l}$B$\end{tabular}}}}%
    \put(0.14127764,0.32555288){\color[rgb]{0,0,0}\makebox(0,0)[lt]{\lineheight{1.25}\smash{\begin{tabular}[t]{l}$A$\end{tabular}}}}%
    \put(0.72153971,0.22440623){\color[rgb]{0,0,0}\makebox(0,0)[lt]{\lineheight{1.25}\smash{\begin{tabular}[t]{l}$\wt{\mathfrak u}$\end{tabular}}}}%
    \put(0.61444738,0.36581329){\color[rgb]{0,0,0}\makebox(0,0)[lt]{\lineheight{1.25}\smash{\begin{tabular}[t]{l}$A'$\end{tabular}}}}%
  \end{picture}%
\endgroup%

%% file: Sections/floerhomo.tex
\section{Floer Cohomology and Mutations}\label{sec:floerhomo}
We begin this section with recalling the notion of Lagrangian Intersection Floer cohomology with local systems. After setting up the conventions and the review of the Floer cohomology, we prove the invariance of Floer cohomology under mutation. This review follows the exposition in Section 2 of \cite{wall} and the local system version follows the exposition in \cite{aurBeg}. 

\subsection{Floer Cohomology with local system} \label{subsec:flohomloc}

We explain the construction of Floer cohomology with local system for monotone Lagrangians. This was first introduced by Oh in \cite{ohMonotone}, \cite{oh:addendum}. Let $(L,K)$ be a monotone pair of compact orientable spin Lagrangians in a compact monotone symplectic manifold $M$ that transversely intersect where $L$ is locally mutable and $K$ is outside the mutation neighborhood. Assume the minimal Maslov number of $L,K$ are at least 2. Let $\mathcal{E}_L,\mathcal{E}_K$ be flat $\C$ line bundles over $L,K$ with holonomy $\rho_L,\rho_K $. We denote the tuple $(L,\rho_L)$ and $(K,\rho_K)$ by $\wt L , \wt K$. We define the chain complex  $CF(\wt L, \wt K)$ as follows.

\begin{equation*}
    CF(\wt L, \wt K) = \bigoplus_{p\in L \cap K} \hom(\mathcal{E}_{\wt {L}}|_p,\mathcal{E}_{\wt {K}}|_p )
\end{equation*}

We explain the relevant moduli spaces of strips that we use to define a coboundary operator. Fix a $t$-dependent almost complex structure $J_t$ such that $J_t=J_{std}$ in the mutation neighborhood and makes all the holomorphic strips $u:\R \times \interval \to M$ with $u(\R \times \{0\})\sub L$ and $u(\R \times \{1\})\sub K$ regular by using classical transversality results such as Proposition 3.2 in \cite{ohMonotone}.  We denote the moduli space of unparameterized (i.e., quotiented by the $\R$ translation of the domain) holomorphic strips by $\mathcal{M}(L,K)$ and we use $\mathcal{M}(L,K)_0$ to denote the moduli space of strips of index one. To ensure that a regular $J$ such that it restricts to the standard $J_{std}$ in the mutation neighborhood, we just note that there are no holomorphic strips or discs contained completely in the mutation with boundary on $L$ or $K$. We can also ensure that $J_0,J_1$ are chosen generically such that the Maslov $2$ discs with boundary on $L,K$ are regular.   The disc potential of $\wt L$ is defined as follows. Let $\beta \in H_2(M,L)$, we denote the moduli space of Maslov $2$ disc in the homology class $\beta$  by $\mathcal{M}(L,
\beta)$.  These  moduli spaces have an evaluation map $ev : \mathcal{M}(L) \to L$ after choosing a boundary point, and $n_\beta$ denotes the degree of the evaluation map. 

\begin{equation} \label{eq:discpot}
    W(\wt L) = \sum_{\beta \in \mathcal{M}(L,\beta)} n_\beta . \rho_L (\del \beta)
\end{equation}

We can similarly define the disc potential $W(\wt K)$ for the Lagrangian $K$. Given a rigid holomorphic strip $u\in \mathcal{M}(L,K)_0$ such that $\lim_{s\to -\infty}u(s,t)=p$ and $\lim_{s\to \infty }u(s,t)=q$. We denote the parallel transport map along $u(s,0)$ on the line bundle $\mathcal{E}_{\wt L}$ by  $\psi_L(u): \mathcal{E}_{\wt L}|_p \to \mathcal{E}_{\wt L}|_q $.  Similarly, define $\psi_K(u)$ to be the parallel transport along $u(s,1)$.

We now define the Floer coboundary  operator. Denote the coboundary operator as $\del : CF(\wt L, \wt K) \to CF(\wt L, \wt K) $, we will explain how  it acts on an element $\eta \in \hom (\mathcal{E}_{\wt {L}}|_p,\mathcal{E}_{\wt {K}}|_p )$ of a summand in $CF(\wt L, \wt K)$ and then extend it linearly.

\begin{equation*}
    \del \eta = \sum_{u\in \mathcal{M}(L,K)_0} \psi_K(u) \of \eta \of \psi_L(u)\inv
\end{equation*}

\noindent We can check that for a generic choice of $J$ we can arrange $\del^2 x = (W(\wt L) - W(\wt K)) x$, thus when the disc potentials are equal, we have that $\del$ is an actual differential.

\begin{Theorem}[\cite{ohMonotone},\cite{oh:addendum}]

If the disc potentials of the Lagrangians in the monotone pair match, i.e.,  $W(\wt L ) = W(\wt K)$, for a generic choice of $J_t$ we have that $\del^2=0$. We define the Lagrangian Intersection Floer homology $HF$ of $(\wt L$, $\wt K)$ as the homology of the chain complex $(CF(\wt L, \wt K),\del)$. $$HF (\wt L, \wt K) = \frac{\ker \del }{\im \del}$$

\end{Theorem}

We explain the properties on the almost complex that we require while counting strips.  We call a $t$-dependent family $J_t$ of almost complex structure \textit{regular for Floer theory} if it satisfies the following properties:
\begin{itemize}
    \item $J_0$,$J_1$ makes the moduli space of Maslov 2 discs with 1 boundary marking with boundary on $L$ or $K$ regular and the evaluation at boundary is transverse to $L \cap K$.
    \item No Chern 1 $J_t$ holomorphic sphere intersects $L\cap K$ for all $t \in \interval$.
    \item All Maslov index 1 $J_t$ holomorphic strips are regular.
    \item There are no negative index $J_t$ holomorphic strips.
\end{itemize}

\noindent Although we don't enforce regularity of index two strips, counting the index one $J_t$ holomorphic strips for a $J_t$ that is regular for Floer theory gives us a correct Floer coboundary operator $\del$. This follows from the fact that the space of regular $\J^{reg}$ of t-dependent almost complex structures is dense, and regularity of index one strips imply by a implicit function theorem argument and Gromov compactness argument that the disc count for $J_t$ is constant for small perturbations. 

\begin{remark}
Since the Lagrangian $L$ is exact in the mutation neighborhood, there are no holomorphic discs or spheres lying completely inside the mutation neighborhood. Thus, we can choose any almost complex structure in the mutation neighborhood, and then we can assume that throughout our perturbation scheme, we don't change the almost complex structure in the mutation neighborhood. We can also ensure that $J_0 = J_1$ since Maslov two discs are simple by \cite{Lazz} and the usual transversality arguments for simple holomorphic maps show that for a Baire dense set one attains regularity.
\end{remark} 

\subsection{Bulk deformed monotone Floer theory with domain dependent perturbations}

We will explain a version of monotone Lagrangian Floer cohomology where we use domain dependent almost complex structure.  In Section \ref{sec:modcons} and Section \ref{sec:reg} we obtain a domain dependent choice of almost complex structure such that the constrained moduli space of rigid holomorphic strips with markings going to a nice divisor $D$ are in bijection with rigid broken strips. We will define a deformation of the Floer coboundary $\del$ by counting strips with interior markings going to a nice divisor $D$ and show that the resulting Floer homology is the same as the usual one which we defined in the previous subsection. For more details on the construction of domain dependent Lagrangian intersection theory, see \cite{charwoodstabil}

We explain describe domain dependent almost complex structure for marked strips. We define $\M^s_k$ to be the moduli space of strips with $k$ interior markings and $\mathcal{U}^s_k \xrightarrow{\wt \pi^s_k} \M^s_k$ be the universal moduli space. A domain dependent choice of almost complex structure $\mathcal {J}_k$ for $k\geq 1$ is a map $\mathcal{J}_k : \mathcal{U}^s_k \to \mathcal{J}(M,\omega,D)$. A collection of domain dependent almost complex structure is said to be coherent if it satisfies certain properties under gluing and breaking as defined in Definition 3.14 in \cite{charwoodstabil}. A map $u:S_{\textbf{z}} \to (M,L,K)$  from a $k$ marked strip $S_{\textbf{z}} = \wt {\pi^s_k} \inv ([\textbf{z}]) $ with boundary lying on $L,K$ and the endpoints going to the intersection points in $L\cap K $, where $[\textbf{z}] \in \M^s_k$ is called  holomorphic with respect to the domain dependent almost complex structure $\mathcal{J}$ if it satisfies the domain dependent version of the Cauchy Riemann equation 
\begin{equation*}
   du(z) + \J_k(z,(u(z)))\of du \of j_S = 0 .
\end{equation*}

\subsubsection{Choosing perturbations}\label{ssec:cp}

We explain how to choose a domain dependent almost complex structure to achieve regularity of strips. We choose perturbation scheme $\mathcal{P}$ inductively on the number of interior marked points. For $k=0$ marked points, we choose a $t-$dependent almost complex structure that is regular for Floer theory and the index 2 strips are also regular. For a disc with one interior and one boundary marking, we set a domain dependent almost complex structure that makes the moduli space of discs with boundary on $L$ and $K$ regular and the evaluation map of interior marked point to $D$ with multiplicity $d$ and boundary marked point going to $L\cap K$ regular for all $d$.  For $k=1$ marked strip, the domain dependent perturbations on the boundary strata  are determined from the ones on the disc with one interior marking and the t-dependent $J_t$ chosen for $k=0$. We extend the choice to the whole of $\mathcal{U}^s_1$ so that near the end points (input and output) of the strip, $\J_1(z)= J_t(z)$ and we require any $\J_1$ holomorphic strip with the interior marking mapping to $D$ with multiplicity $d$ regular in the constrained moduli space. This can be attained by following the treatment of regularity of tangency conditions as done in \cite{CMtrans}. We proceed inductively to get a coherent perturbation data $\mathcal{P}$ that achieves the following regularity results:
\begin{itemize}
    \item  \textbf{(R1)} Any unbroken strip with the $k$ markings meeting $D$ with differing intersection multiplicities is regular in the constrained moduli space where we enforce the intersection numbers at the markings.
    \item \textbf{(R2)}Every rigid strip has no tangency to $D$ at the interior markings and is regular.
\end{itemize}

Using the transversality theorem 4.1 in \cite{charwoodstabil} we can show that there is a comeagre set of coherent perturbation data which achieves the regularity results \textbf{(R1-2)} as given above. Let $\M^{\mathcal{P}}(M,L,K,D)_0$ be the moduli space of rigid marked holomorphic strips with respect to a coherent perturbation data $\mathcal{P}$ that achieves the regularity result \textbf{(R1-2)}. Since $L,K$ is a  monotone pair of Lagrangians, there is an upper bound on the energy of index 1 strips, thus by a Gromov compactness argument we see $\M^{\mathcal{P}}(M,L,K,D)_0$ is compact. Thus, regularity enforces that there are only finitely many rigid strips with interiors marking going to $D$. From \textbf{(R2)} we also get that any rigid strip intersects $D$ transversally.

\subsubsection{The $\del^{\mathcal{P}}_D$ coboundary operator}

We define a coboundary operator by counting rigid strips with marked points mapping to a nice divisor. Assume line bundles with flat connection have been chosen on $L,K$.  We define the chain complex  $CF^{\mathcal{P}}_D(\wt L, \wt K)$ 

\begin{equation*}
    CF^{\mathcal{P}}_D(\wt L, \wt K) = \bigoplus_{p\in L \cap K} \hom(\mathcal{E}_{\wt {L}}|_p,\mathcal{E}_{\wt {K}}|_p ).
\end{equation*}

\noindent It is the same vector space as $CF(\wt L,\wt K)$, we will now define the deformed Floer coboundary as follows
\begin{equation*}
    \del^{\mathcal{P}}_D \eta = \sum_{u\in \M^{\mathcal{P}}(M,L,K,D)_0}  \psi_K(u) \of \eta \of \psi_L(u)\inv.
\end{equation*}

The above expression for $\del^{\mathcal{P}}_D$ is a finite sum because we can apply Gromov compactness to the set of index 1 strips since $(L,K)$ is a monotone pair. We define deformed disc potential map $W_D(\wt L)$ and $W_D(\wt K)$ similar to Equation (\ref{eq:discpot}) by counting Maslov 2 discs with interior marking going to $D$. Since our Lagrangians are monotone, we can show that $W_D = W$. We can do so by choosing a generic path connecting the almost complex structures that gives us a cobordism between the space of rigid disc,  since there is no sphere or disc bubbling owing to monotonicity of the Lagrangians.

\begin{Theorem}

If the disc potentials of the Lagrangians in the monotone pair match, i.e.  $W(\wt L ) = W(\wt K)$, for a generic choice of perturbation scheme domain dependent almost complex structure $\mathcal{P}$, we have that ${\del^{\mathcal{P}}_D}^2=0$. We define the bulk $D$ deformed Lagrangian Intersection Floer homology $HF^{\mathcal{P}}_D$ of $(\wt L$, $\wt K)$ as the homology of the chain complex $(CF(\wt L, \wt K),\del)$. $$HF^{\mathcal{P}}_D(\wt L , \wt K) = \frac{\ker \del^{\mathcal{P}}_D }{\im \del^{\mathcal{P}}_D}$$

\end{Theorem}

\begin{proof}

We study the compactified one dimensional moduli space of strips to show when the deformed disc potentials match, we can define a Floer cohomology.  Consider the one dimensional moduli space $\M^{\mathcal{P}}(M,L,K,D)_1$, the boundary of this moduli space is obtained by strip breaking or Maslov two disc bubbling as in the case of $J_t$ strips. Although we have possible disc or sphere bubbling by interior markings going to boundary, or coming together, these situations can't occur. Indeed, since the nice divisor $D$ is disjoint from $L \cup K$ any disc bubbling caused by interior marking going to boundary has to be non-constant on the disc, thus at least Maslov two from monotonicity assumptions. Since the moduli space of strips with interior markings tangent to $D$ is regular, constant sphere bubbles formed due to collision of interior marking can't occur because of index reasons. Thus, we get the boundary of   $\M^{\mathcal{P}}(M,L,K,D)_1$   consists of Maslov 2 disc bubbles and strip breaking.  The Maslov two disc bubbles in the boundary gives us the usual relation, 
$${\del^{\mathcal{P}}_D}^2 = W_D(\wt L) - W_D(\wt K)=W(\wt L) - W(\wt K).$$ When the disc potentials are equal, we define the bulk $D$ deformed Floer homology as  $HF^{\mathcal{P}}_D(\wt L , \wt K) = \frac{\ker \del^{\mathcal{P}}_D }{\im \del^{\mathcal{P}}_D}$ .

\end{proof}

\subsubsection{Chain homotopy and isomorphism of Floer cohomologies}
We will prove that this deformed Floer cohomology is isomorphic to the usual Lagrangian Floer cohomology. We will show that if we choose a $t$ dependent almost complex structure $J_t$ that makes evaluation from an interior marked strip transverse to $D$ and it's higher jet bundles (equivalent to saying tangency constrained moduli spaces are regular), the usual $\del$ Floer coboundary is chain homotopic to $\del^{\mathcal{P}}_D$. We can choose such a $J_t$ because there are countably many constraints and every non-constant strip is $s-$somewhere injective as shown in \cite{FlHStrans} and then countable intersection of comeagre sets are comeagre we get that all the constraints can be cut out regularly for a comeagre set of $t$-dependent almost complex structure. This choice of $J_t$ can be viewed as a coherent perturbation data $\mathcal{P}_t$ that is translation invariant under the $\R$ action on the strip, thus $\del_{D}^{\mathcal{P}_t} = \del$ since we can just forget the markings that go to $D$ to get a bijection between the moduli space of rigid $\mathcal{P}_t$ strips and $J_t$ strips.

\begin{lemma}\label{lemm:bulkissame}
 The bulk deformed monotone Lagrangian Floer cohomology is isomorphic to monotone Lagrangian Floer cohomology as defined in \cite{ohMonotone} i.e. $HF^{\mathcal{P}}_D(\wt L , \wt K) \cong HF(\wt L,\wt K)$
\end{lemma}

\begin{proof}
The idea of showing chain homotopy is the same as in the proof of independence of Floer homology from the choice of almost complex structure $J_t$. A parameterized strip with $k$ interior markings and one boundary marking which is a ``fake" marking added to parameterize the strip, intuitively this boundary marking on the strip fixes the $0$ in the $s$ coordinate of $(s,t) \in \R \times \interval$. We define a homotopy $\J_h$ between $J_t$ and $\J_k$ such that for $s<-1$, $\J_h = J_t$ and  for $s>1$, $\J_h = \J_k$. We create a map $\phi_h : CF(\wt L,\wt K) \to CF^{\mathcal{P}}_D(\wt L,\wt K)$  by counting rigid (Maslov index 0) parameterized $\J_h$ holomorphic strips between intersection points in $L \cap K$. We denote the moduli space of  rigid $\J_h$ parameterized strips by $\M_{par}^{\J_h}(M,L,K,D)_{0}$. By the usual generic transversality results, we know that a generic homotopy $\J_h$ would be enough to achieve regularity.  

\begin{equation*}
    \phi_h \eta = \sum_{u\in \M_{par}^{\J_h}(M,L,K,D)_{0}}  \psi_K(u) \of \eta \of \psi_L(u)\inv
\end{equation*}

We can check that $\phi_h$ is a chain map i.e. $\del^{\mathcal{P}}_D \of \phi_h = \phi_h \of \del$ by considering the breaking of $1$ dimensional( Maslov index 1) moduli space $\J_h$ strips. The only breaking that occurs is when there is a strip breaking. If the strip breaking occurs near $-\infty$ of, the parameterized strip  contributes to the count in $\phi_h \of \del$ and when strip breaking occurs near $+\infty$ it contributes to the count in $\del^{\mathcal{P}}_D \of \phi_h$.

We define $\J_{h\inv}$ on a parameterized strip by $\J_{h\inv}(s,t) = \J_h(-s,t)$. We can count index 0 solutions to get a map $\phi_{h\inv}$. We can show $\phi_h \of \phi_{h\inv}$ is chain homotopic to identity by considering a twice-parameterized moduli space of holomorphic strips. Twice parameterization refers to choosing two boundary markings $w_1,w_2$ on a strip such that $w_1 \leq w_2$ in the $s-$coordinate. In the universal moduli space of twice parameterized  strips with $k$ markings, we define a complex structure $\J_{ch}$ such when $w_1 = w_2$, $\J_{ch} = J_t$ and at the lower dimensional strata of strip breaking with $w_1$ and $w_2$ lying in different strips, $\J_{ch}= \J_h$ in the strip component with $w_1$ and  $\J_{ch}= \J_{h\inv} $ in the component with $w_2$. We now count rigid $\J_h$ holomorphic twice parameterized strips  (Maslov index -1) to create a chain homotopy map $H:CF(\wt L,\wt K) \to CF(\wt L,\wt K) $. By considering the breaking of index 0 strips, we get that $H\of \del + \del \of H=\phi_h \of \phi_{h\inv} - Id$. Hence, we get that homologies are the same, $HF^{\mathcal{P}}_D(\wt L , \wt K) \cong HF(\wt L,\wt K)$.

\end{proof}

\subsection{Proof of Mutation formula}

We have now all the materials needed to prove the invariance of $HF$ under mutation up to a change in the local system. 

We recall  the notations and conventions we introduced in Section \ref{sec:intro}.
\begin{listC}
    \item $M$ is a compact, monotone, symplectic manifold of dimension $2n$.
    \item $L$ is a \textit{locally mutable},  oriented,  spin, compact, connected, monotone Lagrangian manifold in $M$. 
    \item $K$ is an oriented,  spin, compact, connected,  monotone Lagrangian manifold in $M$. 
    \item The Lagrangian $K$ intersects $L$ transversely but does not intersect a small \textit{mutation neighborhood} $B$ of $L$. 
    \item  $(L, K)$ forms a \textit{monotone tuple}. 
    \item $\rho_L$ and $\rho_K$ are local systems on $L$ and $K$ respectively.
    \item In a mutation neighborhood $B$,  under the identification given by the symplectomorphism $\phi$ ,  $L$ is equal to the \textit{torus segment} $T_\gamma$. 
    \item $x_1, \dots x_{n-1}$ is the image of generators of $H_1(T_\gamma)$ under the inclusion map $i:T_\gamma \hookrightarrow L$.  
    \item $(x_1, x_2, \dots, x_n,  y_1, \dots y_l)$ is a generating set of $H_1(L)$ such that $x_n \cap i_*[T_\gamma]$ where $[T_\gamma]$ is the generator of $H_{n-1}(T_\gamma)$. A Poincar\'{e} duality argument ensures the existence of such an $x_n$.
\end{listC}

Given a flat line bundle over $L$, we denote the holonomies $\rho(x_i)=z_i$, $\rho(y_i)=w_i$. The  \textit{mutated local system} on the mutation $L_\mu$ is equal to $\rho^\mu(x_i)=z_i$ for $i<n$ and $\rho^\mu(x_n)=z_n(1+z_1 + \dots + z_{n-1})$  and $\rho^\mu(y_i)=w_i$ for all $i$. We denote the mutated Lagrangian with this local system as $\wt{ L_{\mu}}$.

\begin{Theorem}[Mutation]\label{premain}
The disc potential is invariant under mutation if the local system $\rho$  is changed to $\rho^\mu$ ie. $W(\wt L) = W(\wt L_\mu)$. Moreover, if $W(\wt L) = W(\wt K)$ then the Floer homology is invariant up to the same change of local system i.e. $HF(\wt L,\wt K) \cong HF(\wt{ L_{\mu}},\wt K)$.
\end{Theorem}

\begin{proof}

We can use a Hamiltonian perturbation on $L$ so that it is identified to the torus segment  $T_{\gamma'}$ in the mutation neighborhood such that near 0, $\gamma' = c\gamma$ for a very small $c$. For any $\e>0$ we can choose $c$ small enough that $T_{\gamma'}$ is cylindrical near a contact sphere $S^{2n-1}_\e$ of radius $\e$. Thus, we can do a neck stretching along the hypersurface $S^{2n-1}_\e$. We choose a perturbation scheme $\mathcal{P}$ as given in Theorem \ref{thm:ultimatetransversality}, thus  from the Lemma \ref{lemm:2level} and Theorem \ref{lemm:onlysimple} we get that the only rigid broken strips or rigid broken disc with point constraint have simple discs. We know gluing results in Corollary \ref{corr:rigidregularbijection} give us a bijection between rigid broken strips and rigid strips of a neck-stretched manifold $M^T$ with boundary on $L^T$ or ${L_\mu}^T$ and $K$. We similarly get a bijection between the rigid broken disc with point constraints and the rigid unbroken disc with point constraints.

Note that the broken discs with boundaries lying entirely on the outside under mutation do not change. Thus, we only need to study the broken discs which have inside pieces with boundaries lying on the Lagrangians.

The classification results in Theorems \ref{thm:classif} and \ref{thm:flattenning} give us that for a fixed Reeb chord $\xi$ of length $\pi/n$, there is  one upper disc(simple discs that are sections over the positive half-space) and $n$ lower discs(simple disc over negative half space) with boundary on $L_{in}$. On the other hand, for ${L_\mu}_{in}$ there are $n$ upper discs and $1$ lower disc going to the Reeb chord $\xi$. See Remark \ref{remark:counting discs}. Thus, there is  a $n\longleftrightarrow1$ and $1\longleftrightarrow n$ correspondence between rigid broken strips with boundary on $\mathbb{L}$ and $\mathbb{L}^\mu$. There is a similar correspondence between broken discs.

Applying the gluing result in Corollary \ref{corr:rigidregularbijection}, we get that there is a $1\leftrightarrow n$ and $n\leftrightarrow1$ correspondence between the rigid strips with respect to the glued perturbation scheme $\mathcal{P}^{gl}$ in the stretched manifold $M^T$. Since all the  $\mathcal{P}^{gl}$ holomorphic rigid strips are regular, and by monotonicity, there is an upper bound on energy, we get that for a generic perturbation data $\wt {\mathcal{P}}$ close to $\mathcal{P}^{gl}$, the count of rigid strips don't change. We choose a domain-dependent perturbation as in subsection \ref{ssec:cp}. A similar argument works for broken discs with point constraints. Thus, by comparing the homology classes under the $1\leftrightarrow n$ and $n\leftrightarrow1$ correspondence we get the invariance of disc potential,

\begin{equation*}
    W(\wt L) = W(\wt L_\mu).
\end{equation*}

Similarly, using the $1\leftrightarrow n$ and $n\leftrightarrow1$ correspondence, and keeping track of the homology classes, we see that the coboundary operators match up exactly,

$$(\del^{\wt {\mathcal{P}}}_D)_{(\wt L,\wt K)}  = (\del^{\wt {\mathcal{P}}}_D)_{(\wt {L_\mu} ,\wt K)}  $$

\noindent where $(\del^{\wt {\mathcal{P}}}_D)_{(\wt L,\wt K)}$ is the bulk $D$ deformed Floer coboundary for the pair $(\wt L,\wt K)$  and $(\del^{\wt {\mathcal{P}}}_D)_{(\wt L,\wt K)}$  is the bulk $D$ deformed Floer coboundary for $(\wt {L_\mu} ,\wt K)$. Thus, we have $ HF^{\mathcal{P}}_D(\wt L , \wt K) \cong HF^{\mathcal{P}}_D(\wt {L_\mu} , \wt K) $. Now using Lemma \ref{lemm:bulkissame} we get that the Floer homologies are the same.

\end{proof}

Now, Theorem \ref{maintheorem} directly follows from Theorem \ref{premain}.

\subsection{Applications}

\begin{figure}
    \centering

        \def\svgscale{0.8}
    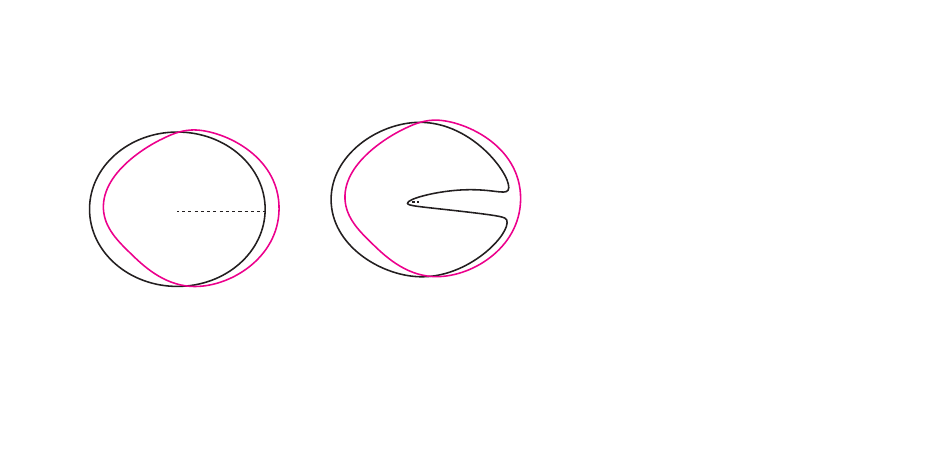

    \caption{Clifford and Chekanov type tori}
    \label{clifchek}
\end{figure}
\begin{proof}[Proof of Corollary \ref{cor:chekclif}]
    We will apply Theorem \ref{maintheorem} to $L,\phi(L)$ where $\phi(L)$ is a Hamiltonian perturbation of $L$ such that  $\phi(L) \cap L$ is away from the mutation neighborhood. Such a perturbation always exists for Clifford type torus in \cite{wall}, see Figure \ref{clifchek}. By applying Theorem \ref{maintheorem} to $L, K = \phi(L)$. We have 
    $$HF(\wt L , \wt \phi(L)) \cong HF (\wt L_\mu , \wt \phi(L)).$$
    Now using invariance of $HF$ under Hamiltonian isotopy we have $$HF(\wt L, \wt L) \cong HF(\wt L_\mu, \wt L).$$
    Thus we have that the Clifford and Chekanov tori are Hamiltonian non-dis placeable.
\end{proof}

\begin{remark}[Non-Displaceability of Chekanov tori and Real Projective spaces in odd dimensions]\footnote{This answers a question raised in Remark 3.2.3  of \cite{kawashelu}.}
    
One straightforward generalization of Theorem \ref{maintheorem} can be obtained by considering Lagrangians with local-systems coming from $\mathbb{F}-$line bundles, where $\mathbb{F}$ is any field. The proof of Theorem \ref{maintheorem} is based on comparing moduli space of rigid broken disks. The field $\mathbb{C}$ was not essential in any part of the proof, thus we can replace the complex local system with $\mathbb{F}$-local systems defined by the holonomy map $\rho_L: H_1(L) \to \mathbb{F}^\times $.

An application of such a generalization would let us show that for odd $n$, the higher mutation $T_{Ch}$ of the Clifford torus $T_{Cl}$ in $\mathbb{P}^n$ is non-displaceable from $\mathbb{RP}^n$. This follows from the computation of $HF(T_{Cl},\mathbb{RP}^n)$ in \cite{alston} and  from using $\Z_2$ coefficients in Theorem \ref{maintheorem} and seeing that the mutation formula 
 fixes the trivial local system, i.e. $(1,1,\dots, 1) = (1,1,\dots, n) \mod 2$ when $n$ is odd. 
\end{remark}

%% file: pics/clifchek.pdf_tex
\begingroup%
  \makeatletter%
  \providecommand\color[2][]{%
    \errmessage{(Inkscape) Color is used for the text in Inkscape, but the package 'color.sty' is not loaded}%
    \renewcommand\color[2][]{}%
  }%
  \providecommand\transparent[1]{%
    \errmessage{(Inkscape) Transparency is used (non-zero) for the text in Inkscape, but the package 'transparent.sty' is not loaded}%
    \renewcommand\transparent[1]{}%
  }%
  \providecommand\rotatebox[2]{#2}%
  \newcommand*\fsize{\dimexpr\f@size pt\relax}%
  \newcommand*\lineheight[1]{\fontsize{\fsize}{#1\fsize}\selectfont}%
  \ifx\svgwidth\undefined%
    \setlength{\unitlength}{450bp}%
    \ifx\svgscale\undefined%
      \relax%
    \else%
      \setlength{\unitlength}{\unitlength * \real{\svgscale}}%
    \fi%
  \else%
    \setlength{\unitlength}{\svgwidth}%
  \fi%
  \global\let\svgwidth\undefined%
  \global\let\svgscale\undefined%
  \makeatother%
  \begin{picture}(1,0.5)%
    \lineheight{1}%
    \setlength\tabcolsep{0pt}%
    \put(0,0){\includegraphics[width=\unitlength,page=1]{clifchek.pdf}}%
    \put(0.11501954,0.18565046){\makebox(0,0)[lt]{\lineheight{1.25}\smash{\begin{tabular}[t]{l}$L$\end{tabular}}}}%
    \put(0.27024007,0.20072586){\makebox(0,0)[lt]{\lineheight{1.25}\smash{\begin{tabular}[t]{l}$\textcolor{red}{\phi(L)}$\end{tabular}}}}%
    \put(0.37576775,0.19011726){\makebox(0,0)[lt]{\lineheight{1.25}\smash{\begin{tabular}[t]{l}$L_\mu$\end{tabular}}}}%
    \put(0.54829702,0.220268){\makebox(0,0)[lt]{\lineheight{1.25}\smash{\begin{tabular}[t]{l}$\textcolor{red}{\phi(L)}$\end{tabular}}}}%
  \end{picture}%
\endgroup%